\documentclass{article}
\usepackage[utf8]{inputenc}
\usepackage{appendix}
\newcommand{\bigzero}{\mbox{\normalfont\Large\bfseries 0}}

\newcommand{\R}{\mathbb{R}}

\newcommand{\cL}{\mathcal{L}}

\newcommand{\cJ}{\mathcal{J}}

\usepackage[a4paper,top=1.45cm,hmargin=1.45cm,bottom=1.45cm,includeheadfoot]{geometry}

\usepackage[most]{tcolorbox}
\definecolor{ballblue}{rgb}{0.74, 0.83, 0.9}

\definecolor{ballblue}{rgb}{0.13, 0.67, 0.8}
\definecolor{darkcerulean}{rgb}{0.03, 0.27, 0.49}
\usepackage{tikz}
\usetikzlibrary{shapes,arrows,arrows.meta}
\usepackage{makecell}
\tikzstyle{bag} = [align=center]
\tikzstyle{box} = [rectangle, draw, text centered, minimum height=3em, text width=5cm]
\tikzstyle{diamond} = [diamond, draw, text centered, minimum height=3em]
\tikzstyle{para}=[trapezium, draw, text centered, trapezium left angle=60, trapezium right angle=120, minimum height=1.5cm]
\tikzstyle{freccia} = [draw, -latex']

\newcommand{\pt}[1]{\partial_{t} #1}

\newcommand{\bs}{\mathbf{s}}
\newcommand{\z}{\mathbf{z}}

\newcommand{\bW}{\mathbf{W}}
\newcommand{\bb}{\mathbf{b}}

\newcommand{\bO}{\mathbf{0}}

\newcommand{\bz}{\mathbf{z}}
\newcommand{\bv}{\mathbf{v}}

\newcommand{\bx}{\mathbf{x}}

\newcommand{\bu}{\mathbf{u}}

\newcommand{\argmin}{\textrm{argmin}}
\usepackage[numbers]{natbib}        
\usepackage{amsmath}
\usepackage{tabularx} 
\usepackage{graphicx} 
\usepackage{float}
\usepackage{multirow}
\usepackage{amssymb}
\usepackage{colortbl}
\usepackage{mathtools}
\usepackage{amsfonts}
\usepackage[ruled,vlined]{algorithm2e}

\SetKwComment{Comment}{}{} 
\usepackage{colortbl}
\usepackage{mathrsfs,bbold}
\usepackage{booktabs}
\usepackage{hhline}
\usepackage{subcaption}
\usepackage{adjustbox}
\usepackage{multirow}
\usepackage{a4wide,color}
\usepackage{enumitem}
\usepackage[normalem]{ulem}
\usepackage{hyperref}
\hypersetup{colorlinks,allcolors=black}
\usepackage{comment}
\usepackage{caption}
\usepackage{subcaption}

\newtheorem{thm}{Theorem}[section] 
\newtheorem{rmk}{Remark}[section] 
\usepackage{amssymb}

\begin{document}\vspace{-1.5cm}
\begin{center}
{\huge Control of high-dimensional collective dynamics by deep neural feedback laws and kinetic modelling\\\vspace{0.5cm}}{\large
Giacomo Albi$^a$, Sara Bicego$^b$, Dante Kalise$^b$\\\vspace{0.5cm}}\footnotesize
    $^a$ \emph{Department of Computer Science, University of Verona, Strada le Grazie 15, Verona, 37134, Italy}\\
    $^b$ \emph{Department of Mathematics, Imperial College London, South Kensington \\
Campus, London, SW72AZ, UK}
\end{center}

\vspace{0.1cm}
\begin{abstract}
Modeling and control of agent-based models is twice cursed by the dimensionality of the problem, as both the number of agents and their state space dimension can be large. Even though the computational barrier posed by a large ensemble of agents can be overcome through a mean field formulation of the control problem, the feasibility of its solution is generally guaranteed only for agents operating in low-dimensional spaces. To circumvent the difficulty posed by the high dimensionality of the state space a kinetic model is proposed, requiring the sampling of high-dimensional, two-agent sub-problems, to evolve the agents' density using a Boltzmann type equation.  Such density evolution requires a high-frequency sampling of two-agent optimal control problems, which is efficiently approximated by means of deep neural networks and supervised learning, enabling the fast simulation of high-dimensional, large-scale ensembles of controlled particles. Numerical experiments demonstrate the effectiveness of the proposed approach in the control of consensus and attraction-repulsion dynamics.
\end{abstract}

\tableofcontents

\section{Introduction}

Collective behaviour in agent-based models (ABMs) is of evergrowing interest across various disciplines, including mathematics, physics, biology, and economics. ABMs enable the description of complex phenomena through a general paradigm that combines endogenous interactions between agents with external influences. Their applicability spans diverse areas, such as social sciences\cite{MAS_animalscollective,MAS_social}, robotics, and computer science  \cite{MAS_robs}. A fundamental topic of interest in ABMs is the study of pattern formation and self-organization \cite{survey_MAS_collective_behaviour,survey_MAS4}. However, beyond self-organization, a fascinating topic arises in relation to the design of external signals or controls to influence a system and inducing a prescribed collective behaviour \cite{survey_MAS3}.

Agent-based models encode pairwise agent-to-agent interactions through a balance of attraction and repulsion forces acting over first or second-order dynamics, while the influence of the external world on the system is expressed as a suitable control signal. To make matters more concrete, let us consider a second-order system with $N$ agents in $\R^d$, where the state of the $i$-th agent is encoded by the pair $s_i = (x_i,v_i)\in\Omega_x\times\Omega_v\subset\R^{2d}$, representing position and velocity, respectively, evolving according to transport-interaction dynamics of the form
\begin{equation}\label{dyn}
\begin{aligned}
    \Dot{x}_i(t) &= v_i(t)\\
    \Dot{v}_i(t) &= \dfrac{1}{N} \sum\limits_{j=1}^{N} P(x_i(t),x_j(t))(v_j(t)-v_i(t)) + u_i(t)\,,\qquad i = 1,...,N.
\end{aligned}
\end{equation}
Here, $P(\cdot,\cdot):\R^d\times\R^d\to\R$ denotes an interaction kernel, while $u_i(t)$ is a control signal influencing agent $i$. The ensemble of control signals is denoted by $\bu(t)=(u_1(t),\ldots,u_N(t))^{\top}$.  The core of the self-organization behaviour of the free dynamics resides in $P$, which can induce clustering, polarization or alignment, among many others. This self-organization behaviour can be modified by the influence of an external control law. In the framework of dynamic optimization, this control is synthesized by minimizing a cost functional which rewards the convergence of the system towards a cooperative goal, e.g. consensus
\begin{equation}\label{Cost}\tag{OCP}
    \tilde \bu(t) = \underset{\bu(\cdot)}{\textrm{argmin}} \bigg\{\mathcal{J}_N(\bx,\bv,\bu) := \dfrac{1}{N} \int\limits_{0}^{+\infty} \sum\limits_{i=1}^{N} \|v_i-\Bar{v}\|^2 + \gamma\, \ell(u_i)\, dt\bigg\}\,, \qquad \gamma>0\,,
\end{equation}
where the first term in the cost is promoting consensus towards the target velocity $\Bar{v}$, while $\ell:\R^d\to\R_{+}\cup\{0\}$ is a convex function penalizing the energy spent by the control $u_i$.

The solution of the optimal control problem defined by minimizing \eqref{Cost} subject to \eqref{dyn} is \textsl{twice cursed} by the dimensionality of the problem: 
solving the OCP becomes prohibitively expensive for large values of $N$ (as in swarm robotics or collective animal behaviour where hundreds or thousands of agents are present), as well as for high-dimensional state spaces (i.e. $d\gg1$, as in portfolio optimization). {  A partial remedy to this problem comes from statistical mechanics by assuming the number of agents $N\to\infty$, so that the dynamics of the individual-based problem can be approximated by a mean field equation 
\begin{equation}\label{mfpde}
\pt f +v\cdot \nabla_x f = - \nabla_v \cdot \bigg((\mathcal{P}[f] + u)f\bigg)\,,
\end{equation}
where $f(t,x,v)$ is the probability density of having an agent with state $(x,v)$ at time $t$, and the mean field interaction force is given by the non-local operator
\begin{equation}\label{operator}
    \mathcal{P}[f](t,x,v) := \int\limits_{\Omega_x\times\Omega_v} P(x,x_*)(v_*-v) f(t,x_*,v_*) dx_*\,dv_*\,.
\end{equation}
A direct transcription of the objective \eqref{Cost} leads to the mean field optimal control problem \cite{mf_ocp}
\begin{equation}\tag{MFOC}\label{mf}
\underset{u}{\textrm{min}}\; \mathcal{J}(f,u) := \int\limits_{0}^{+\infty}\!\!\int\limits_{\Omega_x\times\Omega_v}\!\! \left(\|v-\Bar{v}\|^2 + \gamma\, \ell(u)\right)\,f(t,x,v)\,dx\,dv\,dt\,.
\end{equation}
Even though the formulation using \eqref{mfpde} alleviates the curse of dimensionality with respect to $N$, it leads to a PDE-constrained optimization problem over $2d+1$ dimensions, which becomes prohibitively expensive already for moderate values of $d$. The solution of high-dimensional mean field optimal control problems has been addressed using deep learning techniques in \cite{PINNs_new} in the case of linear dynamics, however, the nonlinear case remains open. In this paper, we propose the synthesis of a feedback control for \eqref{mfpde} which is a suboptimal solution to \eqref{mf}. In order to overcome the obstacle posed by the treatment of \color{black}high-dimensional agents with nonlinear interactions\color{black},} we resort to modeling the evolution of the agents density $f$ from a kinetic viewpoint, reformulating the mean field controlled dynamics as a Povzner-Boltzmann type equation
\begin{equation}\label{pv}
	\partial_t f(t,x,v)  + v\cdot \nabla_x f = \lambda  \mathcal{Q}_{\eta,u}(f,f) (t,x,v)\,,
\end{equation}
where the operator $\mathcal{Q}$ takes into account the gain and the loss of particles in $(x,v)$ at time $t$, due to the motion of individuals via free transport and the velocity changes resulting from the controlled interaction dynamics, see  \cite{Povzner62,Kinetic_MC}. If we consider the interaction dynamics in \eqref{dyn}, reduced to the case of $N=2$ particles with velocities $(v,v_*)$, by denoting as $(v',v_*')$ the update of those velocities after a forward Euler step of length $\eta$, we have
\begin{equation}
		\mathcal{Q}_{\eta,u}(f,f)(t;x,v) = \int\limits_{\Omega_x\times\Omega_v}\dfrac{1}{\mathcal{J}_\eta} f(t,x,'v)f(t,x_*,'v_*) - f(t,x,v)f(t,x_*,v_*)\,dx_*\,dv_*\,,
\end{equation}
where $('v,'v_*)\longmapsto(v,v_*)$ are the pre-interaction velocities that generate the couple $(v,v_*)$, while $\mathcal{J}_\eta$ is the Jacobian of the binary interactions map $(v,v_*)\longmapsto(v',v_*')$. 

{  We will show that the alternative mesoscopic description of ABMs via eq. \eqref{pv} is consistent with a sub-optimal, controlled solution of the constrained mean field PDE \eqref{mfpde}}, when assuming high frequency $\lambda $ and weak interactions $\eta$ between particles, similarly to the {\em grazing-collision} limit in kinetic theory \cite{quasielasticKin}. The convenience of such kinetic formulation relies on the reduced computational cost required for its solution via \emph{direct simulation Monte Carlo} (DSMC) methods, \cite{kinetic_opinion}. {  Such a sampling-based technique approximates the solution of the Boltzmann eq.\eqref{pv}} by computing a collection of binary sub-problems for couples of agents sampled from the population density function $f(t,\cdot)$ \cite{Kinetic_MC}. 
\begin{figure}[t]
    \centering
    \includegraphics[trim={3cm 10cm 7.5cm 0cm}, clip,width =0.9\textwidth]{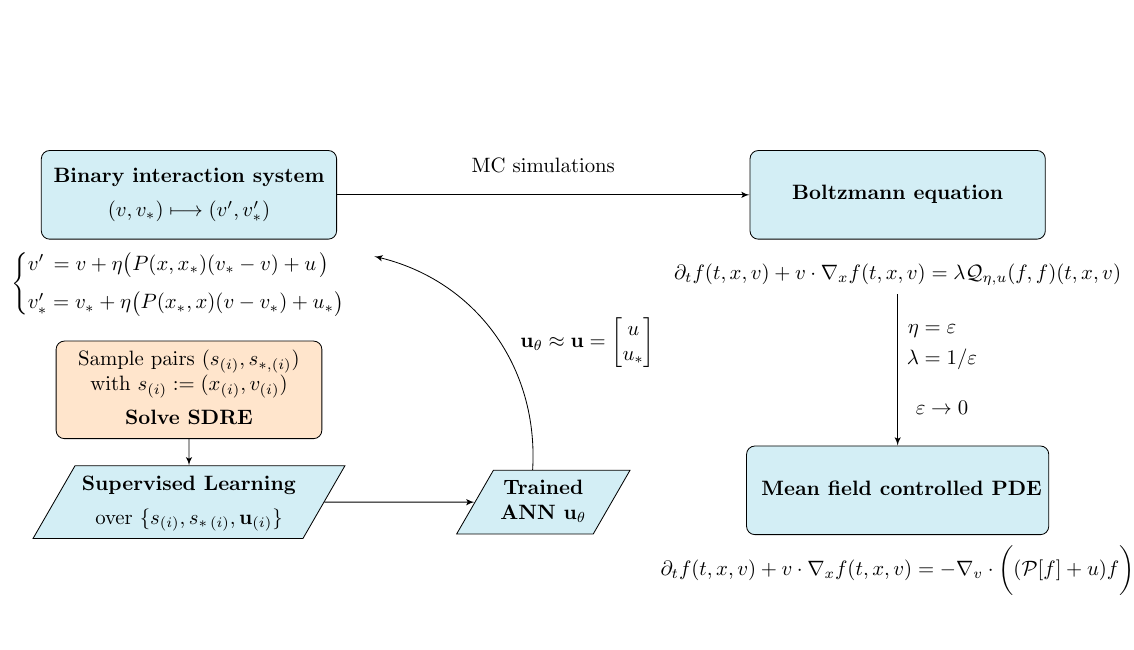}
        \caption{A diagram summarizing the main building blocks of the proposed numerical method. The mean field optimal control problem is approximated by a kinetic equation, which relies on sampling interactions between pairs of agents.  To reconstruct the complete microscopic dynamical system, these interaction equations are coupled with a free transport process (omitted here for simplicity).   Since every sample requires the solution of a binary optimal control problem, we resort to supervised learning techniques to build a feedback control law which can be called in the sampling step without solving online optimization problems.  We denote by $\mathbf{s}_{(i)}:=(x_{(i)},v_{(i)})$ the state of the $i-$th agent.}
    \label{fig:diag}
\end{figure}

Assuming a high interaction frequency between agents requires, at each evolution step, the solution of a large number of reduced 2-agent \eqref{Cost}. For efficiency purposes, we rely on supervised learning to build an \emph{Artificial Neural Networks} (ANN) approximating the solution of the binary OCP, thus circumventing its online solution at every sampling instance. For training, we generate synthetic data via the solution of the binary OCP using a discrete time \emph{state-dependent Riccati equation} (SDRE) approach \cite{SDREapproach,SDRE_survey2}. The overall procedure is outlined in Figure \ref{fig:diag}.

\paragraph{Related literature and contributions}
Simulation and control of high-dimensional ABMs is a longstanding and challenging topic. The modeling of such systems has been extensively studied within the kinetic research community, aiming at reducing the computational complexity of simulations \cite{Bird,Nanbu,plasma_kinetic}. The formulation of a kinetic description for evolutionary models has inspired a flourishing literature in ABMs \cite{CS_Boltzmann,kinetic_wealth,kinetic_uncertanty,Boltzmann_traffic}.
In this direction, the authors of \cite{binary_interactions} proposed a Monte Carlo approach based on binary collision dynamics \cite{Kinetic_MC}, inspired from plasma physics.  The proposed methodology addresses the integration of the mean field formulation of ABMs by means of a Boltzmann scheme, allowing to retrieve the mean field evolution of the system as a limit of 2-agents sub-problems. Such complexity reduction  has motivated the extension of this numerical scheme to fit the optimal control framework \cite{kinetic_opinion,kineitc_control_uncert,kinetic_control_CS}.  

{ In \cite{mfc_hierarchy,kinetic_opinion}, consistency was shown between the feedback controlled, non-local, mean field equation \eqref{mfpde} and the homogeneous non-linear Boltzmann equation \eqref{pv} with controlled binary interactions in the quasi-invariant scaling, and }where the control action has been optimized by means of \emph{model predictive control} (MPC) \cite{MPC_book}, or via \emph{dynamic programming} (DP) \cite{DP_book}. Similar approaches can be found also in \cite{albi2014boltzmann} for leader-follower multi-agent systems. In the MPC case, the controller is to be considered sub-optimal, as it is designed to optimize the OCP up to a reduced horizon. DP instead, leads  to optimal solutions for the binary OCP, but it requires the solution of a first-order nonlinear Hamilton-Jacobi-Bellman (HJB) PDE cast in the state-space of the system, which can be of arbitrarily high dimension. 

Over the last years, the solution of high-dimensional HJB-PDEs has been addressed with a number of different numerical approaches \cite{AKK,SKK,tensor22,visHJB,oster2024comparison}, and a flourishing literature on  Artificial Neural Networks (ANNs) . ANN methods are differentiated between unsupervised learning techniques \cite{dnn_hjb1,meng2024physicsinformed,PINNs_control,PINNs_new,onken} and supervised ones \cite{Kang_algorithms,ABK2021,newNN}. Examples of deep learning algorithms for solving PDEs can be found in \cite{NN_boltzmann} for high dimensional Boltzmann equations, or in \cite{pinns1,dgm} for more general applications. 

The main contributions of our work, inspired by the aforementioned results, can be summarized as follows:
\begin{itemize}
    \item Aiming at a further reduction in computational complexity with respect to \cite{kinetic_opinion,mfc_hierarchy}, we rely on data-driven approximation models fitted in a supervised learning fashion over synthetic data for the reduced binary OCP. For this, we consider candidate approximation models of both feed-forward (FNN) and recurrent (RNN) neural network type. 
    The architecture of such fitted approximation models conveniently allows for batch evaluations of data, meaning that at each time step the controls can be computed for all the sampled 2-agents sub-systems at once. The precision of the fitted models and their efficiency as the dimensionality of the problem increases have been assessed through numerical tests.

    \item The dataset guiding the ANNs training phase is collected from synthetic data, obtained from the solution of the infinite horizon \eqref{Cost} reduced to the binary case. Unlike the approaches taken in \cite{mfc_hierarchy}, we circumvent the solution of the HJB equation associated to the binary OCP, which would be unfeasible even in moderate agent dimensions. Instead, we rely on a procedure for synthesizing nonlinear feedback controls that combines elements from both DP and MPC: the State-Dependent Riccati equation approach. In particular, we will focus on discrete time settings \cite{discrete_sdre, D-SDRE_modelpredictive}.
    
    \item A first, intuitive, choice is to consider the feedback control as target variable of the approximation task, as done in \cite{MPC_book,ABK2021,Kang_algorithms,Kunish_feedbackNN,Data-TT}. Furthermore, we also train networks for approximating the whole of the controlled right hand side of the discrete-time binary controlled dynamics. For uncontrolled systems, approximation techniques has been studied for detection and approximation of interaction kernels $P(\cdot,\cdot)$ \cite{Maggioni_interactionrules,Maggioni_int_kernels}. However, to our knowledge, this work has not been yet extended to the controlled framework. 
\end{itemize}

The rest of the paper has been organized as follows. {  In Section 2 we construct an approximation for the sub-optimal control of the mean field PDE \eqref{mfpde} from a kinetic viewpoint. This leads to a Boltzmann description of the dynamics \eqref{pv}, which is then proven to converge to its mean field counterpart when performing a proper scaling of the frequency and strength of interactions}. In Section 3 we discuss how to approximate the evolution of system according to the Boltzmann dynamics for the distribution of agents via Monte Carlo simulation: from the current system configuration, we sample a pool of agents, which are then randomly coupled according to statistics of interaction. The post-interaction agents are then considered, so that their sampling distribution models the updated population density. Following this scheme, the evolution of the system is tracked at discrete times. This requires the computation of the feedback control in discrete-time settings, which is addressed in Section 4 by means of a discrete time state-dependent Riccati Equation approach. In Section 5, we discuss a supervised learning approximation to the solution of the binary OCP. Numerical tests are presented in Section 6.

\section{A Boltzmann formulation for the control of mean field dynamics}
The accurate modeling of self-organization phenomena and optimal control in ABMs requires a large number of interacting individuals, implying the need for the solution of a very high dimensional optimal control problem, which often comes at a prohibitively expensive computational cost. An alternative way to address this problem is to model instead the distribution function $f(t,x,v)$ describing the density of individuals having state variable  {$(x,v)\in {\R^d\times\R^d}$} at time $t \geq 0$. The evolution of $f(t,x,v)$ can be characterized by a kinetic equation accounting for the motion of individuals undergoing pairwise interactions, as modeled in \eqref{boltzmann}. Thus, the mean field dynamics can be retrieved by suitable scaling of the interactions, also referred as {\em quasi-invariant scaling}, or {grazing collision limit} \cite{Kinetic_MC, villani1998new}.  In particular, the quasi-invariant limit consists of considering an interaction regime where low intensity interactions occur with high frequency. {  In this regime, the density $f(t,\cdot,\cdot)$ is expected to converge pointwise in time to the solution of the mean field controlled PDE \eqref{mfpde}, which corresponds to a sub-optimal solution of the mean field optimal control problem \eqref{mf}.}

\subsection{Binary controlled dynamics}
We denote  {by $(x,v),(x_*,v_*)\in {\R^d\times\R^d}$} the position and velocity states of two agents in the population, and we assume that they modify their velocity states according to binary interaction maps $(v,v_*)\mapsto(v',v_*')$  as follows
\begin{equation}\label{post_int_dyn}
	\begin{aligned}
		v' &= v + \eta P(x,x_*)(v_*-v) \,+ \eta u,\\
		v_*' &= v_* + \eta\, P(x_*,x)(v-v_*) +\eta u_*\,,
	\end{aligned}
\end{equation}
where $\eta$ is the strength of interaction, $\bu = (u,u_*)^{\top}$ is the forcing term associated to interaction between agents. 
The goal of such external influence is to steer a couple of agents toward consensus. In particular, we will address the optimal control $\bar u$ as the solution of an infinite horizon binary optimal control problem, as follows
\begin{equation}\label{cost}
    \tilde \bu = \underset{\bu\in U_{\textrm{adm}}}{\argmin} \bigg\{\int\limits_{0}^{+\infty} \|v-\Bar{v}\|^2 + \|v_*-\Bar{v}\|^2 + \gamma\,\big(\mathcal{G}(u) + \mathcal{G}(u_*)\big)\, dt\bigg\},\,
\end{equation}
 {where $U_{adm}=L^{\infty}([0,+\infty[;\R^{2d})$ is the space of admissible controls.} The numerical procedure for the feedback control synthesis will be addressed in section \ref{discrete feedback}.

Furthermore, we model the evolution in time of $f(t,x,v)$ with a kinetic integro-differential equation of Boltzmann type \cite{mfc_hierarchy,MC_boltzmann,Povzner62}
\begin{equation}\label{boltzmann}
	\partial_t f(t,x,v)  +v\cdot \nabla_x f(t,x,v) = \lambda  \mathcal{Q}_{\eta,u}(f,f) (t,x,v)\,,
\end{equation}
where the parameter $\lambda $ encodes the interaction frequency, and the operator $Q_{\eta,u}(f,f)$ accounts for the gain and loss of particles with state $(x,v)$ at time $t$
	\begin{align}\label{eq:QBoltz}
		\mathcal{Q}_{\eta,u}(f,f) = \mathcal{Q}^{+}_{\eta,u}(f,f) - \mathcal{Q}^{-}_{\eta,u}(f,f). 
	\end{align}
In particular, 
we can express respectively the gain and loss operators as follows
\begin{subequations}
	\begin{align}\label{eq:QBoltz+}
		\mathcal{Q}_{\eta,u}^+(f,f)(t,x,v) &= \int\limits_{\R^d\times\R^d}\dfrac{1}{\mathcal{J}_{\eta,u}} f(t,x,'v)f(t,x_*,'v_*)\,dx_*\,dv_*\,,\\\label{eq:QBoltz-}
		\mathcal{Q}_{\eta,u}^-(f,f)(t,x,v) &= \int\limits_{\R^d\times\R^d}f(t,x,v)f(t,x_*,v_*) \,dx_*\,dv_*\,,
	\end{align}
 \end{subequations}
where $('v,'v_*)\longmapsto(v,v_*)$ are the pre-interaction velocities that generate the couple $(v,v_*)$, while $\mathcal{J}_{\eta,u}$ is the Jacobian of the binary interactions map \eqref{post_int_dyn}. 
 {To avoid the presence of the jacobian we can
introduce a test function $\varphi(x,v)\in C^2_0(\R^d\times\R^d)$ and express the Povzner-Boltzmann operator in the weak form as 
\begin{align}\label{eq:weak_boltz}
\lambda\langle \mathcal{Q}_{\eta,u}(f,f),\varphi\rangle = \lambda\iint_{\R^{2d}\times\R^{2d}}\left(\varphi(x,v')-\varphi(x,v)\right)ff_*\,dx_*\,dv_*\,dx \,dv,
\end{align}
where $v'$ denotes the post-interaction velocity as in \eqref{post_int_dyn}, and and we adopted the shorten notation $f = f(t,x,v)$, and $f_* = f(t,x_*,v_*)$. 
}
\begin{rmk}
 {
Here we considered the unbounded domain in space and velocity, where the operator  \eqref{eq:QBoltz} accounts interactions with constant collision frequency $\lambda$, analogously to the Boltzmann equation for Maxwell molecules. In general, the interaction frequency is ruled by a non-linear kernel, where in the Povzner-Boltzmann model it reads as follows
\begin{align}\label{eq:QBoltz_Povz}
		\langle\mathcal{Q}_{\eta,u},\varphi\rangle = \iint_{\R^{2d}\times\R^{2d}} B(x,x_*,v,v_*)\left(\varphi(x,v')-\varphi(x,v)\right)f\,f_*\,dx_*\,dv_*\,dx\,dv,
	\end{align}
 where the kernel $B(\cdot)$ in the Povzner approach \cite{Povzner62} considers non-local interactions among particles.}
 
  {
If we consider a bounded domain, further conditions should be included in \eqref{boltzmann}. Of particular interest is the case, when the velocity space is bounded, e.g. $v,v_*\in\Omega_v\subset \R^d$. In this situation, we can directly impose in \eqref{eq:QBoltz_Povz} that the boundary is satisfied introducing a non-linear interaction kernel $B$ such as $B\equiv\chi(v'\in\Omega_v)\chi(v'_*\in\Omega_v)$. This latter choice automatically satisfies the  boundaries, however such nonlinear kernel causes major difficulties when we are interested in studying asymptotic properties of the model \eqref{boltzmann}, such as mean field approximations. Alternatively, following the approach proposed in \cite{Kin_Toscani}, one can ask that the discrete interaction \eqref{post_int_dyn} preserves the boundary. In this setting, such requirement is possible by properly designing the admissible space of controls $U_{\textrm{adm}}$ in such a way that the boundary cannot be violated. 
We also observe that the control acts by forcing the state within an admissible position in the bounded domain. Thus, to preserve numerically the bounds a possible strategy is to require a sufficiently small penalty parameter $\gamma>0$ in the cost functional \eqref{cost}, as for example shown in \cite{Boltzmann_traffic}, where an explicit form of the control is obtained.
}
\end{rmk}


\subsection{Consistency with the mean field formulation}
Here we focus on the consistency of the Boltzmann operator \eqref{eq:QBoltz} with a mean field controlled dynamics {  of the type \eqref{mfpde}}, in particular introducing a {\em quasi-invariant optimality limit} we can regularize such operator, considering a regime where interactions strength is low and frequency is high. This technique, analogous to the grazing collision limit in plasma physics, has been thoroughly studied in \cite{villani1998new} and specifically for first order models in \cite{cordier2005kinetic,opinion_Toscani}, and allows to pass from Boltzmann equation \eqref{boltzmann} to a mean field equation \cite{kinetic_opinion,albi2014boltzmann}.
In what follows we consider the change of notation for the control in the binary dynamics \eqref{post_int_dyn}
$$u\to u_{\eta}(x,v,x_*,v_*),\quad u_*\to u_{\eta}(x_*,v_*,x,v)$$
to give explicit dependence on the parameters and the state variables, since we focus on feedback type controls, and we introduce the following  assumptions
\begin{itemize}
\item[$(i)$]  the system \eqref{post_int_dyn} constitutes an invertible changes of variables from $(v,w)$ to $(v',w')$;
\item[$(ii)$] there exists an integrable function $u(x,v,x_*,v_*)$ such that the following limit is well defined
\begin{align}\label{eq:limK}
\lim_{\eta\to 0}u_\eta(x,v,x_*,v_*) = u(x,v,x_*,v_*).
\end{align}
\end{itemize}
Hence, in accordance with \cite{CS_Boltzmann} and with \cite{mfc_hierarchy}, where an analogous of this argument is furnished for the controlled dynamics in stochastic settings,  we state the following theorem

\begin{thm}\label{thm:grazing}
{ 
Consider the Boltzmann-type equation \eqref{boltzmann}, with $\eta,\lambda>0$ and the control $u_{\eta}\in U_{\textrm{adm}}$, where $U_{\textrm{adm}}$ is the class the admissible controls.
Furthermore, assume the kernel function $F_\eta(\cdot)\in L^2_{loc}(\R^{2d}\times\R^{2d})$ for all $\eta> 0$, where $$F_\eta(x,v,x_*,v_*)=P(x,x_*)(v_*-v) + u_\eta(x,v,x_*,v_*),$$
and define the parameter $\varepsilon >0$ to introduce the quasi-invariant scaling as follows
\begin{equation}\label{eq:scaling}
\eta = \varepsilon, \qquad \lambda ={1/\varepsilon},
\end{equation}
which links the strength and  the frequency of the interactions in the Boltzmann-type equation \eqref{boltzmann}.
}
Thus, if $f^\varepsilon(x,v,t)$ is a solution for the scaled equation  \eqref{boltzmann},
for $\varepsilon\to0$  $f^\varepsilon(t,x,v)$ converges pointwise, up to a subsequence, to  $f(t,x,v)$  where $f$ satisfies  the following {  controlled} mean field equation,
\begin{align}\label{eq:FP}
\partial_t f +v\cdot \nabla_x f = -\nabla_v \cdot\left((\mathcal{P}[f] + \mathcal{U}[f])f\right)
\end{align}
with initial data $f^0(x,v)=f(0,x,v)$ and 
where $\mathcal{P}[\cdot]$ represents the interaction kernel  \eqref{operator} and the control is such that
\begin{align}\label{eq:kernelK}
\mathcal{U}[f](t,x,v)  = \int_{\R^d\times\R^d}u(x,v,x_*,v_*)f(t,x_*,v_*)\,dx_*\,dv_* 
\end{align} 
where $u(x,v,x_*,v_*)$ is defined as the limiting value in \eqref{eq:limK}. 
\end{thm}
We report the proof of this result in the \ref{AppendixA}, as a reformulation of the result proposed in \cite{mfc_hierarchy}.
\subsection{Full state controlled binary dynamics}
A further generalization consists in introducing a binary interaction dynamics where, differently from the previous section, we consider a binary exchange of information for the full states of the agents $ {s} =(x,v)^\top,~ {s}_* = (x_*,v_*)^\top \in \R^{2d}$.
The interacting process defining 
$({{s}},{{s}}_*)\to ({s}',{{s}'_*)}$ is described by a binary rules of the following type
\begin{equation}\label{eq:genbinary}
\begin{aligned}
{s}' &= {s}+ \eta G({s},{s}_*) + \eta  { H} u\cr
{s}_*' &= {s}_* + \eta G({s}_*,{s}) + \eta  { H} u_{*}
\end{aligned}
\end{equation}
where the controls are respectively $u = u({s},{s_*}), u_{*} = u_*(s,s_*)$, and the operators $G$ and $H$ are defined accordingly to \eqref{post_int_dyn} as follows 
\begin{align}\label{eq:ABmatrix}
G({s},{s}_*) =
\begin{pmatrix}
\mathbb{0}_{d} &        \mathbb{I}_{d}\\
\mathbb{0}_{d}&        -P(x,x_*)\mathbb{I}_{d}\end{pmatrix}
\begin{pmatrix}
x
\\
v
\end{pmatrix}
+
\begin{pmatrix}
\mathbb{0}_{d} &        \mathbb{0}_{d}\\
\mathbb{0}_{d} &        P(x,x_*)\mathbb{I}_{d}\end{pmatrix}
\begin{pmatrix}
x_*
\\
v_*
\end{pmatrix},\qquad  {H} = \begin{pmatrix}\mathbb{0}_{d}\\ \mathbb{I}_{d}\end{pmatrix}.
\end{align}

 {Hence,  the kinetic density $f(t,{s})$ evolves according to a Boltzmann-type model of the following type
\begin{align}\label{boltzmann_fullstate}
\partial_t f(t,s) = \lambda \hat{\mathcal{Q}}_{\eta,u}(f,f)(t,s),
\end{align}
with frequency $\lambda>0$, and where, differently from the formulation in \eqref{boltzmann}, the transport term is now encoded in the interaction operator $\hat Q(f,f)(t,s)$.
The interaction operator $\hat Q(f,f)(t,s)$ is now defined as follows
\begin{align}\label{boltzmann_fullstate_operator}
\hat{\mathcal{Q}}_{\eta,u}(f,f)(t,s) = \int_{\mathbb{R}^{2d}}\left(\frac{1}{\hat{\mathcal{J}}_{\eta,u}} f(t,'s)f(t,'s_*) - f(t,s)f(t,s_*)\right) \, ds,
\end{align}
where $\hat{\mathcal{J}}_{\eta,u}$ is the jacobian of the binary rule \eqref{eq:genbinary}, with $('s,'s_*)$ the precollisional states.} 

 {
In this framework, we can retrieve the consistency with mean field model \eqref{eq:FP} via Theorem \ref{thm:grazing}, considering now as a kernel function the following 
\begin{equation}\label{eq:Feta_full}
F_\eta(s,s_*)=\eta G({s},{s}_*) + \eta  { H} u(s,s_*),
\end{equation}
and assuming $F_\eta(\cdot)\in L^2_{loc}(\R^{2d}\times\R^{2d};\R^{2d})$.
Then, in the quasi-invariant scaling \eqref{eq:scaling}, and taking the limit $\varepsilon\to 0$ \color{black} we expect pointwise convergence of the kinetic model \eqref{boltzmann_fullstate} to the following mean field equation \color{black} 
\begin{align}\label{eq:FP_full}
      \partial_t f(t,s) =- \nabla_{{s}}\cdot \left(f(t,s)\int_{\R^{2d}}\left( G({s},{s}_*)+ H u({s},{s}_*)\right)f(t,{s}_*)\,d{s}_*\right).
\end{align}
Recalling that $\nabla_{{s}} = (\nabla_x,\nabla_v)^\top$ and that the operators $G,H$ are defined as in \eqref{eq:ABmatrix}, we have that \eqref{eq:FP_full} is equivalent to the mean field equation \eqref{eq:FP}.}

 {
This result follows the same steps as the proof of Theorem \eqref{thm:grazing}, in this regard 
further details can be found in \ref{AppendixA}.
} {
We stress that the main differences, with respect to the results of the previous section, consist in treating the transport term as part of the interaction dynamics, and requiring that \eqref{eq:Feta_full} is $L^2_{loc}$. This last requirement, useful for the validity of the consistency Theorem \ref{thm:grazing}, is in general more restrictive. Nonetheless the full state binary interaction \eqref{eq:genbinary}, in a reformulated version, allow to provide an efficient control synthesis for the second-order dynamics. This aspect will be discussed in more details in Section \ref{discrete feedback}.
}

\section{Asymptotic Monte Carlo methods for constrained mean field dynamics}
We provide a fast simulation numerical scheme for the Povzner-Boltzmann-type controlled equation, in the asymptotic regime, reminiscent of Direct Simulaton Monte Carlo methods (DSMCs), used in plasma physics, \cite{Nanbu,Bird,binary_interactions} and later adapted to collective dynamics  \cite{binary_interactions,Kin_Toscani}. 

First of all, considering a splitting of the transport and collisional term of the Boltzmann-type equation \eqref{boltzmann} in the asymptotic regime in two different steps \cite{binary_interactions}: 
\begin{equation}\label{splitting}
    \begin{cases}
    \partial_t f = - v\cdot \,\nabla_x f \qquad \qquad &\text{\emph{transport}}\\[0.5ex]
    \partial_t f = {\varepsilon}^{-1} \,\mathcal{Q}_{\eta,u}(f,f). \qquad \qquad &\text{\emph{interaction}}
    \end{cases}
\end{equation}

The purpose of the splitting scheme \eqref{splitting} is to focus the discussion on the convergence of the collisional term, since the free transport process coincides with the mean field formulation and the Boltzmann one.

In the previous section, consistency between the time evolution of the agents' population has been proven between the mean field model for the population density function, and a scaled Boltzmann description for the dynamics. The latter modeling of $f(t,x,v)$ is guided by the reduced microscopic binary interactions \eqref{post_int_dyn} between agents. This motivates the resorting to Monte Carlo simulation techniques for the approximation of the population density function under the Boltzmann formulation \cite{Kinetic_MC,MC_boltzmann,binary_interactions,mfc_hierarchy}.

The evolution \eqref{splitting} of $f$ in $[0,T]$ can be modeled by means of a Forward Euler method, for which we define a time step $\Delta t$, and discrete times $t_n = n\Delta t$ for $n=0,\ldots,N_T$
\color{black}
\begin{equation}
  {\partial_t f} (t_{n},x,v)\, \approx  \frac{f(t_{n+1},x,v) - f(t_{n},x,v)}{\Delta t}.
\end{equation}
\color{black}

In the spirit of DSMC methods, we can design a stochastic simulation scheme where we sample $N_s$ agents from $f^0$ and we approximate the time variation of $\pt f$ as the evolution of the sampled distribution according to the post transport/interaction states.

Concerning the free transport process, the time evolution amounts to the exact free flow at time $t_n$ of sample particles 
$\left\{(x^{n}_i, v^{n}_i)\right\}_{i=1}^{N_s}$ and evolves as follows
\begin{equation}\label{transport}
    x^{n+1}_{i} = x^{n}_i + \Delta t v^{n}_i\, \qquad i=1,\ldots,N_s.
\end{equation}

The collisional term in the Boltzmann-type equation \eqref{splitting}, can be rewritten in dependence of the gain and loss components of the operator $\mathcal{Q}_{\varepsilon,u}$:
\begin{equation}\label{interaction}
    \partial_t f(t,x,v)   = \frac{1}{\varepsilon} \left(\mathcal{Q}_{\varepsilon,u}^{+}(f,f) (t,x,v) - \rho f(t,x,v)\right)\,,
\end{equation}
where $\rho>0$ represent the total mass
\begin{equation}
    \rho = \int\limits_{\R^{2d}}f(t,x_*,v_*)\,dx_*\,dv_*\,.
\end{equation}
Assuming $f$ to be a probability density function we will consider $\rho=1$. The Forward Euler scheme for \eqref{interaction}, with the notation $f^{n}=f(t_{n},x,v)$ reads
\begin{equation}\label{f_evolution}
    f^{n+1} = \bigg(1-\frac{\Delta t}{\varepsilon}\bigg)f^{n} + \frac{\Delta t}{\varepsilon}\,\mathcal{Q}_{\varepsilon,u}^+(f^n,f^n)\,,
\end{equation}
since $f^{n}$ is a probability density, thanks to mass conservation  $\mathcal{Q}_{\varepsilon,u}^+$ is again probability density function. Moreover, under the restriction $\Delta t \leq \varepsilon$, also  $f^{n+1}$ is also a probability density.  Equation \eqref{f_evolution} can be interpret as follows: an agent with state $(x,v)$ has probability $(1-\Delta t/\varepsilon)\in[0,1]$ to avoid collision with other agents in each time interval $[t,t+\Delta t]$. When, instead, the collision does happen (event with probability $\Delta t/\varepsilon \in [0,1]$),  the evolution follows the interaction law $\mathcal{Q}_{\varepsilon,u}^+$ described by the scaled binary interaction, then a sampled pair of agents at time $t_n$ 
$$(x,v):=(x_{i}^{n},v_{i}^{n}),\qquad (x_*,v_*):=(x_{j}^{n},v_{j}^{n})$$ 
evolves according to \eqref{post_int_dyn} as follows
\begin{equation}\label{binary_euler}
    \begin{aligned}
		v^{n+1}_{i} &= v^{n}_{i} + {\varepsilon}\,P\left(x^{n}_{i},x^{n}_{j}\right)\left(v^{n}_{j}-v^{n}_{i}\right) +  \varepsilon\, u\left(x^{n}_{i},v^{n}_{i},x^{n}_{j},v^{n}_{j}\right),\\
		v^{n+1}_{j} &= v^{n}_{j} + {\varepsilon}\,P\left(x^{n}_{j},x^{n}_{i}\right)\left(v^{n}_{i}-v^{n}_{j}\right) +  \varepsilon\, u\left(x^{n}_{j},v^{n}_{j},x^{n}_{i},v^{n}_{i}\right).
	\end{aligned}
\end{equation}

\color{black} In what follows we consider the asymptotic regime with $\Delta t = \varepsilon$, where the binary interactions \eqref{binary_euler} are equivalent to a Forward Euler scheme for the two agent dynamics, and the particles that interacts at each iteration are maximized according to the scheme \eqref{f_evolution}. The choice of $\varepsilon$ is of paramount importance to approximate consistently the mean-field model \eqref{eq:FP} according to Theorem \ref{thm:grazing}. Nonetheless, there is a trade-off between the number of samples $N_s$ and the size of $\varepsilon$, indeed it is possible to see that below a certain threshold $\varepsilon^\star_{N_s}$ there is no further improvement in approaching the mean field model. We refer the reader to \cite{binary_interactions,Kinetic_MC} for a detailed discussion of these methods, and to \cite{borghi2025wasserstein} for a rigorous convergence analysis of stochastic particle dynamics in this framework.


The asymptotic Monte-Carlo algorithm can be formalized as the procedure outlined in Algorithm \ref{alg:pseudoDSMC}, which is comparable to Algorithm 1 in \cite{mfc_hierarchy}. This approach differs from the conventional Nanbu scheme in plasma physics, which was originally developed for free particle systems lacking external control inputs, and it does not require the introduction of a mesh.
\color{black}

Up to this moment, the discussion overlooked the derivation of the forcing terms $u$ in the interaction \eqref{binary_euler}. These control variable is meant to be of feedback type, as they only depend on the current agents' states. Nevertheless, the discrete-time nature of the interaction map  \eqref{binary_euler}, embeds within the control variable a dependency on the time-step $\Delta t$, which in turn is related to the parameter $\varepsilon$, see for example \cite{mfc_hierarchy,kinetic_opinion}. In the following section, we address the solution of the discrete time formulation \eqref{transport}\eqref{binary_euler} of the binary interaction control problem.

\section{Infinite horizon optimal control of binary dynamics}\label{discrete feedback}
In this section we study the solution of the infinite horizon optimal control problem for the reduced 2-agent system. 

We recall some fundamental notions on optimal control for discrete-time systems \cite{OCS_book}, to then present its numerical approximation using a discrete-time State-Dependent Riccati Equation approach. 

\color{black}We begin by reformulating \color{black} \eqref{binary_euler} in a general control-affine discrete time systems of the form:
\begin{equation}\label{dynamicsc}
\bz(t_{n+1})=\mathcal A(\bz(t_n))+\mathcal B(\bz(t_n))\bu(t_n)\,,
\end{equation}
where $\bz(t_n)\in\R^\kappa$ and $\bu(t_n)\in\R^\mu$ denote the state of the system and the control signal at time $t_n=n\Delta t$, such that { $\bz(t_n) :=(s^n,s^n_*)^\top\in\R^{\kappa},$ with $\kappa=d$ or $2d$ depending on whether the agent dynamics are of first or second order, respectively. The control vector is defined as $\bu(t_n):=(u^n,u_*^n)^\top\in\R^{\mu}\,,$ with $\mu = 2d$.}

The state-to-state map $\mathcal A(\bz):\R^\kappa\to\R^\mu$ and the control operator $\mathcal B(\bz):\R^\kappa\to\R^{\kappa\times \mu}$ are assumed to be $C^1(\R^\kappa)$, satisfying $\mathcal A (\bO_\kappa)=\bO_\kappa$ and $\mathcal B(\bO_\kappa)=\bO_{\kappa\times \mu}$. Note that binary systems of the form \ \eqref{binary_euler} fit this setting.
Given $Q\in\R^{\kappa\times \kappa}$, $Q\succeq 0$ and $R\in\R^{\mu\times \mu}$, $R\succ 0$, we are interested in the infinite horizon optimal control problem 
\begin{equation}\label{ocp}
	\underset{\bu(\cdot)}{\min\,}\cJ(\bu;\bs):=\sum\limits_{n=0}^{+\infty} \bz(t_n)^{\top}Q\bz(t_n)\,+\, \bu(t_n)^{\top}R\bu(t_n)\,,
\end{equation}
subject to the dynamics \eqref{dynamicsc} with initial state $\bz(t_0)=\bs$. We look for a solution $\tilde \bu$ to \eqref{ocp} in feedback form, that is, an optimal control map $\tilde \bu:\R^\kappa\to\R^\mu$ which is expressed a function of the state, $\tilde \bu(t_n)=\tilde \bu(\bz(t_n))$. The computation of an optimal feedback law follows a dynamic programming argument, for which we define $V(\bs):\R^\kappa\to \R$ as the optimal cost-to-go departing from $\bs$:
\begin{equation}
	V(\bs)=\underset{\bu(\cdot)}{\min}\;\cJ(\bu;\bs)\,, 
\end{equation}
where $V$ satisfies the Bellman equation 
\begin{equation}\label{bellman}
    V(\bs) = \underset{\bu\in\R^\mu}{\textrm{min}} \bigg( \bs^{\top}Q \bs + \bu^{\top}R \bu + V\big(\mathcal A(\bs)+\mathcal B(\bs)\bu\big)\bigg)\,, \qquad \text{for all }\;\bs\in\R^\kappa\,.
\end{equation}
From this point onwards, as we work globally in the state space, $\bz$ and $\bs$ are treated indistinctly.

In order to better describe the difficulties related to the solution of this optimal control problem, we first focus on the linear quadratic case, where the optimal solution is computed via the discrete-time linear quadratic regulator (LQR). The optimality of the feedback law is ensured by a direct link between the LQR solution and the dynamic programming one, obtained from equation \eqref{bellman}. In non-linear settings, the parallelism between the nonlinear QR and DP is broken, together with the optimality of the solution. 

\subsection{Discrete-time Linear Quadratic Regulator}
The discrete time linear quadratic problem is a particular instance of the optimal control problem \eqref{ocp} when $\mathcal A(\bs) = A\bs,\; A \in \R^{\kappa\times \kappa}$, and $\mathcal B(\bs) = B \in \R^{\kappa\times \mu}$. Under these assumptions, we make the ansatz $V(\bs)=\bs^{\top}\Pi \bs$, so that the Bellman equation \eqref{bellman} becomes
\begin{equation}\label{lqr_hj}
    \bs^{\top}\Pi\,\bs =  \underset{\bu\in\R^\mu}{\textrm{min}} \bigg( \bs^{\top}Q\,\bs + \bu^{\top}R\,\bu + (A\bs + B\bu)^{\top}\Pi\,(A\bs + B\bu)\bigg)\qquad \text{for all }\;\bs\in\R^\kappa\,.
\end{equation}
Solving the equation above leads to an optimal feedback of the form
\begin{equation}\label{u_dare}
    \tilde \bu= \tilde \bu(\bs) = -(R + B^{\top}\Pi\,B)^{-1} B^{\top}\Pi\,A\bs\,,
\end{equation}
where $\Pi\in\R^{\kappa\times \kappa}$ is the unique positive definite solution of the discrete-time algebraic Riccati equation (DARE):
\begin{equation}\label{dare}
    \Pi = Q + A^{\top}\Pi\,A - A^{\top}\Pi\,B(R + B^{\top}\Pi\,B)^{-1}B^{\top}\Pi\,A\,.
\end{equation}

\begin{rmk} 
\color{black}
Consider a discrete first-order binary system as in \eqref{binary_euler} with constant interaction kernel $P(x,x_*)=p$, with individual states $v,v_*\in\R^d$, target $\tilde v=0$ and $\varepsilon=\Delta t$. In this case, the optimal control problem \eqref{ocp} reads
\begin{equation}\label{binary_ex}
\underset{\bu(\cdot)}{\min}~ \Delta t\sum\limits_{n=0}^{+\infty} \|v^n\|^2 + \|v_*^n\|^2 + \gamma\,(\|u^n\|^2+\|u_*^n\|^2)
\end{equation}
\begin{equation}\label{binary2disc}
      \text{s.t. }\;\begin{cases}
    v^{n+1} &= v^n \,+ {\Delta t}\big(p(v_*^n-v^n) + u^n\big)\\
    v_*^{n+1} &= v_*^n + {\Delta t}\big(p(v^n-v_*^n) + u_{*}^n\big)\,.
    \end{cases}  
\end{equation}
This formulation corresponds to a linear-quadratic control problem by setting $\bs=(v,v_*)^{\top}$, $\bu=(u,u_{*})^{\top}$ and
\[
A =\mathbb{I}_{2d}+\Delta t \tilde A\,,\quad \tilde A= \begin{pmatrix}
-p &  p\\
        p & -p
\end{pmatrix}\otimes \mathbb{I}_{d},\qquad
B =  \Delta t \mathbb{I}_{2d}
\qquad
Q=\Delta t \mathbb{I}_{2d}, \qquad R=\Delta t \gamma \mathbb{I}_{2d}.
\]
Similarly as in Theorem \ref{thm:grazing}, we are interested in the asymptotic limit for $\Delta t\to0$, in which case \eqref{dare} becomes
\begin{equation}\label{dare_ex}
 2\tilde A\Pi-\frac{1}{\gamma}\Pi^2  +\mathbb{I}_{2d}=\mathbb{0}_{2d}.
\end{equation}
Exploiting to the symmetric structure of the two agent interaction, $\Pi$ can be reduced to diagonal and off-diagonal components, see e.g. \cite{herty2015mean,albi2022moment}, obtaining the following

\begin{equation}\label{riccati_ex}
  \Pi = \begin{pmatrix}
\pi_{\textsc d } & \pi_{\textsc o }\\
 \pi_{\textsc o } & \pi_{\textsc d }
\end{pmatrix}\otimes \mathbb{I}_{d}\,,\qquad \pi_{\textsc d }=\frac{\sqrt{\gamma}}{2}\left(1-2\sqrt{\gamma}p+\sqrt{1+4\gamma p^2}\right)\,,\quad\pi_{\textsc o }={\sqrt{\gamma}}-\pi_{\textsc d }\,.
\end{equation}

The limiting optimal feedback law $\tilde\bu(\bs)=(\tilde u,\tilde u_{*})$ is given by 
\begin{equation}\label{eq:control_ex}
\tilde \bu(\bs)=-\frac{1}{\gamma}\Pi \bs=-\frac{1}{\gamma}\begin{pmatrix}
\pi_{\textsc d } v+\pi_{\textsc o }v_{*}\\
 \pi_{\textsc o } v + \pi_{\textsc d }v_{*}
\end{pmatrix}\,,
\end{equation}
from where it follows that the individual feedback laws are given by 
$$\tilde u(v,v_{*})=-\frac{1}{\gamma}(\pi_{\textsc d } v+\pi_{\textsc o }v_{*})=\tilde u_{*}(v_*,v),$$
Equivalently, notice that the discretized control computed in\eqref{u_dare} in this case reads as follows
 \[
\tilde\bu_{\Delta t}(\bs) = -(\Delta t \gamma \mathbb{I}_{2d} + \Delta t^2 \Pi)^{-1}\Delta t \Pi(\mathbb{I}_{2d}+\Delta t \tilde A)\bs
 \]
and in the limit for $\Delta t\to 0$ we obtain exactly the same expression as in \eqref{eq:control_ex}.
 \color{black}
In the limiting kinetic equation this corresponds to the mean field control model 
\begin{align}\label{eq:FP_ex}
\partial_t f +v\cdot\nabla_x f = -\nabla_v \cdot\left(f\int_{\R^{2d}}(p(v_{*}-v)-\frac{1}{\gamma}(\pi_{\textsc d } v+\pi_{\textsc o }v_{*}))f(t,x_*,v_*)\ dx_*\ dv_*\right)\,.
\end{align}
Notice that the structure of  $~\Pi\in\mathbb{R}^{2d\times 2d}$ in \eqref{riccati_ex} implies isotropic action of the control in the $d$-dimensional state space.
\end{rmk}

\subsection{Discrete-time Non-Linear Quadratic Regulator}
It is feasible to apply LQR to a linearization of the system \eqref{binary2disc}. However, if we want to compute a law accounting for nonlinearities, then we compromise the convenient connection between the DARE and the Bellman equation. Moreover,  the computational cost required for the solution consistently increases.

In order to extend LQR to the nonlinear case, we begin by re-arranging the difference equation \eqref{dynamicsc} in semilinear state-dependent form 
\begin{equation}\label{semilinear}
    \bz(t_{n+1})=A(\bz(t_n))\bz(t_n)+B(\bz(t_n))\bu(t_n)\,,
\end{equation}
where $\mathcal{A(\bz)}=A(\bz)\bz$  and $\mathcal B(\bz) = B(\bz)$. 
Note that this semi-linearization is not unique for systems of order greater than $1$. In what follows, we further assume pointwise controllability, i.e. $\forall \bz\in\Omega$ the pair $(A(\bz),B(\bz))$ is controllable. 

Similarly to the LQR design, the feedback control policy can be calculated as
\begin{equation}\label{u_sdre}
    \begin{aligned}
            \tilde\bu(\bz(t_n)) &= -K(\bz)\bz\\
            K(\bz) &= \big(B(\bz)^{\top}\Pi(\bz) B(\bz) + R\big)^{-1}B(\bz)^{\top}\Pi(\bz) A(\bz)
    \end{aligned}
\end{equation}
where the argument $\bz$ for the operators denotes dependency of the current state $\z(t_n)$  with discrete time scale $\Delta t$, and $\Pi(\bz)$ is the solution of a DARE with state-dependent coefficients (DSDRE):
\begin{equation}\label{dsdre}
    \Pi(\bz) = Q + A(\bz)^{\top}\Pi(\bz)\,A(\bz) - A(\bz)^{\top}\Pi(\bz)\,B(\bz)(R + B(\bz)^{\top}\Pi(\bz)\,B(\bz))^{-1}B(\bz)^{\top}\Pi(\bz)\,A(\bz)\,.
\end{equation}
\color{black} As discussed in \cite{ChangBentsman2013}, similarly to the continuous-time case, the discrete-time SDRE is derived by assuming a quadratic ansatz for the value function in the Bellman equation and by freezing the matrix $\Pi(\bz)$ at the current state, thus neglecting its variation between time steps.  \color{black}
The state dependency of the feedback operator $K(\bz) = K(\bz(t_n))$ suggests the need for consecutive sequential solutions $\Pi(\bz) = \Pi(\bz(t_n))$ of \eqref{dsdre} at every discrete time $t_n$ along a trajectory. Thus, for high dimensions of the state space $\Omega\subset\R^n$, the exact solution of the DSDRE \eqref{dsdre} comes at a cumbersome computational cost. Between the several numerical approaches that have been proposed to address this, we cite Taylor series method \cite{SDRE_taylor} and interpolation for $\Pi(\bz)$ \cite{SDREapproach}, while we refer the interested reader to \cite{SDRE_tutorial} for a more exhaustive review.  In what follows, we rely on a discrete version of the SDRE approach proposed in \citep{SDREapproach}.

\subsection{Discrete-time SDRE approach}
We aim at circumventing the computational challenge of solving \eqref{dsdre} at each time step, by realizing the DSDRE feedback law in a model predictive control fashion: given the current state $\bz(t_n)$ of the system, we assume the operator $\Pi(\bz)$ to be a positive definite matrix in $\Pi\in\R^{n\times n}$, meaning that \eqref{dsdre} reduces to its algebraic form \eqref{dare}, where the state dependencies in \eqref{dsdre} are neglected by accordingly freezing all the operators at the current configuration. 
The resulting feedback variable leads to a suboptimal approximation of the controlled trajectory between $t_n$ and $t_{n+1}$, after which the procedure is repeated by freezing the system at the updated state $\bz(t_{n+1})$. 

\vspace{0.5cm}
\begin{algorithm}[H]
\caption{MPC-DSDRE approach}\label{alg1}
$\bz(0) \gets \bs$\Comment*[r]{\small{initial condition}}
\For{$n=0,...,N_{T}$}{
$\bs\gets \bz(t_n)$\Comment*[r]{\small{read current state}}
$A \gets A(\bs)$, $\;B \gets B(\bs)$\Comment*[r]{\small{freeze the system accordingly}}
solve \eqref{dsdre} for ${\Pi}$\;
$\bu \gets -(R + B^{\top}\Pi B)^{-1} B^{\top}\Pi A\bs$\Comment*[r]{\small{approximated feedback control}}
$\bz(t_{n+1}) \gets \bs + \Delta t(A\bs + B\bu)$\Comment*[r]{\small{the control system evolves for $\Delta t$}}}
\end{algorithm}

\vspace{0.5cm}
The main computational bottleneck still persisting with this approach is the availability of a sufficiently fast solver for \eqref{dare}. Efficiency in retrieving the DARE solution is key in the settings under consideration.   The consistency between the mean field and the kinetic dynamics only holds for $\varepsilon\ll1$ which in turn is associated with high frequency sampling, i.e. large numbers of couple of interacting particles in the Monte Carlo simulation.   At each time step, we aim at conveniently approximate the controlled post-interaction states for a large number of paired agents, living in an arbitrarily high dimensional state space $\Omega\subset\R^n$. We address this task by means of supervised learning approximation, relying on models within the family of Artificial Neural Networks.

\section{Neural networks and supervised learning approximation of feedback laws}
The computational method proposed in this paper models the time evolution of the population density $f(t,x,v)$ according to the Forward Euler scheme \eqref{f_evolution}, where the updated density is to be approximated by the sampling distribution of \emph{controlled} post-interaction  states in a Monte Carlo fashion. 

In this section, we provide an efficient approximating feedback map for the controlled binary dynamics via Feedforward and Recurrent Neural Networks (FNN/RNN). We briefly define these architectures and supervised learning framework for training \cite{NN_approximation}. After that, we will focus on the problem of interest, discussing synthetic data generation for feedback laws.

\subsection{Feedforward Neural Networks}
\begin{figure}
\centering
    \includegraphics[width = 0.7\textwidth]{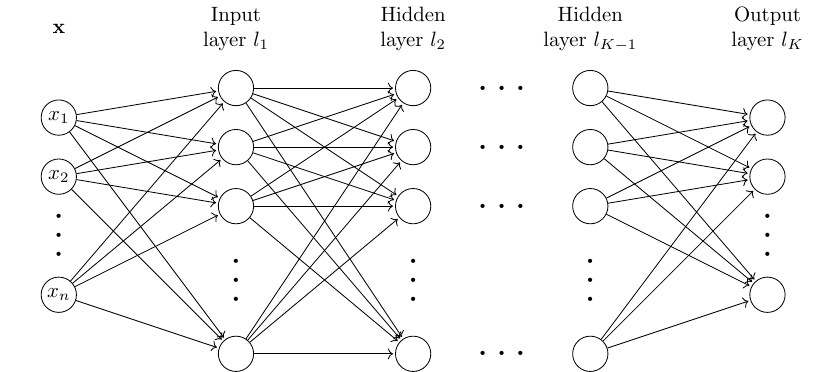}
    \caption{}
    \label{fig:DAG}
\end{figure}

These architectures approximate general functions $\psi$ through a sequence of layered transformations $l_1, \ldots, l_K$, resulting in $\tilde{\psi} = \psi_\theta \approx \psi$ where
\begin{equation}
\psi_\theta(\bx)=l_K \circ \ldots \circ l_2 \circ l_1(\bx),\qquad l_{k}(\bx) = \sigma_k(\bW_k \bx_k + \bb_k),
\end{equation}
with information flowing from the input nodes to the output ones in a unidirectional path, avoiding any cycles or loops, as shown in Figure \ref{fig:DAG}. Each layer applies a nonlinear activation function $\sigma_k$ element-wise to a linear transformation of its input $\bx_k$. Assuming the $k$-th layer to have $n_{k}$ neurons, the parameters $\bW_k\in\mathbb{R}^{n_{k-1}\times n_k}$ and $\bb_k\in\mathbb{R}^{n_{k}}$ represent the weight matrix and bias vector for the $k$-th layer, respectively. We use the same activation function for all neurons within a given layer, with $\sigma_1 = \sigma_K$ being the identity function.

In a supervised learning environment, the trainable parameters $\theta = \{\bW_k,\bb_k\}_{k=1}^K$ are then to be computed as minimizers of a suitable \emph{loss function} measuring the approximation error within a set of sampled data (\emph{training set}) $\mathcal{T} = \{\bx_{(i)},\psi_{(i)}\}_{i=1}^{N_s}$ where $\psi_{(i)} := \psi(\bx_{(i)})$:
\begin{equation}\label{training}
\begin{aligned}
        \theta^* = \underset{{\theta}}{\argmin}\;\cL(\psi, \psi_{{\theta}})\,,&\qquad \cL(\psi, \psi_{{\theta}})=\underset{\mathcal{T}}{\sum}l(\psi_{(i)},\psi_\theta(\bx_{(i)}))\,.
\end{aligned}
\end{equation}    
In particular, we consider 
\begin{equation}
    \mathcal{L}(\psi, \psi_{{\theta}}) = MSE(\psi, \psi_{{\theta}})\,, \qquad l(\psi_{(i)},\psi_\theta(\bx_{(i)})) = \dfrac{\|\psi_{(i)}-\psi_\theta(\bx_{(i)})\|^2}{N_s}\,.
\end{equation}
\color{black}
where we denoted by MSE the mean squared error as follows
$$
MSE(\psi, \psi_\theta) = \dfrac{1}{N_t} \underset{\mathcal{T}_v}{\sum} \dfrac{\|\psi_{(j)} - \psi_\theta(\bx_{(j)})\|^2}{\|\psi_{(j)}\|^2}.$$
\color{black}
The number of layers $K$, the number of neurons per layer, the set of activation functions $\{\sigma_k\}_{k=2}^{K-1}$ in the hidden layers, and eventual other parameters represent further degrees of freedom defining the \emph{architecture} $\Theta$ of the FNN. They are to be computed as result of the \emph{hyper-parameter tuning phase}, which amounts to choose, between trained models $\psi_\theta[\Theta]$ having different architectures, the one performing better to the error measure
\begin{equation}\label{hyp_tuning}
    \Theta^* = \underset{\Tilde{\Theta}}{\argmin}\;\sqrt{MSE(\psi, \psi_{\theta}[\Theta])},
\end{equation}
where $\mathcal{T}_v = \{\bx_{(j)},\psi_{(j)}\}_{j=1}^{N_v}$ is the \emph{validation set} such that $\mathcal{T} \nsubseteq\mathcal{T}_v$. 

\subsection{Recurrent Neural Networks}
Recurrent neural networks are a class of models which can be obtained from FNNs by allowing loop connections between layers.   This architecture is particularly well-suited for processing sequential data. However, given that our input lacks temporal structure, we focus on \emph{one-to-one} architectures. 

\begin{figure}
\centering
    \includegraphics[width = 0.4\textwidth]{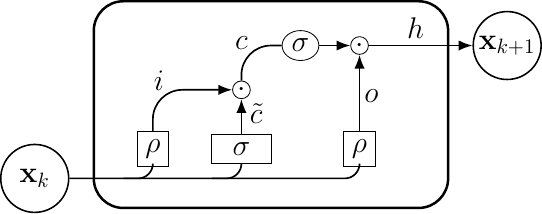}
    \caption{One-to-one LSTM cell, visualization of equation \eqref{LSTM_eq}}
    \label{fig:LSTM}
\end{figure}

Here we extend FNNs into RNNs by allowing hidden layers to be one-to-one \emph{Long Short Term Memory} (LSTM) cells which -- for every input $\bx$ -- generate output $h$  according to the following system of equations  
\begin{equation}\label{LSTM_eq}
    \begin{aligned}
    i &= \rho\big(W_i\, \bx + b_i\big)\\
    \Tilde{c} &= \sigma\big(W_c\, \bx + b_c\big)\\
    c &= i \odot \Tilde{c}\\
    o &= \rho\big(W_o\, \bx + b_o\big)\\
    h &= o \odot \sigma\big(c\big)
    \end{aligned} \qquad
    \begin{aligned}
    &\text{\emph{input gate}}\\
    &\text{\emph{candidate value}}\\
    &\text{\emph{cell value}}\\
    &\text{\emph{output gate}}\\
    &\text{\emph{final output}}
    \end{aligned} 
\end{equation}
where  (assuming $\bx_k\in\R^{n_k}$),   the weight matrices $ W_i,\,W_c,\,W_o \in \mathbb{R}^{n_{k}\times n_{k+1}}$ and the bias vectors $ b_i,\,b_c,\,b_o\in\R^{n_{k+1}}$ are the trainable parameters, $\sigma(\cdot)$ is the activation function, and $\rho(\cdot)$ is the recurrent activation function. The number $\kappa$ denotes the number of neurons in the cell. A visual representation of the flux of information through a LSTM cell can be found in Figure \ref{fig:LSTM}.

The training and architecture selection for this kind of RNNs consists of the same procedure described in \eqref{training} and \eqref{hyp_tuning}, with the exception that the $k$-th layer is replaced by a LSTM cell, mapping input $\bx_k$ to output $h=\bx_{k+1}$, with parameters $\bW_k = (W_i,W_c,W_o)^\top$ and $\bb_k = (b_i,b_c,b_o)^\top$.
\subsection{Synthetic data generation}
\color{black}
Both the training and hyper-parameter tuning phases for the neural network approximation rely on the availability of datasets containing input data points coupled with the corresponding target function evaluations. Since our goal is to approximate the feedback law for the reduced binary controlled dynamics, we consider two distinct approaches:

\begin{enumerate}
    \item[] \textbf{Control approximation}: The first approach treats the feedback control law $\bu$, i.e. we build a network $\bu_{\theta}(\bs, \Delta t) \approx \tilde{\bu}(\bs, \Delta t)$, where $\tilde{\bu}$ is the optimal control obtained from the SDRE approach described in Section~\ref{discrete feedback}.

    \item[] \textbf{State update approximation}: In this case, the neural network directly learns the mapping from current states to their controlled updates, i.e., $\bs'_{\theta}(\bs, \Delta t) \approx \mathcal{A}(\bs) + \mathcal{B}(\bs)\tilde{\bu}(\bs, \Delta t)$, bypassing the explicit computation of the control law.
\end{enumerate}

The training loss for both approaches uses the mean squared error (MSE) as defined in equation \eqref{training}. In Section~\ref{num_section}, we compare these two approaches in numerical tests, examining their precision and computational efficiency.
\color{black}

To create training data and establish a benchmark for precision, we solve the optimal control problem using the discrete-time SDRE approach.    The difference equations for the binary system at the discrete time $t_n$ hold as follows:
\begin{equation}\label{tilde_din}
\begin{aligned}
    &\begin{cases}
    x^{n+1} \,= x^{n} + \Delta t \,v^{n}\\
    x_*^{n+1} \,= x_*^{n}+ \Delta t \,v_*^{n}\end{cases}  \\
        &\begin{cases}
    v^{n+1} \,= v^{n} \,+ {\Delta t}\big[P(x^{n},x_*^{n})(v_*^{n}-v^{n}) + u(x^n,x_*^n,v^n,v_*^n)\big]\\
    v_*^{n+1} = v_*^n + {\Delta t}\big[P(x_*^{n},x^{n})(v^n-v_*^n) + u_{*}(x^n,x_*^n,v^n,v_*^n)\big]
    \end{cases}
\end{aligned}
\end{equation}
being $(x,v),(x_*,v_*)\in \Omega\subset\R^{2d}$ the states for a couple of interacting agents.  
We introduce the notation $\bs = (x,x_*,v,v_*)^{\top}\in\R^{4d}$, $\bu = (u,u_{*})^{\top}\in\R^{2d}$, for which the difference equation reads as in \eqref{dynamicsc}
\begin{equation}\label{dyn_s}
    \bs^{n+1} = \mathcal{A}(\bs^n) + \mathcal{B}(\bs^n)\bu^n\,,
\end{equation}
for $\mathcal{A}(\bs) = (\mathbb{I}_{4d}+\Delta t \,G(\bs))\bs$ and $\mathcal{B}(\bs) = \Delta t\,H(\bs)$, with  $G$ and $H$ defined as in \eqref{eq:ABmatrix}.
From this, comes the $\Delta t$-dependency of the feedback law $\bu$: the DSDRE parameters $\mathcal{A},\,\mathcal{B}$ are dependent of the time-step, as will be the DARE solution associated to the system frozen at the current configuration.

\begin{algorithm} 
\caption{Synthetic data generation}\label{alg_data}
\nl$\{\bs_{(i)},\Delta t_{(i)}\}_{i=1}^{N_s} \;\;i.i.d.$\Comment*[r]{\small{$N_s$ uniform samples in $\Omega^2\times[0,T]$}}
\nl\For{$i=1,...,N_s$}{
\nl $A \gets \mathbb{I}_{4d}+\Delta t_{(i)} \,G(\bs_{(i)})$;\\
\nl $B \gets \Delta t_{(i)}\, H(\bs_{(i)})$\Comment*[r]{\small{$G$, $H$ as in \eqref{eq:ABmatrix}}}
\nl solve \eqref{dsdre} for ${\Pi}$\;
\nl $\bu_{(i)} \gets -(R + B^{\top}\Pi B)^{-1} B^{\top}\Pi A\bs_{(i)}$\Comment*[r]{\small{feedback control}}
\nl$\bs'_{(i)} \gets \bs_{(i)} + \Delta t_{(i)}(A\bs_{(i)} + B\bu_{(i)})$\Comment*[r]{\small{post-interaction state}}
}
\nl $\mathcal{T}_\bu = \{\bs_{(i)},\Delta t_{(i)},\bu_{(i)}\}_{i=1}^{N_s}$\Comment*[r]{\small{dataset for $\bu_{\theta}$}}
\nl $\mathcal{T}_{\bs'} =\{\bs_{(i)},\Delta t_{(i)}, {\bs}_{(i)}'\}_{i=1}^{N_s}$\Comment*[r]{\small{dataset for $\bs'_{\theta}$}}
\end{algorithm}
 
 The generation of datasets $\mathcal{T}_{\bu},\,\mathcal{T}_{\bs'}$ for the supervised learning approximation of the feedback law $\bu_{\theta}$, and post-interaction states $\bs'_{\theta}$ is summarized in Algorithm \ref{alg_data}.
The number $N_s$ of samples in the datasets is discussed in Section \ref{num_section} specifically for each numerical test, whilst the reliability and generalizability of the sampled data are ensured through cross validation of the trained models over an unseen \emph{test set} $\mathcal{T}_t \nsubseteq \mathcal{T}\cup\mathcal{T}_v$ also generated through Algorithm \ref{alg_data}. The addition of the time step $\Delta t$ to the models' input allows to train the networks to provide accurate approximations even when applied to adaptive time-step integration techniques.
 
\begin{rmk}
	  As the controlled density evolution is carried by a Monte Carlo type method, the computational efficiency of evaluating a batch of controlled particle interactions is crucial. Our choice of control/state architectures is motivated by its favorable computational complexity. The evaluation cost of a single point through FNNs is $O(n)$, $n$ being the number of neurons in the network. Furthermore, in our non-sequential framework, the evaluation cost of a one-to-one LSTM cell with $n$ neurons is comparable to the one of a FNN with $4n$ neurons. This linear scaling in both architectures ensures efficient batch evaluation, as the computational cost grows proportionally with network complexity. This efficiency is particularly advantageous for large-batch evaluation, associated with $\varepsilon\ll1$.
\end{rmk}


\section{Numerical Tests}\label{num_section}
In this section we will assess the proposed methodology over consensus control problems for two high dimensional (both in $d$ and in $N$) ABMs. We aim at modeling the time evolution of the agents' distribution $f(t,x,v)$ via Monte Carlo simulation of the approximated binary post interaction states. In this section we will compare two different approaches: the approximation of the feedback control map $\bu_\theta$, acting in the binary system \eqref{tilde_din}, and the direct approximation of the controlled dynamics $\bs_\theta'$.  Since the approximation of the post-interactions positions $\bs'_x$ is direct, the neural network model is restricted to post-interaction velocities $\bs'_v$, that is $\bs_\theta' = (\bs_x',\bs_{v,\theta}')^\top$.$\,$  
Algorithm \ref{alg2} \cite{binary_interactions} summarizes the proposed numerical procedure for the density evolution once an approximated optimal feedback law $\bu_\theta \approx \tilde{\bu}$ has been constructed.
\begin{algorithm}
\caption{Monte Carlo simulation for Boltzmann dynamics}\label{alg2}
\nl$\{(x,v)_{i}\}_{i=1}^{N} \sim f(t_0,x,v)\;\;i.i.d.$\Comment*[r]{\small{$N$ samples from the initial distribution}}
\nl\For{$n=0,...,T$, $t_n = n\cdot \Delta t$}{
\nl select $\{(i_k,j_k)\}_{k=1}^{N//2}\;$ \Comment*[r]{\small{random pairs of agents without repetitions}}
\nl$\{\bs_{k}^{n} \}_{k=1}^{N//2},\;\; \bs_{k}^{n} = \big(x_{i_k}^{n},x_{j_k}^{n},v_{i_k}^{n},v_{j_k}^{n}\big)$\Comment*[r]{\small{$N//2$ couples of agents}}
\nl\For{$k = 1,...,N//2$}{
\nl$s^{n+1}_{k} \gets \mathcal{A}(\bs^{n}_k)\bs^{n}_k + \mathcal{B}(\bs^{n}_k)\Tilde{\bu}(\bs^{n}_k,\Delta t)$\Comment*[r]{\small{controlled dynamics \eqref{tilde_din}}}
\nl$\big(x_{i_k}^{n+1},v_{i_k}^{n+1},x_{j_k}^{n+1},v_{j_k}^{n+1}\big)\gets \bs^{n+1}_k$;}
}
\end{algorithm}
When a neural network $\bs_\theta'$ is built, we replace line $ 6$ in the algorithm with
\begin{equation}
    \bs^{n+1}_k\; \gets \;\bs_\theta'(\bs^{n}_k)\,.
\end{equation}
The assumption of constant time-step $\Delta t$ can be easily relaxed, as both the approximants $\bu_\theta$ and $\bs'_\theta$ take $\Delta t$ as input. 
In the following numerical tests, the NN training has been done via \emph{Adam} optimizer \cite{adam} with a learning rate $\alpha = 0.01$  over batches of $100$ samples.
\color{black}
\subsection{Test 1: Sznajd model}
To illustrate the relationship between the control strategy derived at the supervised learning kinetic approximation level and its mean-field counterpart, we consider a preliminary numerical example in a simplified setting. We focus on consensus control for a one-dimensional first-order Sznajd model for opinion dynamics \citep{Sznajd}. At the mean-field level, the control problem reads as a first-order version of \eqref{mf}: 
\begin{equation}\label{test0mf}
    \begin{aligned}
    \min_u \;&\mathcal{J}(f,u):=\int\limits_{0}^{T}\int\limits_\Omega\big(\|x-\bar{x}\|^2 + \gamma \|u\|^2\big)\,f(t,x_*)\,dx\, dt\,,\\
    &\text{subject to }\;    \partial_tf + \partial_x\bigg[(\mathcal{P}[f]+u)f\bigg] = 0\,,
    \end{aligned}
\end{equation}
where $\bar{x}$ denotes the target configuration, and the interaction operator is defined as
\begin{equation}
    \mathcal{P}[f](t,x) = \int\limits_\Omega P(x,x_*)(x_*-x) f(t,x)dx_*\,,\qquad P(x,x_*) = \beta (1-x^2)\,,
\end{equation}
for $\beta\in\R$. In particular, we consider $\Omega = [-1,1]$, representing the opinion space of a population of voters ranging between two extreme positions represented by $\{-1,1\}$, and a negative interaction coefficient $\beta = -1$, which induces polarization behavior. The control energy penalization coefficient is fixed to $\gamma = 0.05$.

This simplified framework enables the numerical solution of the mean-field optimal control problem, which would be computationally prohibitive in the higher-dimensional, second-order models considered in the subsequent numerical tests. 

At the binary interaction level, the problem reads
\begin{equation}
    \min_{\bu(\cdot)} \dfrac{\Delta t}{2} \sum_{n=0}^{+\infty}(\|\bx^n\|^2+\gamma\|\bu^n\|^2) \,,\,\, \text{s.t. }\,\,\begin{cases}
    {x}^{n+1} = x^n+\Delta t({\beta}(1-{(x^n)}^2)(x_*-x) + u^n)\\[0.5em]
    {x}_*^{n+1} = x^n+\Delta t({\beta}(1-(x_*^n)^2)(x-x_*) + u^n_*)
    \end{cases}
\end{equation}
for $\bu^n:=(u^n,u^n_*)^\top$ and $\bx = (x^n,x^n_*)^\top$. The synthetic data generation procedure associates $N_s=10^5$ uniform samples of interacting couples of agents $\bx^n_{(i)}\in\Omega$ to their associated dSDRE feedback control $\bu^n_{(i)}$ and the controlled state update $\bx^{n+1}_{(i)}$ obtained with interaction strenght $\varepsilon=0.01$, for  $i=1,\dots, N_s$. With a ratio of $80/20$, we designate those samples to form the training set $\mathcal{T}$ and the validation set $\mathcal{T}_v$ respectively. We test the following architectures:\\

    $u_\theta^{FNN}:$ $K = 1$ hidden layers with $100$ neurons, and  $\sigma(\bx) =log(exp(\bx) + 1)$\\[-10pt]
    
    $u_\theta^{RNN}:$ $K = 1$ LSTM cell with $100$ neurons, $\sigma(\bx) = tanh(\bx)$, $\rho(\bx) = (1 + exp(-\bx))^{-1}$\\[-10pt]
    
    $\mathbf{x'}_{\theta}^{FNN}:$  $K = 2$ hidden layers with $60$ neurons per layer, and $\sigma(\bx) = tanh(\bx)$\\[-10pt]
    
\begin{figure}
    \centering
    \includegraphics[trim=0cm 4cm 0cm 4cm, clip, width=0.45\linewidth]{ 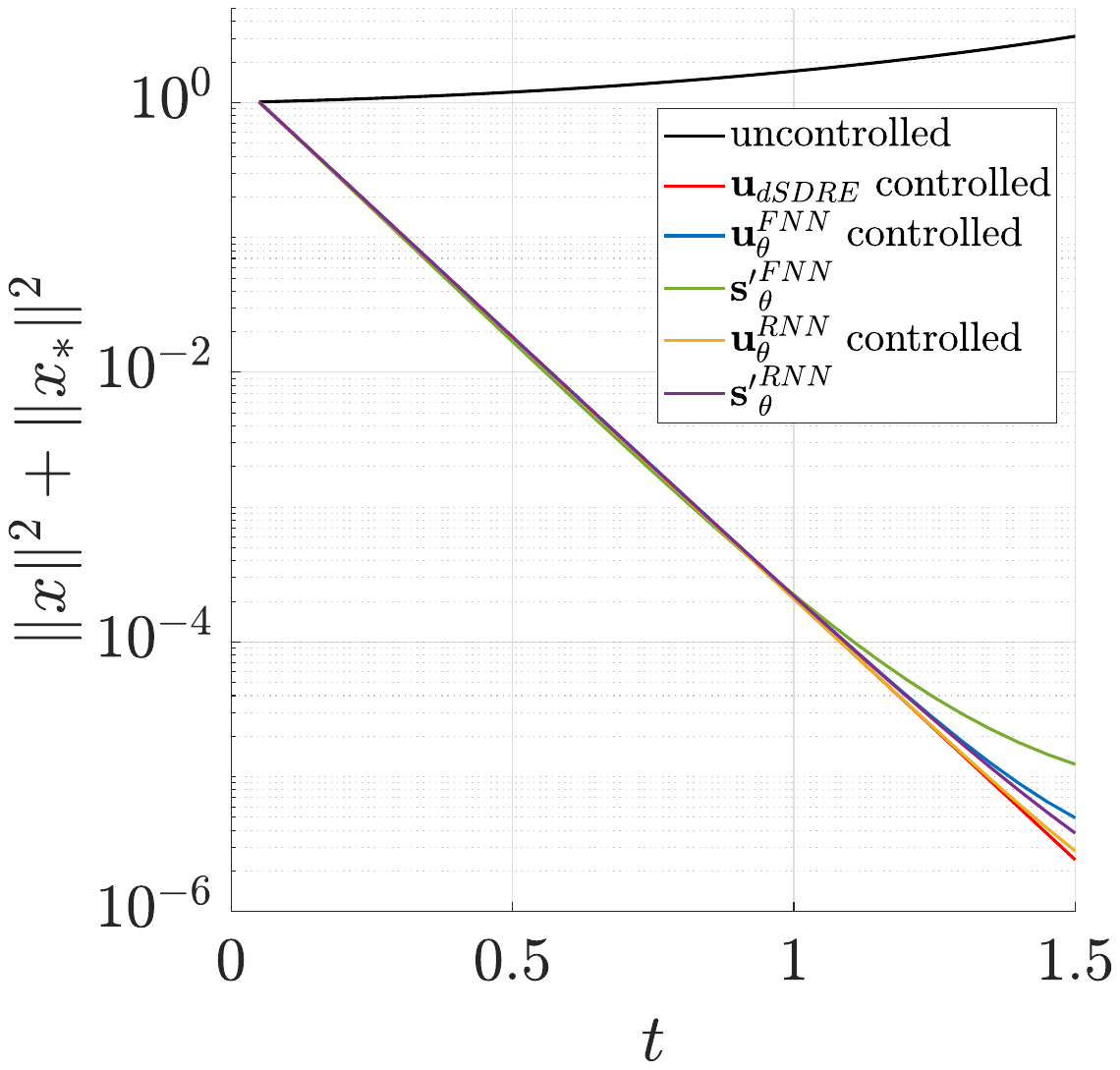}
    \caption{\color{black}Convergence to consensus for two randomly sampled interacting agents with states $x,x_*\in\Omega$, when controlled via dSDRE feedback law vs. the approximation models.}
    \label{fig:consensus_conv_test0}
\end{figure}

    $\mathbf{x'}_{\theta}^{RNN}:$ $K=1$ LSTM cell with $100$ neurons, $\sigma(\bx) = tanh(\bx)$,  $\rho(\bx) = log(exp(\bx) + 1)$\\[10pt]
The error comparison for the binary-controlled evolution of a randomly selected pair of interacting agents is presented in Fig.~\ref{fig:consensus_conv_test0}, where all approximation models demonstrate similar performance to the dSDRE control.

We conclude this preliminary example by quantifying the performance gap between the solution of \eqref{test0mf} and its neural network-accelerated kinetic approximation, both implemented with a time step $\Delta t = \varepsilon = 0.01$. This comparison offers a more concrete understanding of the suboptimality of the control strategy at the kinetic level with respect to the MFOC. The approximation of the mean field control law under the numerical procedure proposed in this paper is retrieved as
\begin{equation}\label{u_mf0}
    \int\limits_\Omega \bu_\theta(x,x_*)f(t,x_*)dx_*\approx u(t,x)\,.
\end{equation}
In Table~\ref{tab:mf_comparison}, we compare the MPC solution of problem~\eqref{test0mf} with the forward integration of the mean-field dynamics controlled via~\eqref{u_mf0}, where the control input $\bu_\theta$ is given by the previously trained neural networks $\bu_\theta^{FNN}$ and $\bu_\theta^{RNN}$. The optimality conditions arising from problem \eqref{test0mf} are approximated with a first-order semi-Lagrangian scheme, as done in \cite{mfc_hierarchy}, but neglecting diffusion. Moreover, to mirror the feedback nature of the control resulting from the numerical routine discussed so far, the solver is applied in an MPC fashion, where the optimality conditions are solved iteratively over a reduced time horizon $T=0.5s$ \cite{MPC_book}. Using the same semi-Lagrangian scheme, we simulate the forward dynamics controlled via the kinetic approximation models through~\eqref{u_mf0}.

\begin{figure}
    \centering
    \includegraphics[trim=0cm 3cm 0cm 3cm, clip, width=0.45\linewidth]{ 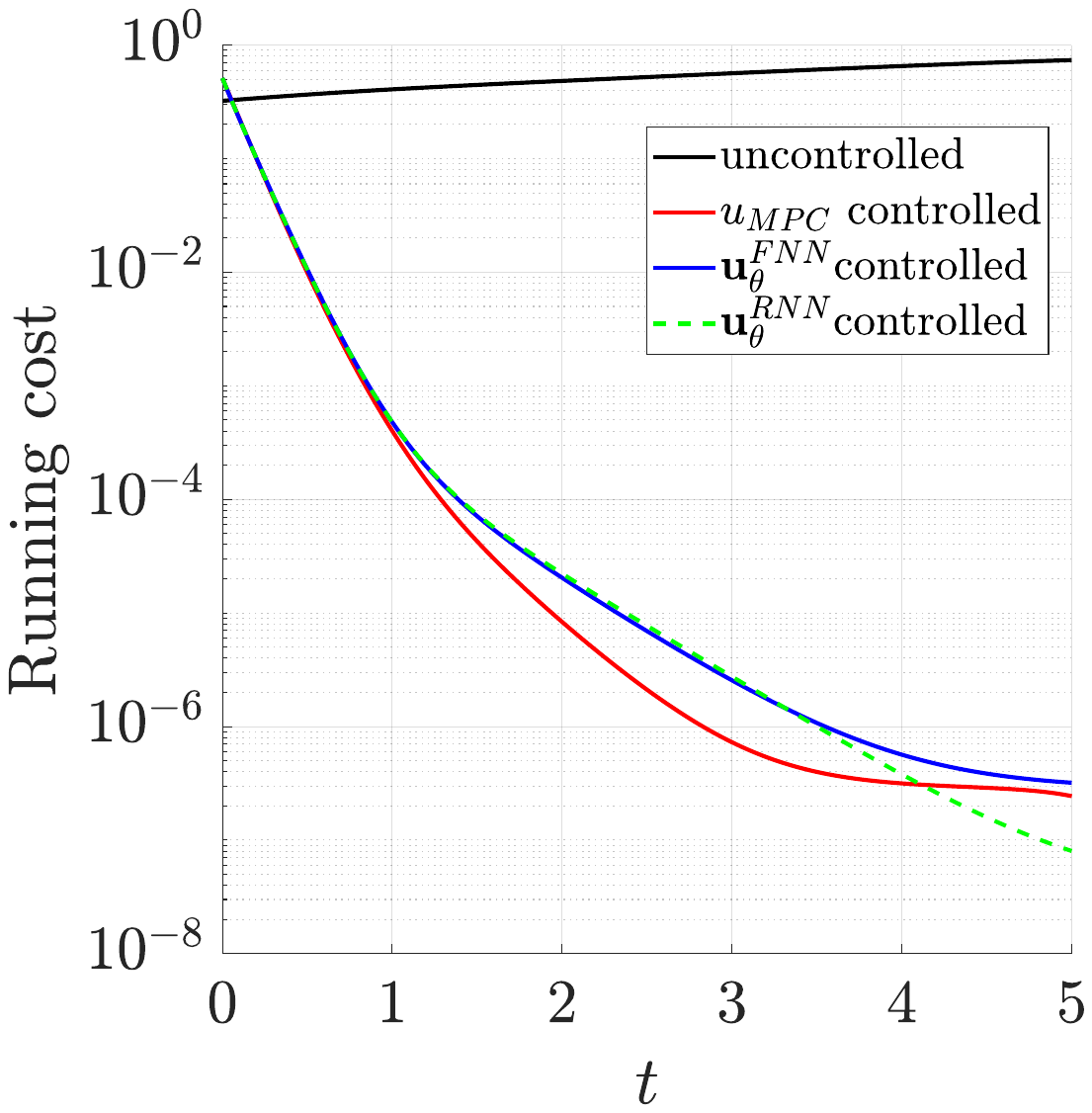}
    \caption{\color{black}Comparison of the time evolution of the running cost for the MFOC \eqref{test0mf}. }
    \label{fig:fullcostMF}
\end{figure}

The table reports comparable outcomes in terms of cost functionals, density evolutions, and control profiles across the three approaches. The mean-field control, denoted here by $u_{\text{MPC}}$, retains the best performance. It is worth noting that even $u_{\text{MPC}}$ should be regarded as close-to-optimal for an infinite horizon cost. This is because it is computed as the numerical solution of the first-order necessary optimality conditions, implemented in a receding horizon fashion over a sequence of finite-horizon MFOCs. In Figure \ref{fig:fullcostMF}, we compare the running cost evolution, i.e. the integral w.r.t. time in \eqref{test0mf}, which is consistent with the results reported in Table~\ref{tab:mf_comparison}.
\begin{table}[h!]\color{black}
\centering
\begin{tabular}{lccc}
& \textbf{$\bu_\theta^{RNN}$} & \textbf{$\bu_\theta^{FNN}$} & \textbf{$u_{MPC}$} \\
\midrule
&$\mathcal{J} = 0.063444$&$\mathcal{J} = 0.063504$&$\mathcal{J} = 0.063230$\\
\rotatebox{90}{\makebox[4.5cm][c]{$f(t,x)$}} &\includegraphics[trim=1cm 3cm 1cm 3cm, clip, height=.275\textwidth]{ 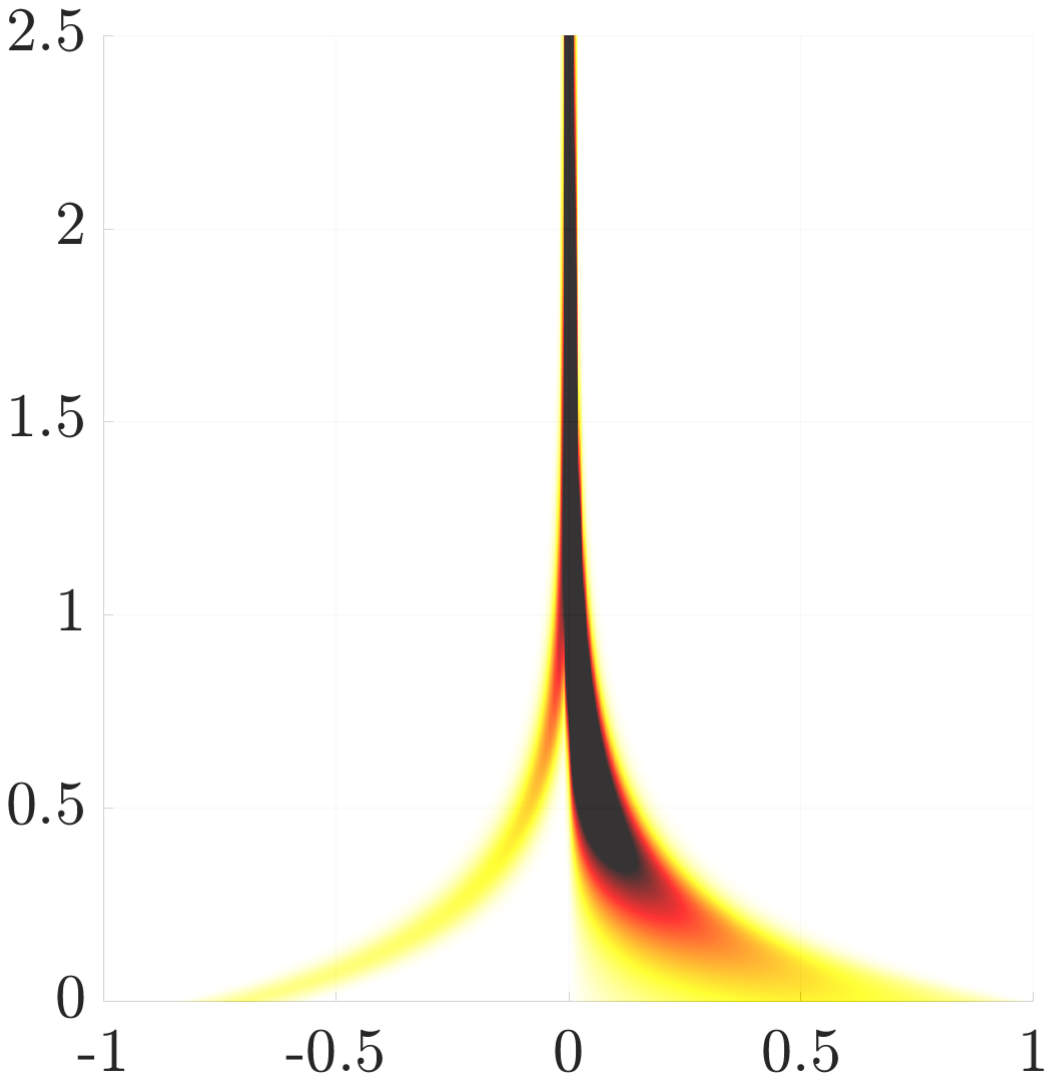} & \centering
\includegraphics[trim=1cm 3cm 1cm 3cm, clip, height=.275\textwidth]{ 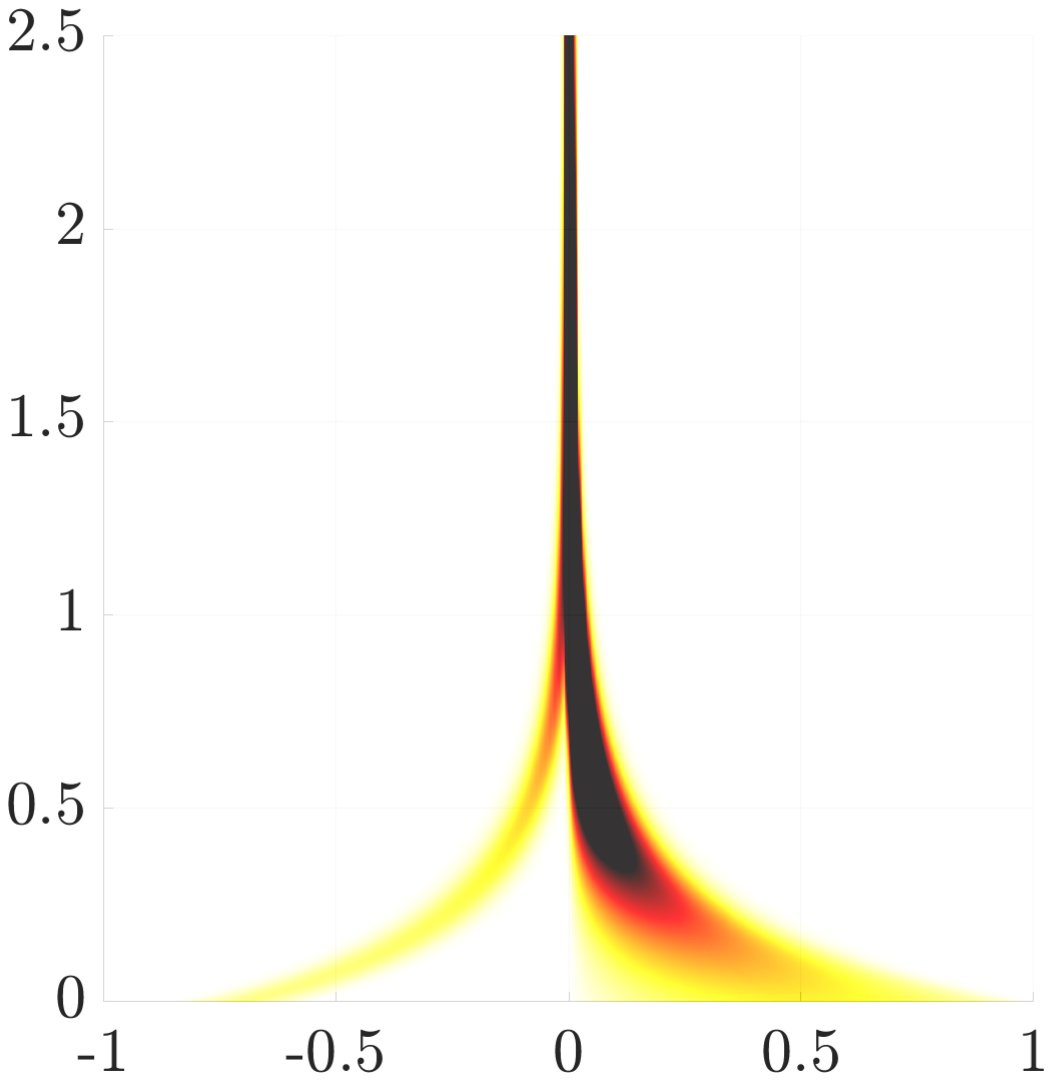} & 
\includegraphics[trim=1cm 3.5cm 1cm 3cm, clip, height=.275\textwidth]{ 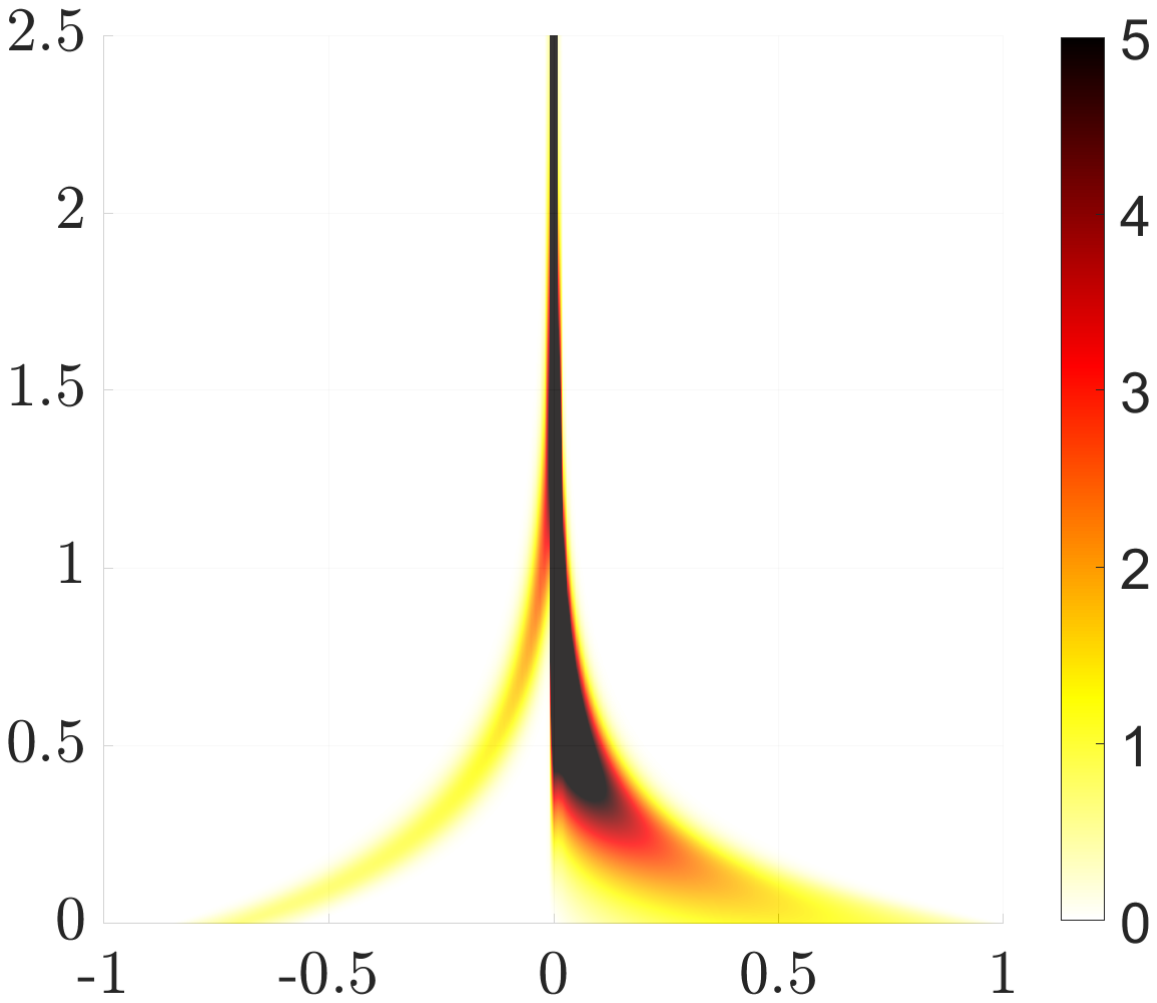} \\ 
\rotatebox{90}{\makebox[4.5cm][c]{$u(t,x)$}} & \centering \includegraphics[trim=1cm 3cm 1cm 4cm, clip, height=.275\textwidth]{ 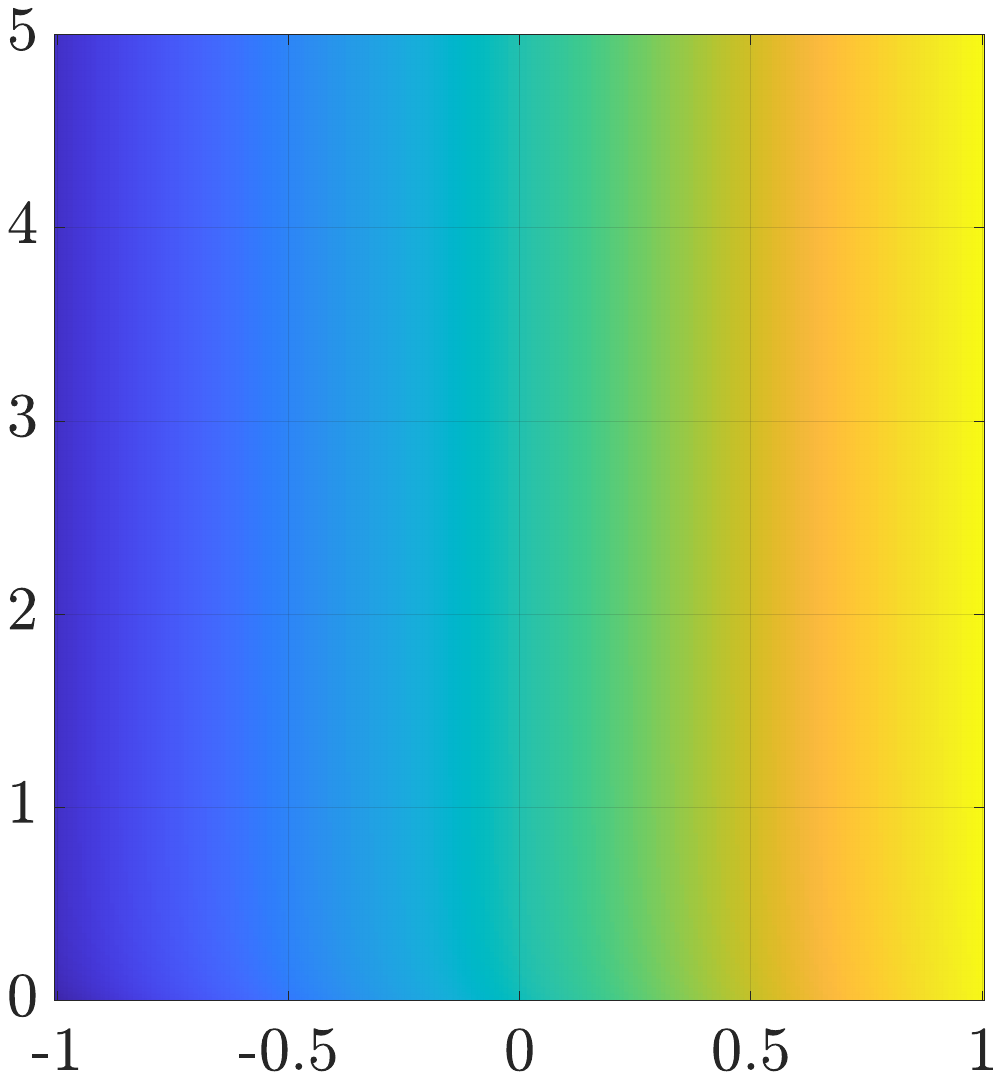} & \centering
\includegraphics[trim=1cm 3cm 1cm 4cm, clip, height=.275\textwidth]{ 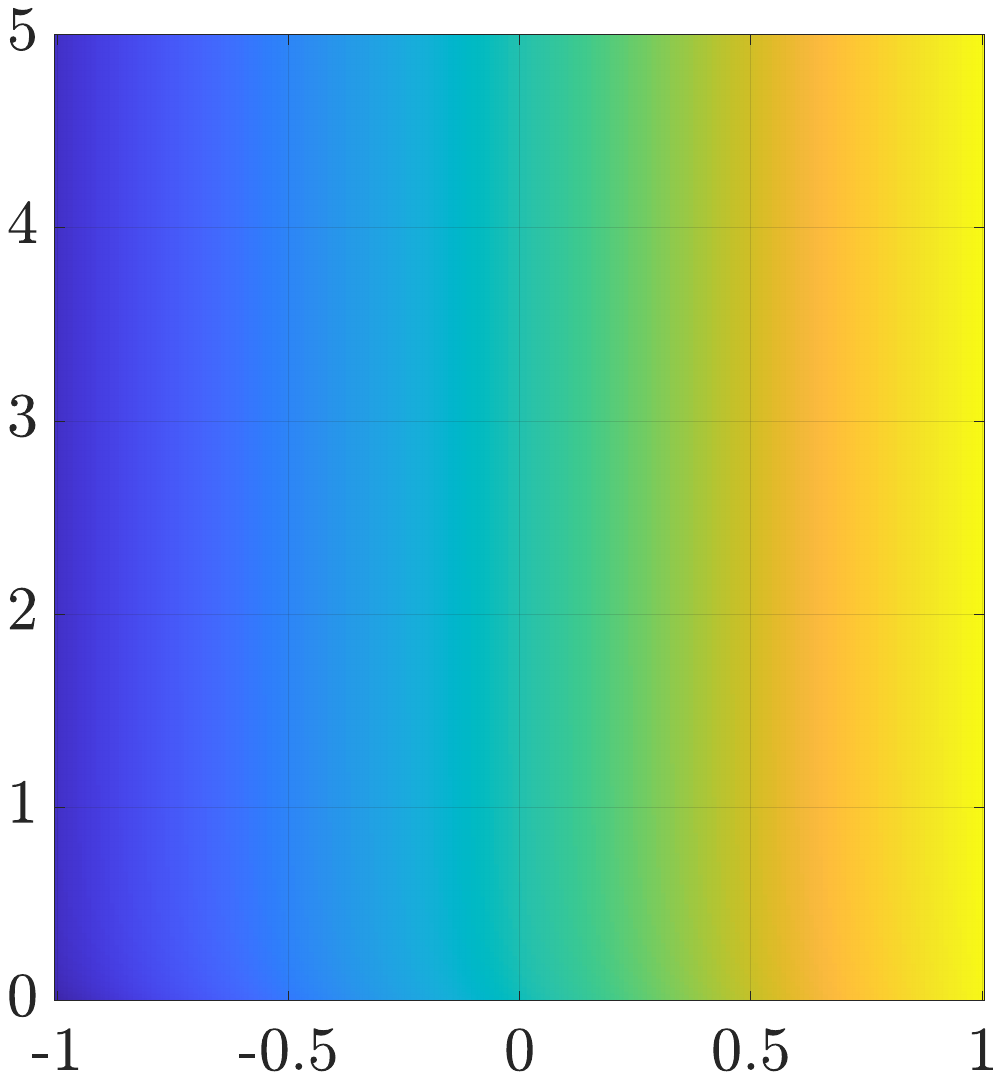} & 
\includegraphics[trim=1cm 3.5cm 1cm 4.5cm, clip, height=.275\textwidth]{ 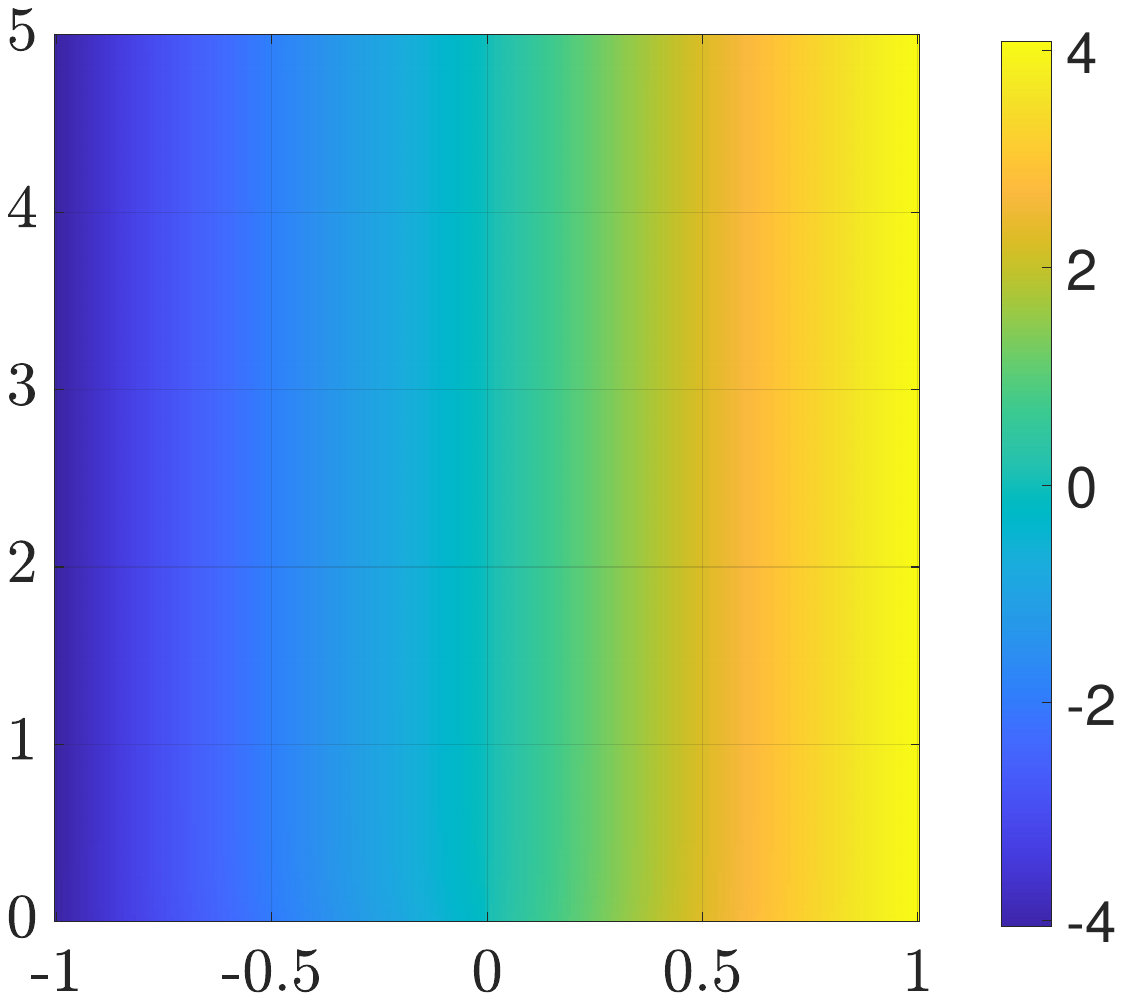} \\
\end{tabular}
\caption{\color{black}Comparison of the mean field optimal control problem solved via the proposed numerical procedure against the MPC routine relying on the first-order semi-Lagrangian scheme discussed in \cite{mfc_hierarchy}. For the three cases, we display the resulting total cost (top), controlled evolution (middle), and mean field control (bottom), with $u_{MPC}$ retaining the best performance.}\label{tab:mf_comparison}
\end{table}
\color{black}
\subsection{Test 2: Cucker-Smale model}
We proceed by testing the proposed methodology in a consensus control problem for non-linear and non-local dynamics governing the evolution of a $N$-agent system. In its reduced binary
semilinear formulation \eqref{tilde_din}, the model is written for the state $ \bs = (x,x_*,v,v_*)^{\top}  \in\R^{2d}\times\R^{2d}$,
with an interaction kernel given by
\begin{equation}
    P(x,x_*) = \dfrac{1}{1 + \|x-x_*\|^2}\,.
\end{equation}

The control variable $\bu = (u,u_*)\in\R^{2d}$ is here computed as
\begin{equation}\label{cscost}
    \Tilde{\bu} = \underset{u(\cdot)}{argmin}\; \dfrac 12\int\limits_0^{+\infty} \|v-\Bar{v}\|^2 + \|v_*-\Bar{v}\|^2 + \gamma \|\bu\|^2 dt\,,
\end{equation}
with $d=15$, $\gamma = 0.01$, and target velocity $\Bar{v} = \dfrac{v+v_*}{2}$. 

As discussed before, the key ingredient of the proposed methodology is a NN approximation of either the control $\bu_\theta$ or the (controlled) update $\bs'_\theta$ of the reduced binary problem. In both cases, the synthetic data have been generated from $N_s = 10^5$ uniform samples of interacting couples of agents $\{\bs_{(i)}\}_{i=1}^{N_s}$ within $\Omega\times\Omega = [-5,5]^{4d}$, together with their associated (sub)optimal DSDRE feedback control for the discrete-time infinite horizon OCP \eqref{u_sdre} starting from each $\bs_{(i)}$ for the sampled time-step $\Delta t_{(i)}\in[0,1]$.   We highlight that while the dynamics and controls are defined in the entire space, we consider a finite computational domain $\Omega$ for sampling.

Furthermore, we notice that we can write the state penalty term as
\begin{equation}
     \|\bv - \Bar{v}\|^2 = \|\bv-M\bv\|^2 
\end{equation}
for   $\bv=(v,v^*)^\top$ and   a suitable block matrix $M\in\R^{2d\times 2d}$ defined as 
\begin{equation} M = 
    \begin{bmatrix}
        \begin{matrix}
            \frac12&\frac12\\\frac12&\frac12
        \end{matrix} & & \bigzero\\
        &\ddots&\\
        \bigzero&&\begin{matrix}
            \frac12&\frac12\\\frac12&\frac12
        \end{matrix}
    \end{bmatrix}
\end{equation}
It follows that
\begin{equation}
\begin{aligned}
            \|\bv-M\bv\|^2 &= \langle \bv - M\bv ,  \bv - M\bv \rangle = \bv^{\top}\mathbb{I}_{2d}\bv + \bv^{\top}M^{\top}M\bv - 2\bv^{\top}M\bv\\
            &= \bv^{\top}\bigg(\mathbb{I}_{2d} + M^{\top}M - 2M\bigg)\bv\,.
\end{aligned}
\end{equation}
Thus, we can write the cost \eqref{cscost} in quadratic form \eqref{ocp} w.r.t. linear operators 
\begin{equation}\label{cost_operators}
    Q = \mathbb{I}_{2d} + M^{\top}M - 2M\,,\qquad R = \dfrac{\gamma}{2}\cdot \mathbb{I}_{2d}\,.
\end{equation}

Once the problem has been written in semi-linear form, we rely on the DSDRE approach (described in algorithm \ref{alg1}) for the generation of a dataset collecting $N_s$ samples of coupled states, time-steps, associated feedback laws, and controlled state updates.

We test the following architectures:\\

    $u_\theta^{FNN}:$ $K = 1$ hidden layers with $100$ neurons, and  $\sigma(\bx) =log(exp(\bx) + 1)$\\[-10pt]
    
    $u_\theta^{RNN}:$ $K = 1$ LSTM cell with $100$ neurons, $\sigma(\bx) = max(0,\bx)$, $\rho(\bx) = (1 + exp(-\bx))^{-1}$\\[-10pt]
    
    $\mathbf{s'}_{v,\theta}^{FNN}:$  $K = 3$ hidden layers with $100$ neurons per layer, and $\sigma(\bx) = max(0,\bx)$\\[-10pt]
    
    $\mathbf{s'}_{v,\theta}^{RNN}:$ $K=1$ LSTM cell with $100$ neurons, $\sigma(\bx) = max(0,\bx)$,  $\rho(\bx) = (1 + exp(-\bx))^{-1}$\\

\begin{figure}
    \centering
    \includegraphics[height = 0.21\textheight]{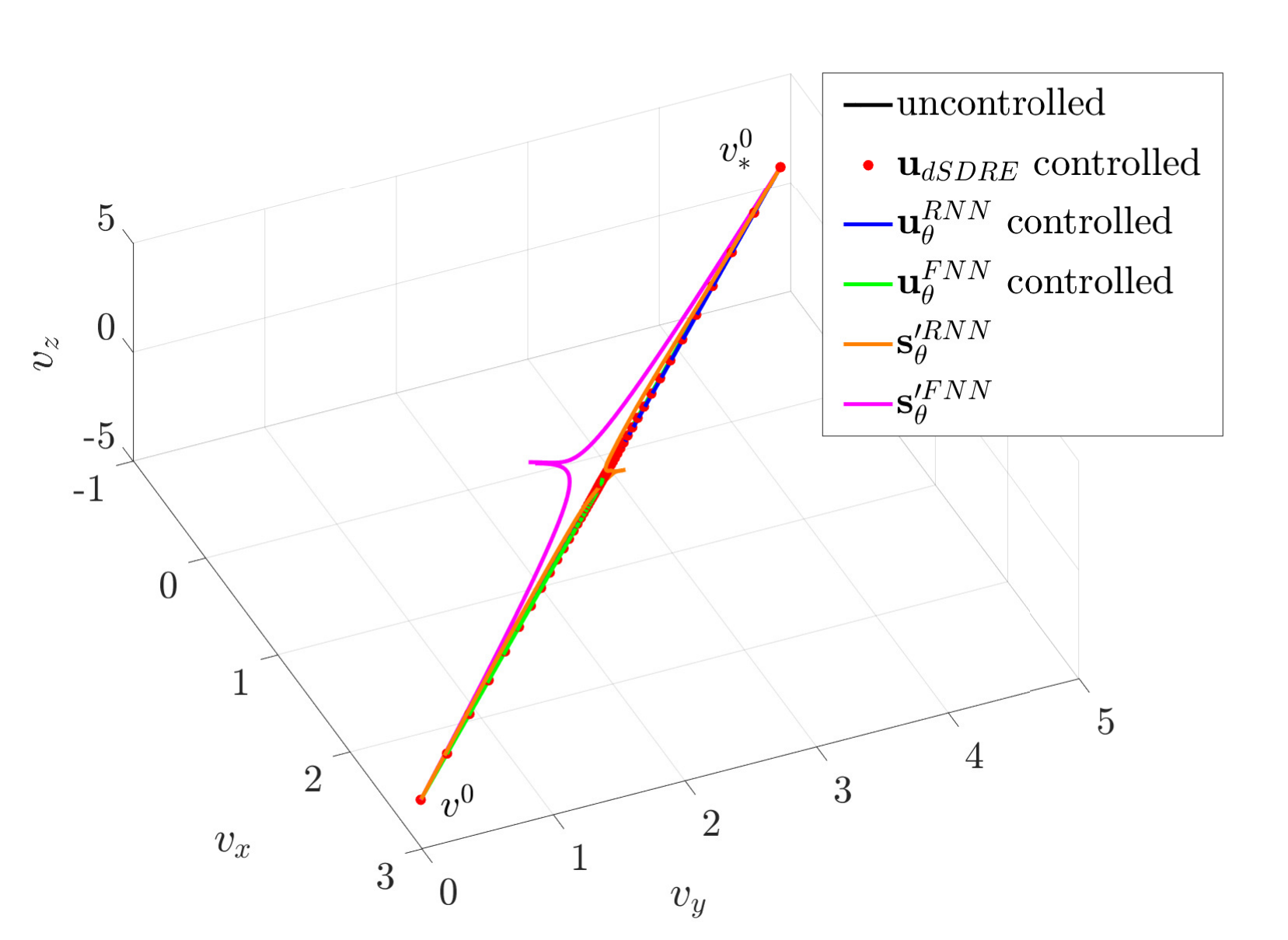}
    \includegraphics[height = 0.21\textheight]{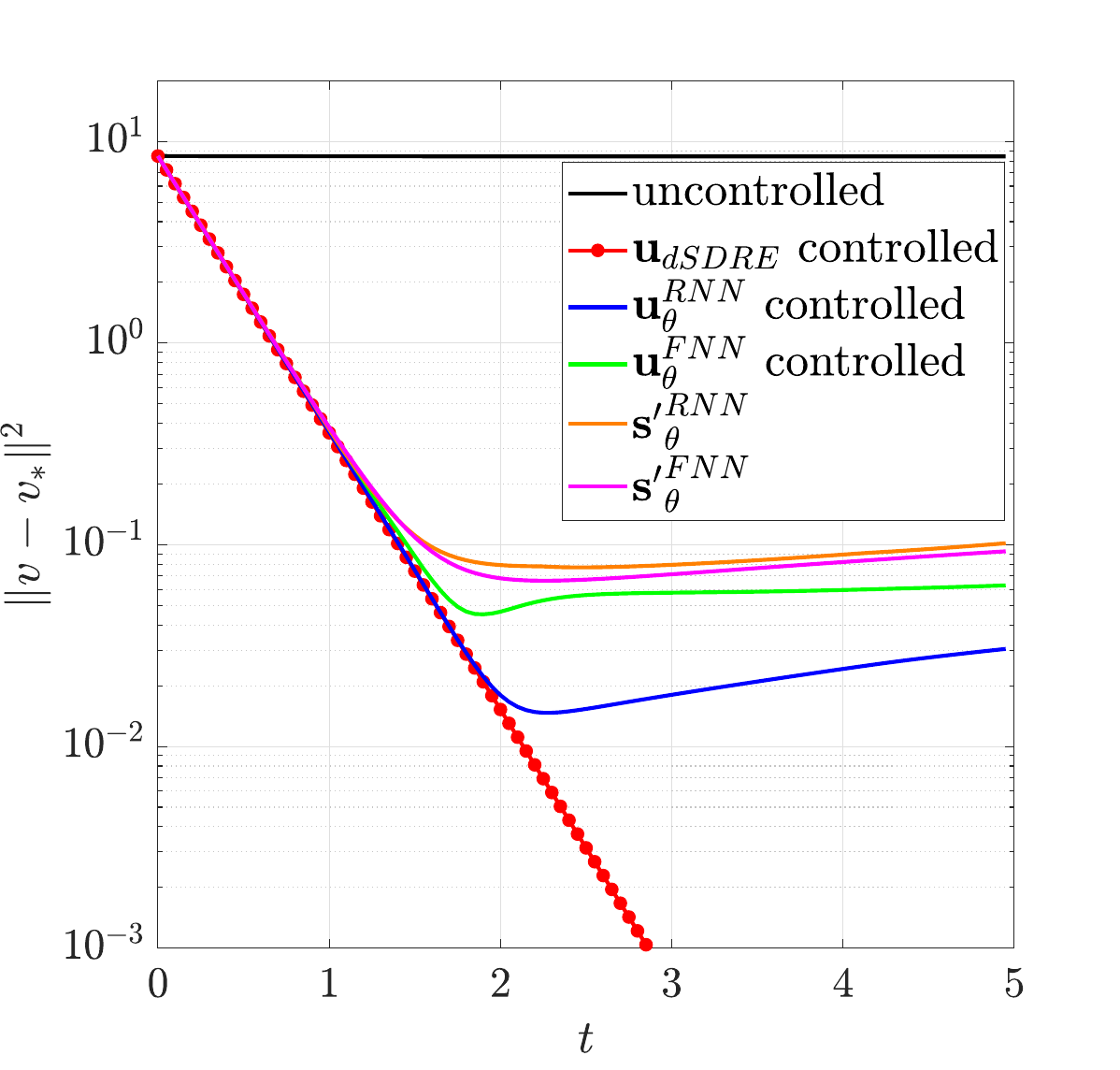}
    \includegraphics[height = 0.21\textheight]{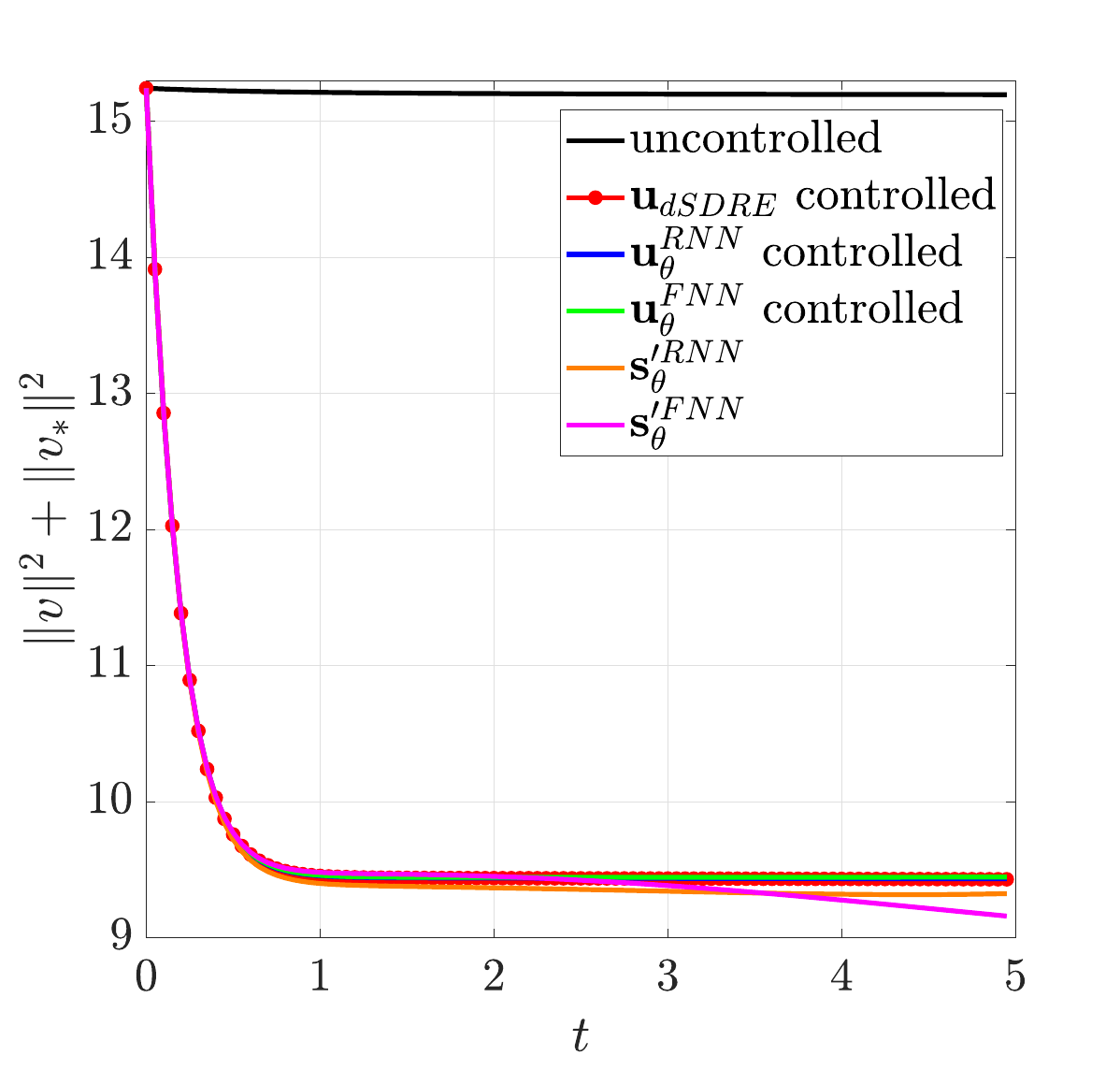}
    \caption{Evolution of the two interacting agents' in the velocity space restricted to the first $3$ dimensions (left); convergence to consensus in logarithmic scale (centre); consensus configurations for the different approximators (right).} 
    \label{fig:cs_binary}
\end{figure}

The behaviour of the approximated controlled binary dynamics for a single couple of agents   is   shown in Figure \ref{fig:cs_binary}, where we compare the true DSDRE solution with our NN approximations. As the trajectories evolve, the approximated states deviate from the DSDRE closed-loop. However, this does not significantly affect the MC simulation, as a sampled pair of agents will only interact for the duration of a single time step.
The goodness of fit of the trained models measured in a \emph{test set} $\mathcal{T}_t \nsubseteq \mathcal{T}\cup\mathcal{T}_v$ is presented in Table \ref{tab:goodness}. 

\begin{figure}
    \centering
    \begin{subfigure}[b]{0.35\textwidth}
    \includegraphics[width =\textwidth]{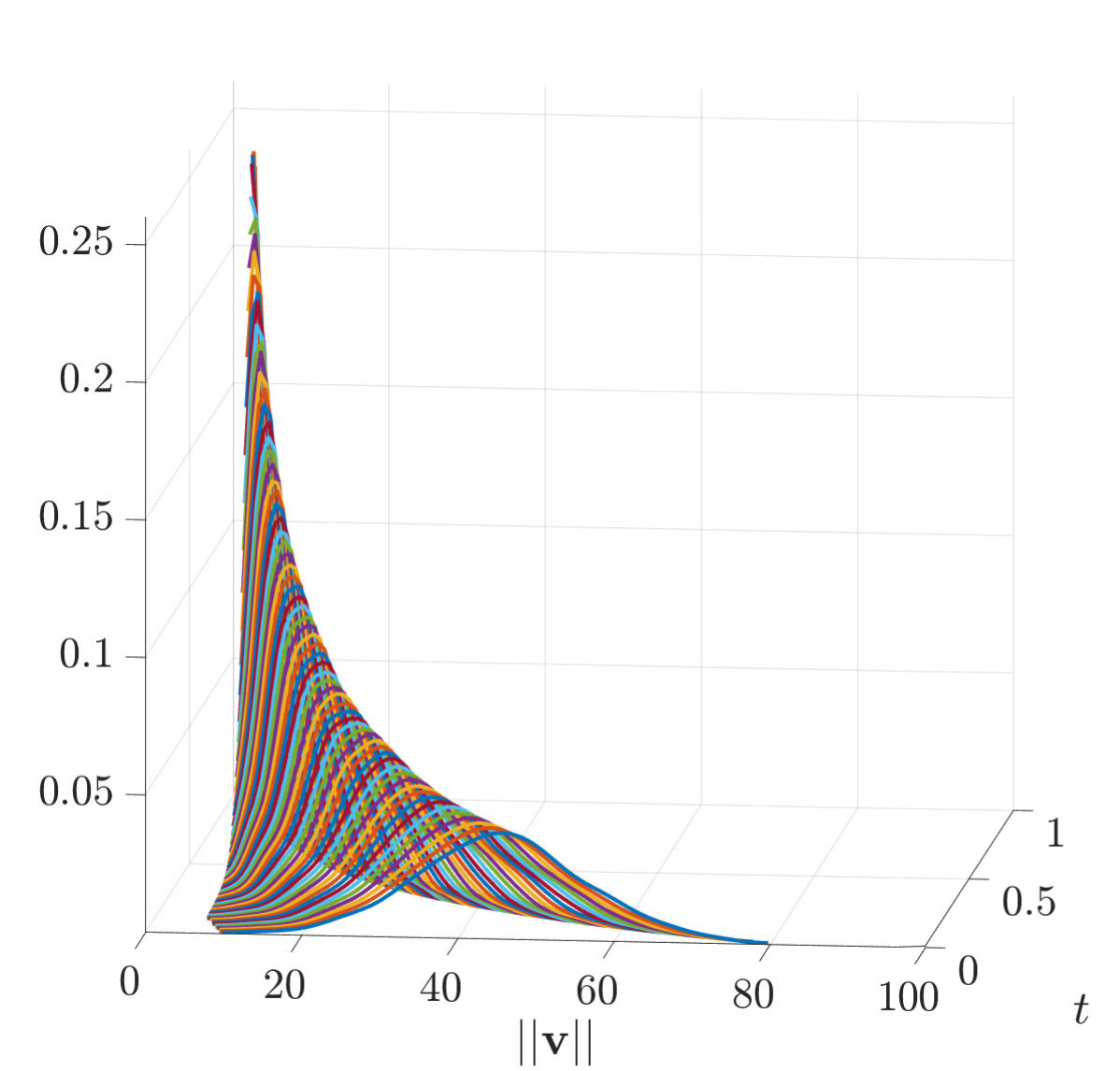}
    \caption{$\bu^{RNN}_\theta$}
    \end{subfigure}
     \begin{subfigure}[b]{0.35\textwidth}
    \includegraphics[width =\textwidth]{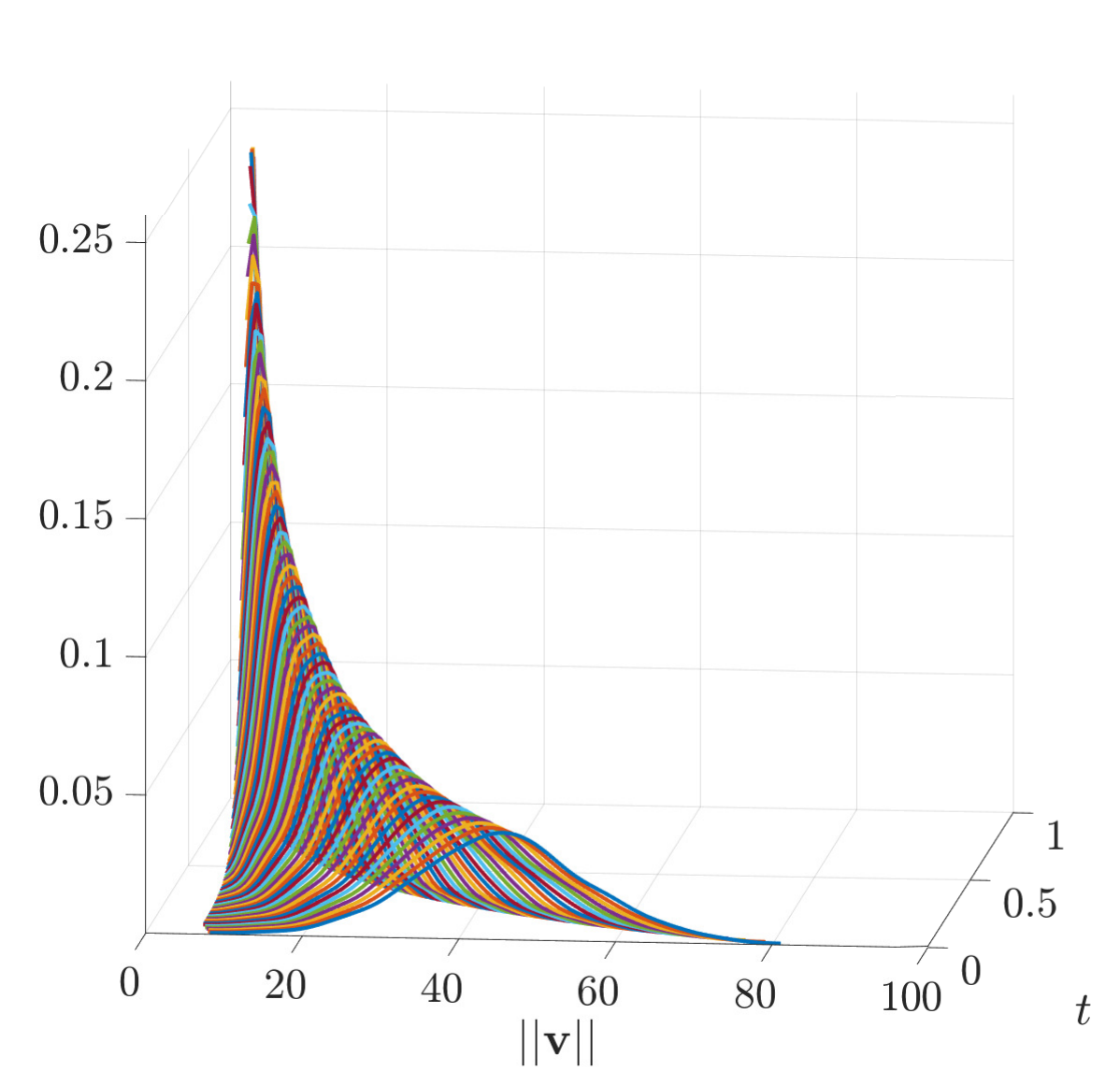}
    \caption{$\bu^{FNN}_\theta$}
    \end{subfigure}\\
     \begin{subfigure}[b]{0.35\textwidth}
    \includegraphics[width =\textwidth]{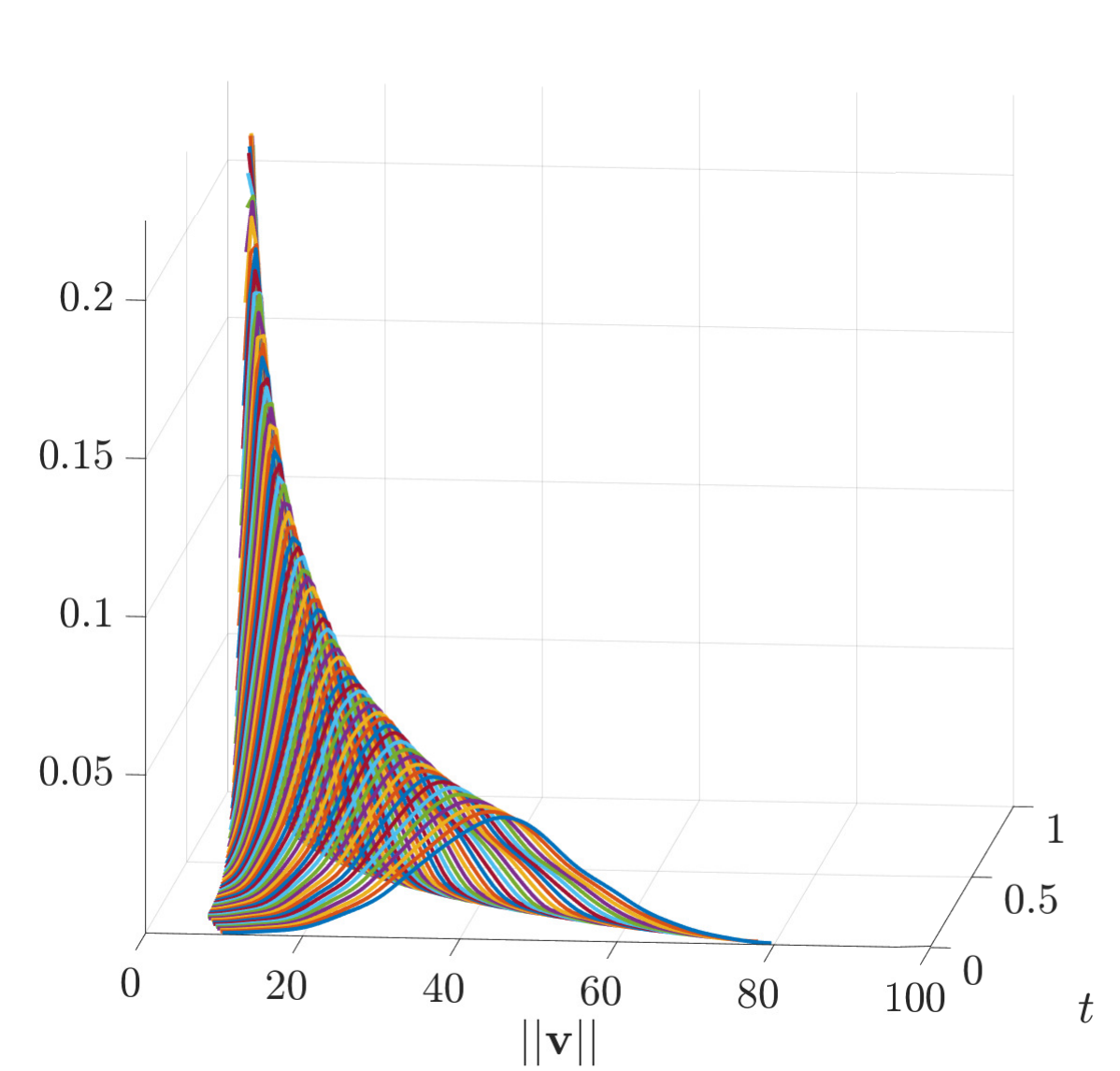}
    \caption{${\bs'}^{RNN}_\theta$}
    \end{subfigure}
     \begin{subfigure}[b]{0.35\textwidth}
    \includegraphics[width =\textwidth]{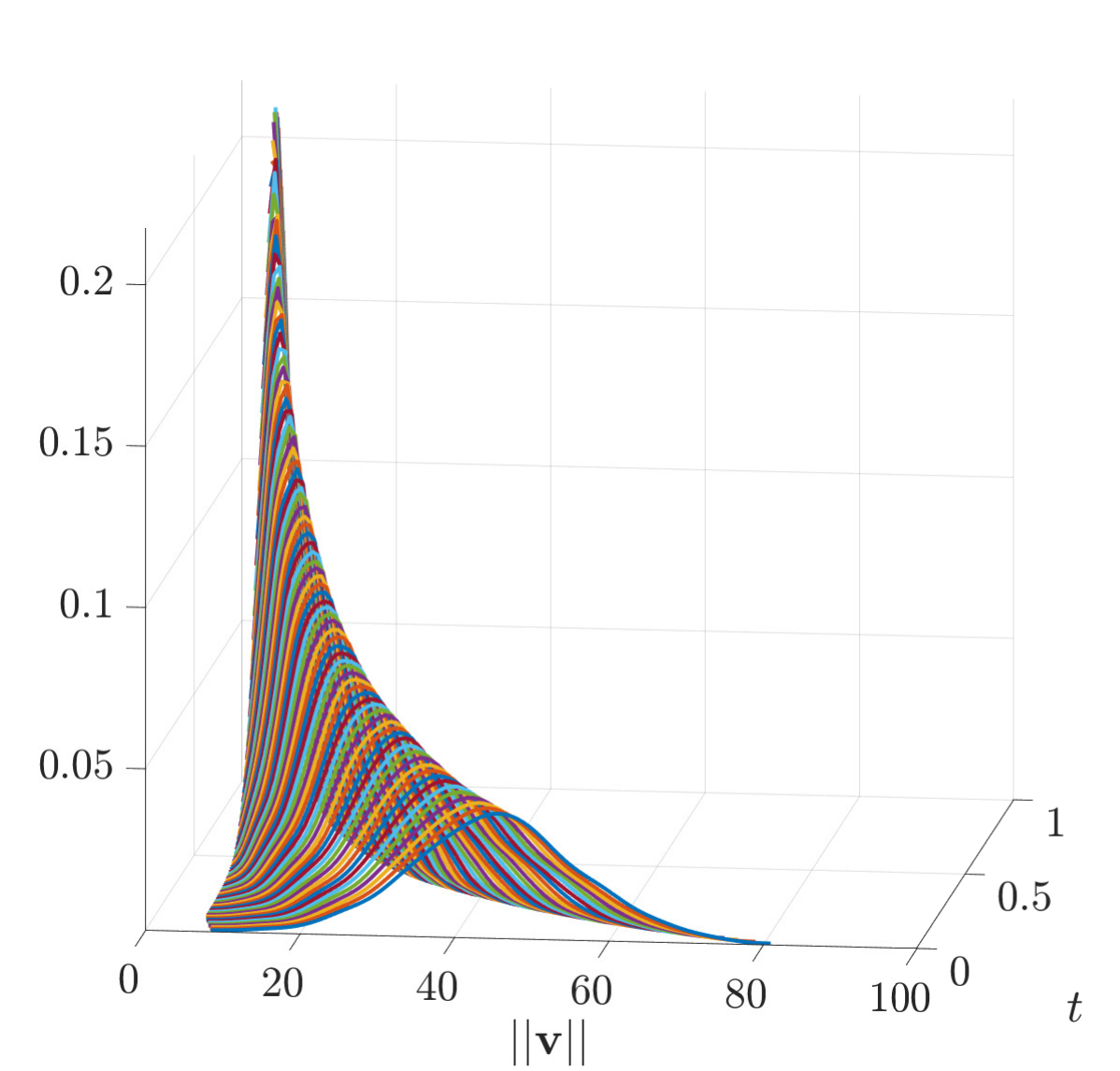}
    \caption{${\bs'}^{FNN}_\theta$}
    \end{subfigure}
    \begin{subfigure}[b]{0.35\textwidth}
    \includegraphics[width =\textwidth]{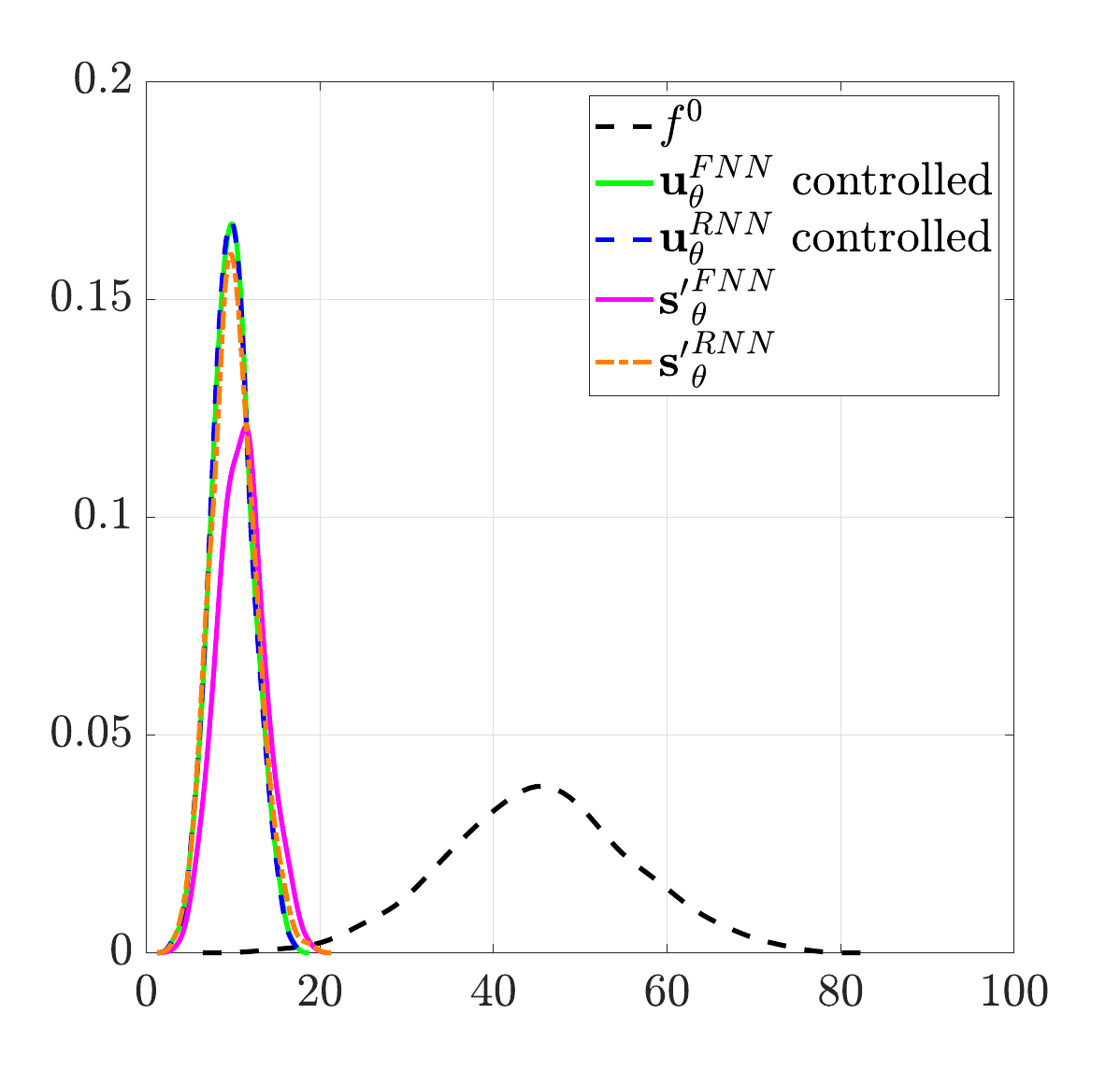}
    \caption{$t_{50} = 0.5s$}
    \end{subfigure}
    \begin{subfigure}[b]{0.35\textwidth}
    \includegraphics[width =\textwidth]{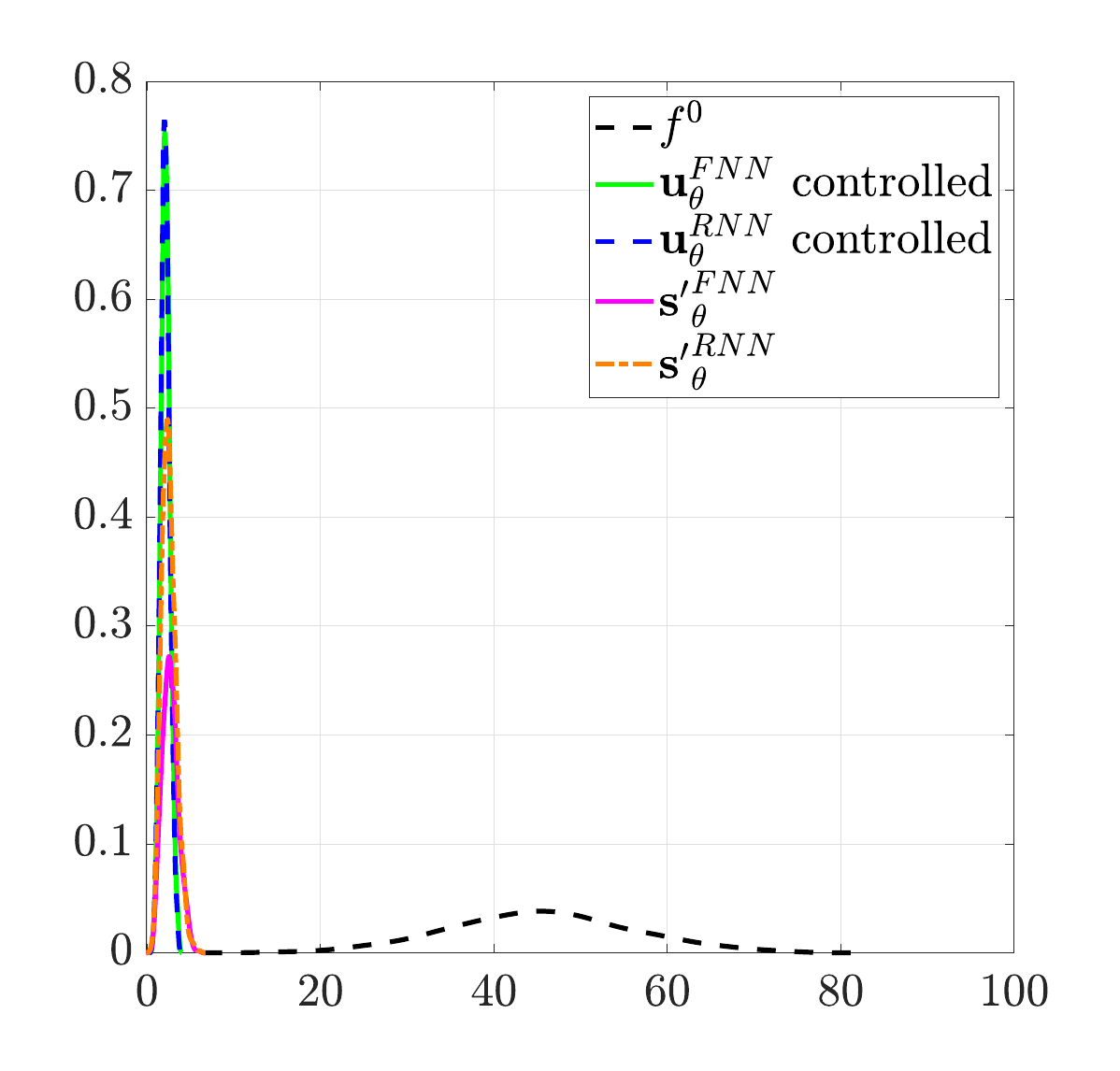}
    \caption{$t_{100} = 1s$}
    \end{subfigure}
    \caption{Density of agents' velocities in norm for Test 1. Comparison of the time evolution of the MC pdf obtained considering different approximation approaches for the controlled binary interactions. In the bottom row, we display a comparison of the different approximated densities at two different discrete times.}
    \label{fig:cs_simulation}
\end{figure}

In Figure \ref{fig:cs_simulation}, we show the evolution of the density of agents' velocities according to a Monte Carlo simulation of binary controlled dynamics, as described in Algorithm \ref{alg2}. All the approximation models equally succeed in steering the density distribution to concentration profiles. Snapshots of the mean field distributions at specific time intervals show how similar the consensus configurations are across the trained models. 
\begin{table}[ht]
    \centering
    \begin{tabular}{l c c c c }
        $d = 15$  &  $N_s= 10^2$ & $N_s=10^3$ &  $N_s=10^4$ & $N_s=10^5$\\[2pt]
        \hline\\[-7pt]
        {$\;{\bs'}^{FNN}_\theta$} &  {$0.045515$} &  {$0.227243$} &  {$2.097497$}&  {$24.060892$}\\
        {$\;{\bs'}^{RNN}_\theta$} &  {$0.067627$} &  {$0.283442$} &  {$3.210856$}&  {$34.087541$}\\
        {$\;\bu_\theta^{FNN}$} &  {$0.293578$} &  {$2.447205$} &  {$21.754862$} &  {$2.2594\times10^2$}\\
        {$\;\bu_\theta^{RNN}$} &  {$0.340493$} &  {$2.225319$} &  {$22.486736$}&  {$2.2559\times10^2$}\\
        {$\;\bu$} &  {$2.3866\times10^2$} &  {$2.3738\times10^3$} &  {$-$} & $-$\\
    \end{tabular}
    \centering
    \caption{CPU times (seconds) for coupled agents in $\R^{4d},\,d=15$, when considering different number of MC samples. The omitted records exceeded a time threshold $t_{max} = 24h$.}
    \label{tab:cpu_ns}
\end{table}

So far, the reliance on approximants for the binary controlled state update has been motivated in terms of efficiency. In Tables \ref{tab:cpu_ns} and \ref{tab:cpu_d}, we present CPU times\footnote{The experiments have been executed in MATLAB R2022a installed on a machine with Intel Core i7-10700 processor running at 2.90GHz.} (in seconds) for the MC simulation with the different models, compared with the DSDRE solution $\bu$.   The improved efficiency achieved through approximation is crucial, as the alignment between the mean field dynamics and their kinetic approximation relies on high-frequency sampling ($\varepsilon\ll1$).     The evolution of the agents' distribution is approximated along a sequence of discrete times $t_n= n \Delta t$, $\Delta t = 0.05$, $n=1,\dots,100$. The tables address the computational cost associated to sampling a number of controlled binary interactions in the MC simulation, and the dimensionality $d$ of the agents physical space, for a total of $4d$ dimensions. The reliance on NN approximation models allows for a speedup of $2$ to $3$ orders of magnitude. For $d\geq7$ and $N_s\geq 10^4$, the computational cost resulting from the use of $\mathbf{u}_\theta$ displays a linear growth, whilst $\mathbf{s'}_\theta$ performs even better. CPU times for the true DSDRE controlled dynamics exceed $24$ hours. 

\begin{table}[ht]
    \centering
    \begin{tabular}{l c c c c c}
        $N_s=10^4$  &  $d=3$ & $d=7$ &  $d=10$ & $d=15$ & $d=30$\\[2pt]
        \hline\\[-7pt]
        {${\bs'}^{FNN}_\theta$} & {$1.048226$} & {$1.212633$} & {$1.390813$} & {$2.142041$} & {$2.617840$}\\
         {${\bs'}^{RNN}_\theta$} & {$2.033726$} & {$2.243084$} & {$2.493256$} & {$3.210856$} & {$3.893368$}\\
          {$\bu_\theta^{FNN}$} & {$7.712628$} & {$11.006977$} & {$15.041731$} & {$21.754862$} & $70.172311$\\
          {$\bu_\theta^{RNN}$} & {$7.734224$} & {$11.224325$} & {$15.991421$} & {$22.486736$} & {$70.564372$}\\
          {$\bu$} & {$1.1979\times10^3$} & {$5.2136\times10^3$} & {$-$} & {$-$} & {$-$}\\[7pt]
    \end{tabular}
    \centering
       \caption{CPU times (seconds)  when considering $10^4$ MC samples of coupled agents in $\R^{4d}$, with varying $d$. The omitted records exceeded a time threshold $t_{max} = 24h$.}
    \label{tab:cpu_d}
\end{table}

\begin{table}[t]
	\centering 
	\begin{tabular}{c c c c c c c }
		& \multicolumn{3}{c}{Test 1}&\multicolumn{3}{c}{Test 2}\\
		\cmidrule(lr){2-4}\cmidrule(lr){5-7}\morecmidrules\cmidrule(lr){2-4}\cmidrule(lr){5-7}
		\emph{Model} & $r^2$ & \emph{MSE} & \emph{MRE $\%$}& $r^2$ & \emph{MSE} & \emph{MRE $\%$} \\ 
		\cmidrule(lr){1-1}\cmidrule(lr){2-4}\cmidrule(lr){5-7}\\[-0.7ex]
		{$\;{\bs'}^{FNN}_\theta$} &  {$0.99998$} &  {$ 0.075252$} &  {$0.38964$} &
        {$  0.9994 $} & {$ 0.038274 $} &  {$ 0.60308$}\\
		{$\;{\bs'}^{RNN}_\theta$} &  {$0.99999$} &  {$ 0.0069192$} &  {$0.3739$} &
		{$ 0.9998 $} & {$ 0.012784$} &  {$ 0.20956 $}\\
        {$\;\bu_\theta^{FNN}$} & {$0.99996$} &  {$0.045596$} &  {$0.63555$} &  {$ 0.9997 $} & {$ 7.6882 $} &  {$ 3.4018 $}\\
        {$\;\bu_\theta^{RNN}$} &  {$0.99998$} &  {$0.018018$} &  {$0.40136$} &  {$ 0.99979 $} & {$ 5.5492 $} &  {$ 2.7483 $}
	\end{tabular}
	\caption{Goodness of fit for Tests 1 and 2 in terms of: \emph{coefficient of determination} $r^2$, \emph{mean squared error} and \emph{mean percentage error}. $\mathcal{T}_t$ is a collection of sampled states $\bs^{(i)}\in(\Omega\times\Omega)$ and related target DSDRE values for $i = 1,...,10^5$.}\label{tab:goodness}
\end{table}

\subsection{Test 3: quasi-Morse potential}
As a second numerical test, we consider the consensus control problem analysed in \cite{controlling_swarming} for the interacting particles of a second order system of agents in the physical space $\R^3$. The velocities are here governed by both a self-propulsion force,  expressed in the $i-th$ agent by the term $(\alpha - \beta \|v_{i}(t)\|^2)v_{i}(t)$, for fixed $\alpha\geq0,\;\beta>0$, and an attraction-repulsion force acting though the pairwise interaction potential $W$. We consider a radial potential of the form
\begin{equation}\label{potential}
    W(x) =  V\big(\|x\|\big) - C\,V\bigg(\dfrac{\|x\|}{l}\bigg),\qquad V(r) := -exp\bigg\{\dfrac{-r^p}{p}\bigg\}
\end{equation}
with $C = 0.6$, $p = 1.5$, $l=0.5$, $\alpha = 2$, $\beta=1.5$, in a similar configuration as in \cite{controlling_swarming}.
Accordingly, for the couple of interacting agents $\bs=(x,x_*,v,v_{*})$ we define
\begin{equation}
                P_x(x,x_{*}) = \|x-x_{*}\|^{(p-2)}\bigg(\dfrac{C}{l^p}\,e^{-\dfrac{\|x-x_{*}\|^p}{p\,l^p}}-e^{-\dfrac{\|x-x_{*}\|^p}{p}} \bigg)
\end{equation}
\begin{equation}
                P_v(v) =  (\alpha - \beta \|v(t)\|^2)\,,
\end{equation}
with which we write the reduced binary dynamics in semilinear form as
\begin{equation}
    A = \begin{bmatrix}
    \mathbb{0}_{2d} & \mathbb{I}_{2d}\\
    A_x &  A_v
    \end{bmatrix}\qquad \bigg[A_x\bigg]_{i,j} = \begin{cases}
    - P_x(x,x_{*}) & \text{if } i=j\\
    P_x(x,x_{*}) & \text{if $m$ odd, } j = i+1\\
    P_x(x,x_{*}) & \text{if $m$ even, } j = i-1\\
    0 & \text{otherwise}
    \end{cases}
\end{equation}
where $A_v$ is the diagonal matrices with diagonal vector the component-wise application of $P_v(\cdot)$ to the vector with $d$ repetitions of $v,v_*$. Moreover, the cost operators for the consensus goal  $Q, R$ hold as in \eqref{cost_operators}.

\begin{figure}
    \centering
    \includegraphics[height = 0.21\textheight]{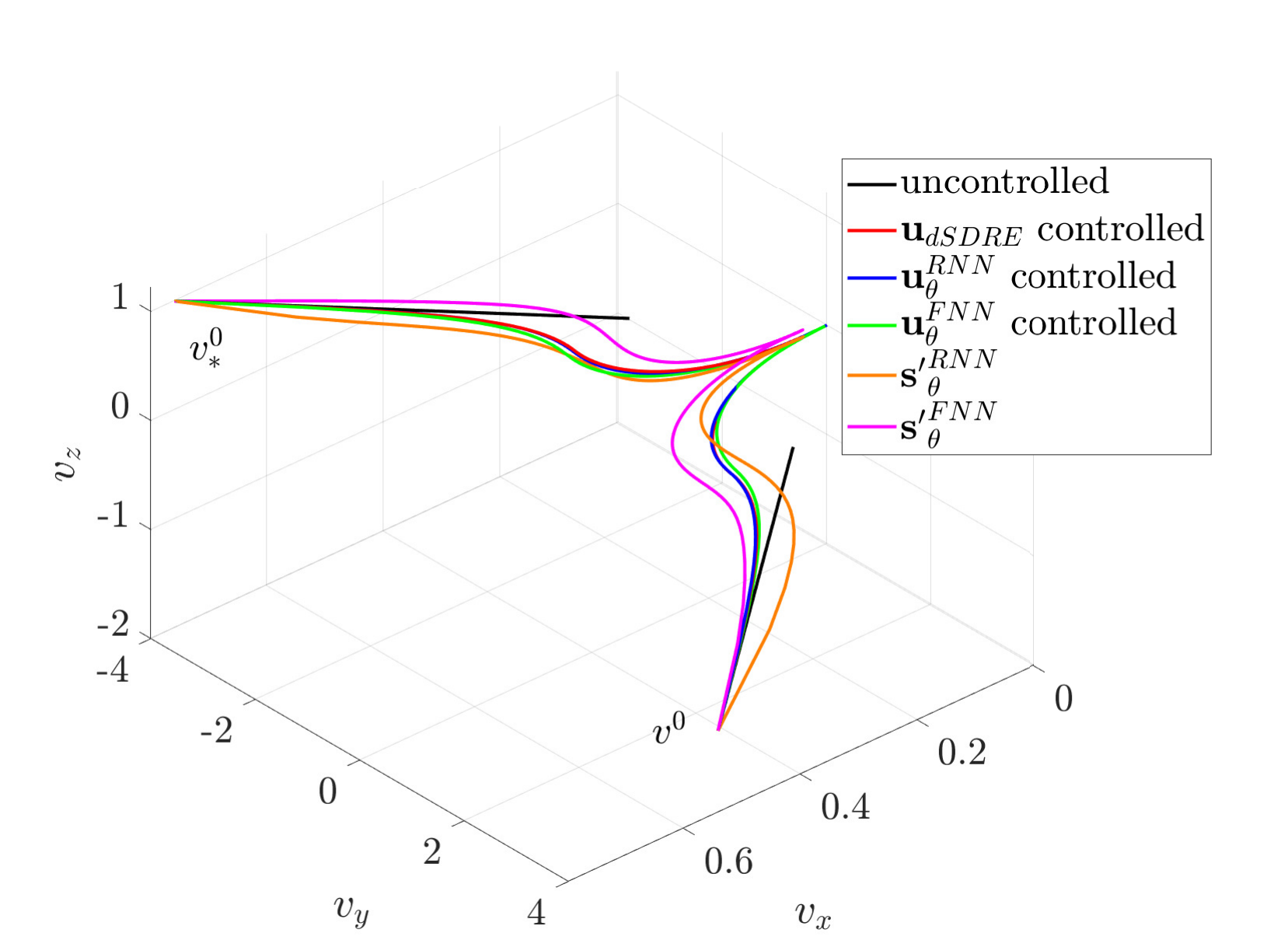}
    \includegraphics[height = 0.2\textheight]{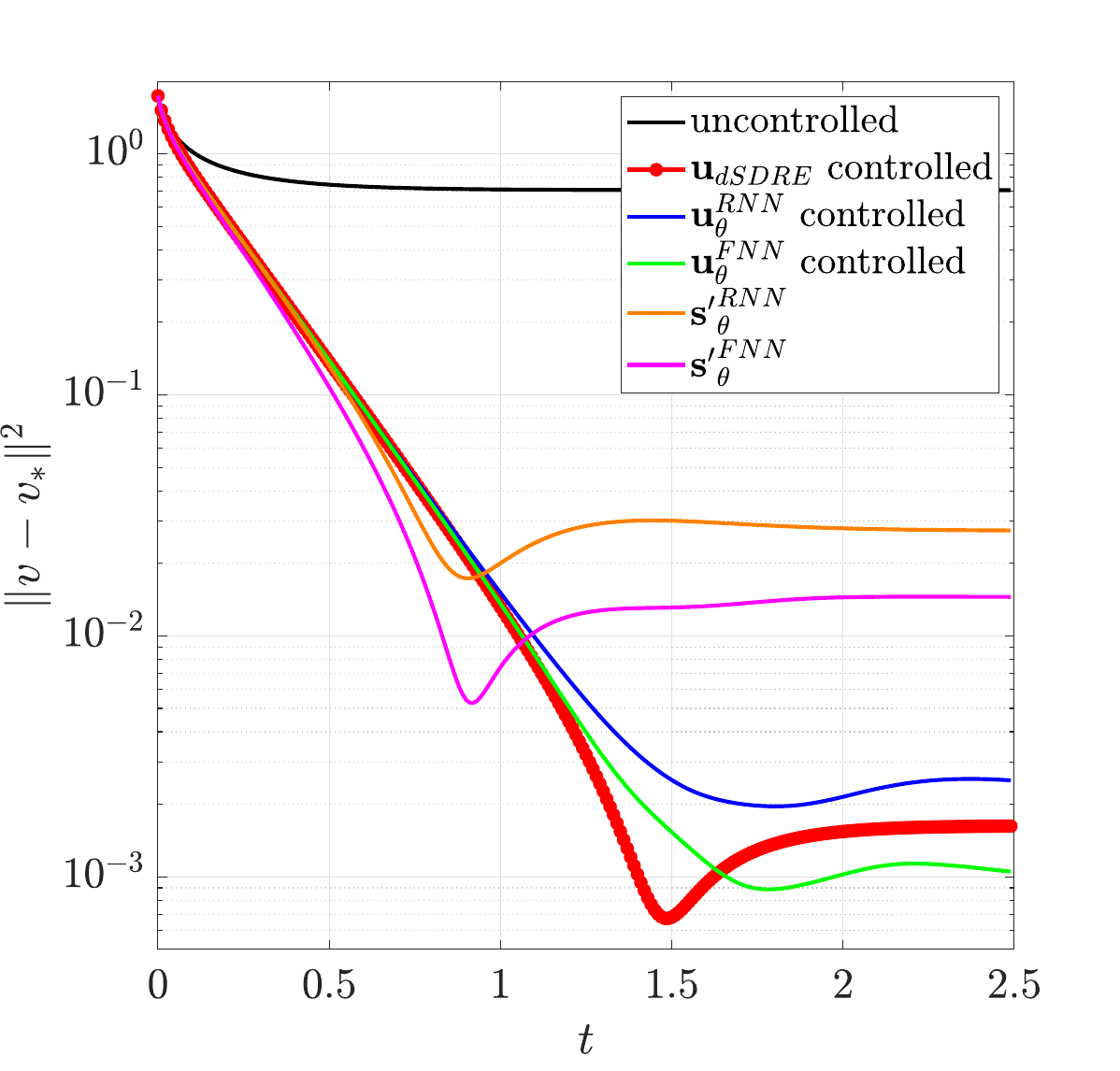}
    \includegraphics[height = 0.2\textheight]{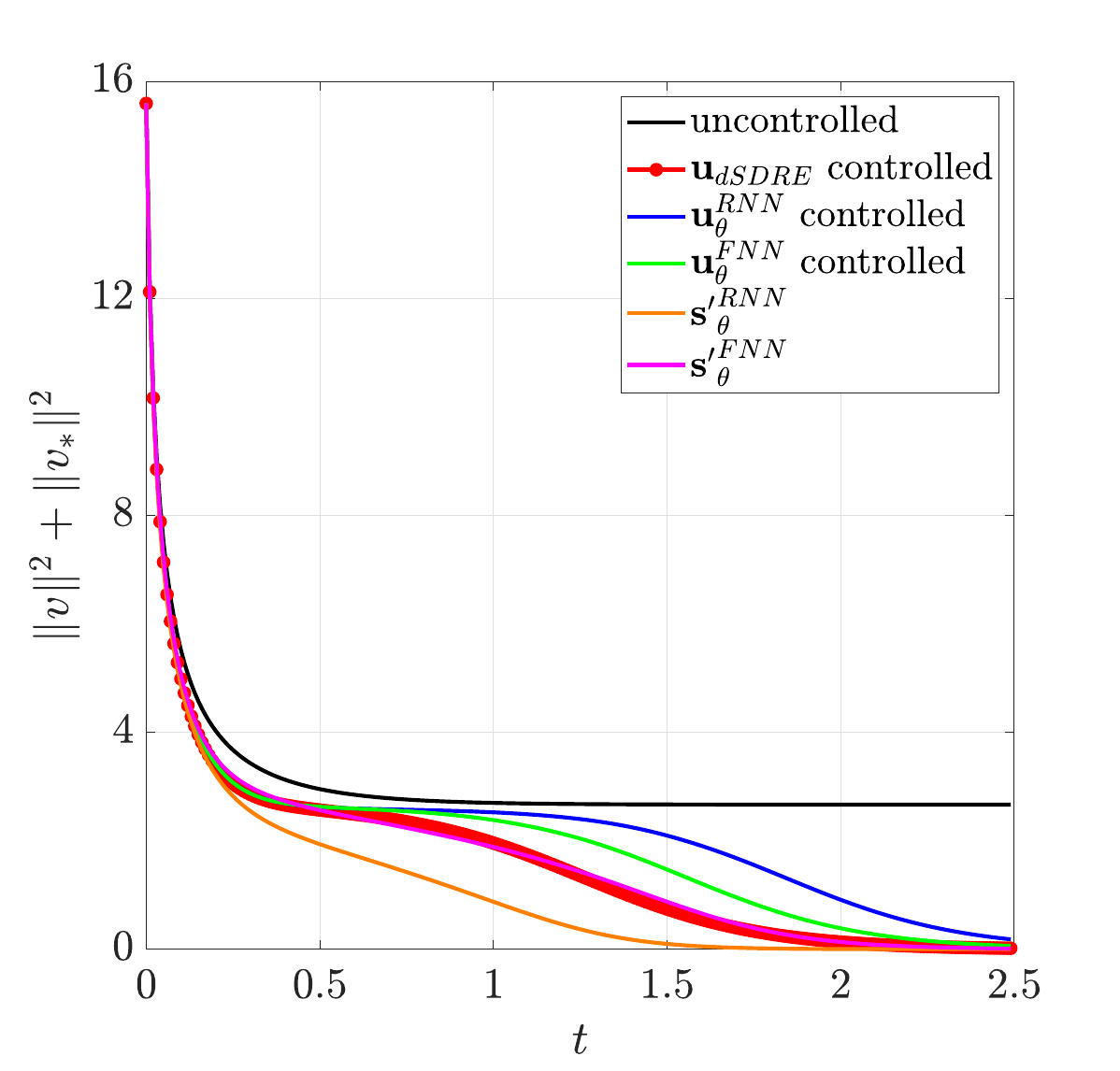}
    \caption{Evolution (in velocities) of a couple of approximately controlled interacting agents under the different models, versus the uncontrolled dynamics and the dSDRE-controlled ones (left); convergence to consensus in logarithmic scale (centre); consensus configurations for the different approximators (right).} 
    \label{fig:morse_binary}
\end{figure}

As discussed for the previous numerical example, we aim at approximating the binary control $\bu$, resulting from the DSDRE approach, and the related controlled state update $\bs'$. We train the following models:\\

    $\bu_\theta^{FNN}:$ $K = 3$ hidden layers ($100$ neurons per layer), $\sigma(\bx) =  log( 1 + e^\bx ) $\\[-10pt]
    
    $\bu_\theta^{RNN}:$ $K = 5$ ($100$ neurons per layer) with $l_1$ LSTM cell, $\sigma(\bx) = \rho(\bx) = log( 1 + e^\bx )$\\[-10pt]
    
    ${\bs'}_{v,\theta}^{FNN}:$  $K = 4$ ($100$ neurons per layer), $\sigma(\bx) = log( 1 + e^\bx )$\\[-10pt]
    
    ${\bs'}_{v,\theta}^{RNN}:$ $K = 3$ ($100$ neurons per layer), $l_1$ LSTM cell,  $\sigma_1(\bx) = \rho(\bx) = log( 1 + e^\bx )$,  
    $$\sigma_2(\bx) = \sigma_3(\bx) =elu(\bx)=\begin{cases}
    e^\bx -1 & \text{if } \bx\leq0\\
    \bx & \text{if } \bx>0
    \end{cases}$$

\begin{figure}[ht]
    \centering
    \includegraphics[height = 0.18\textheight]{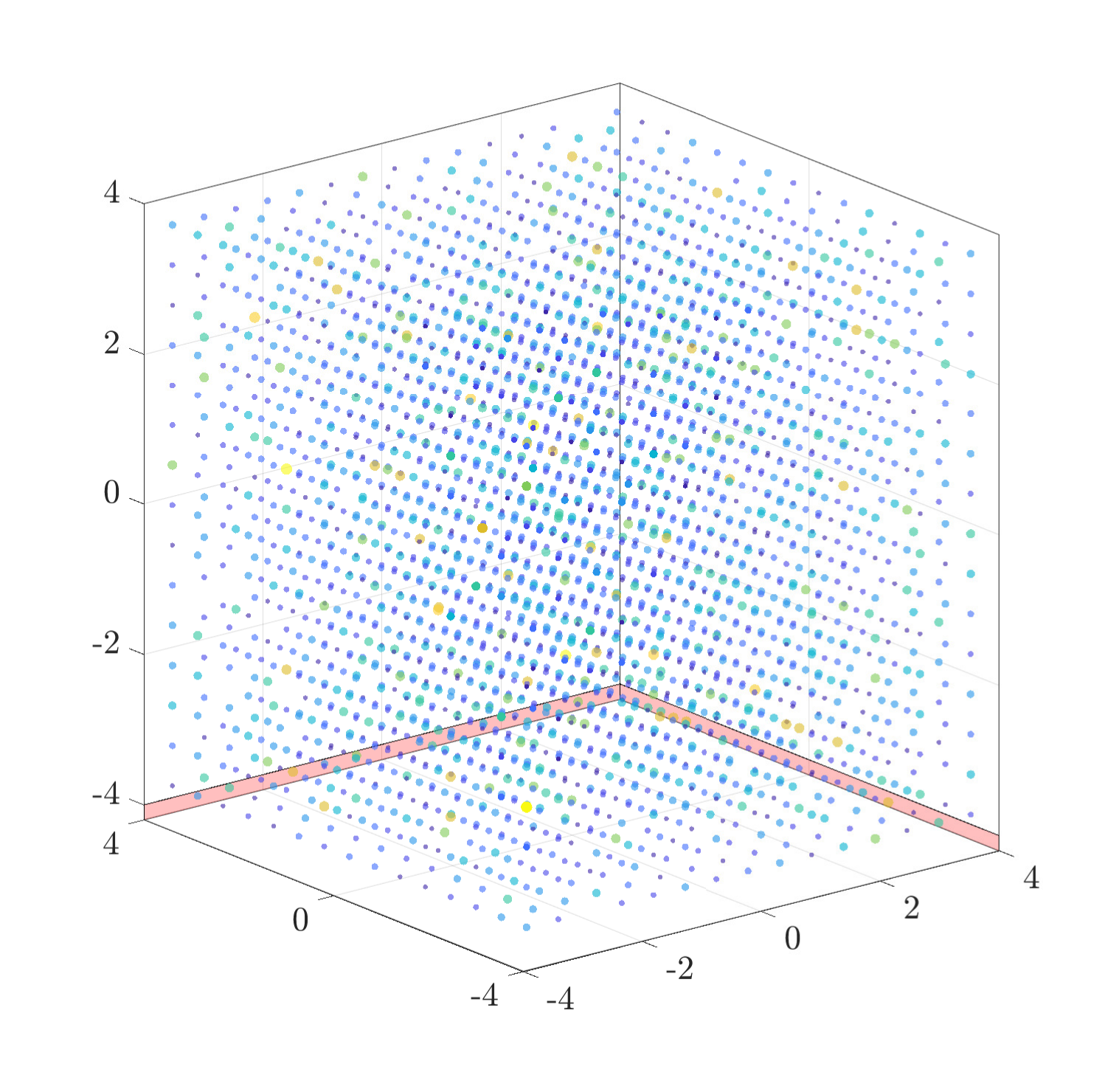}
    \includegraphics[height = 0.18\textheight]{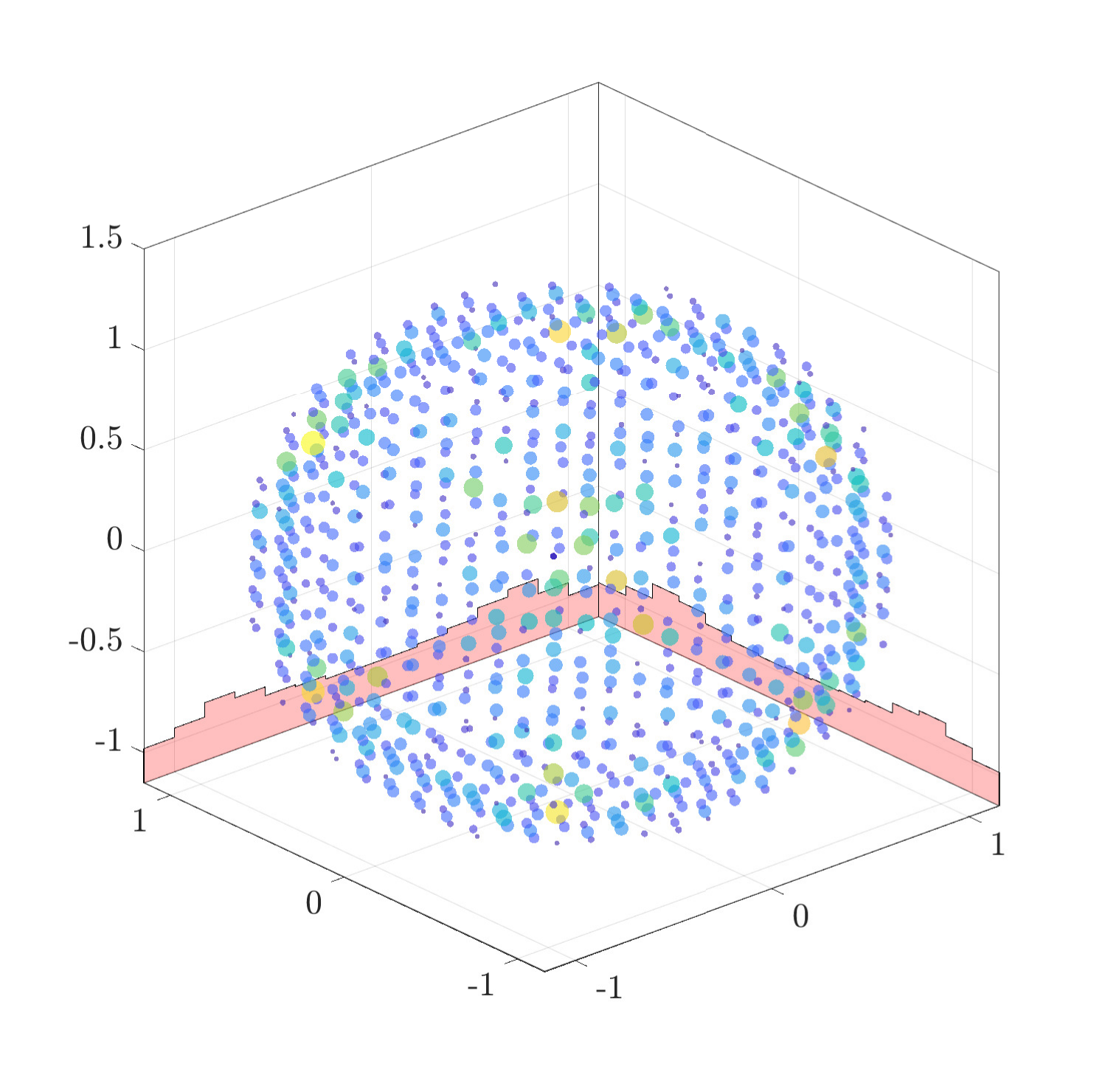}
    \includegraphics[height = 0.18\textheight]{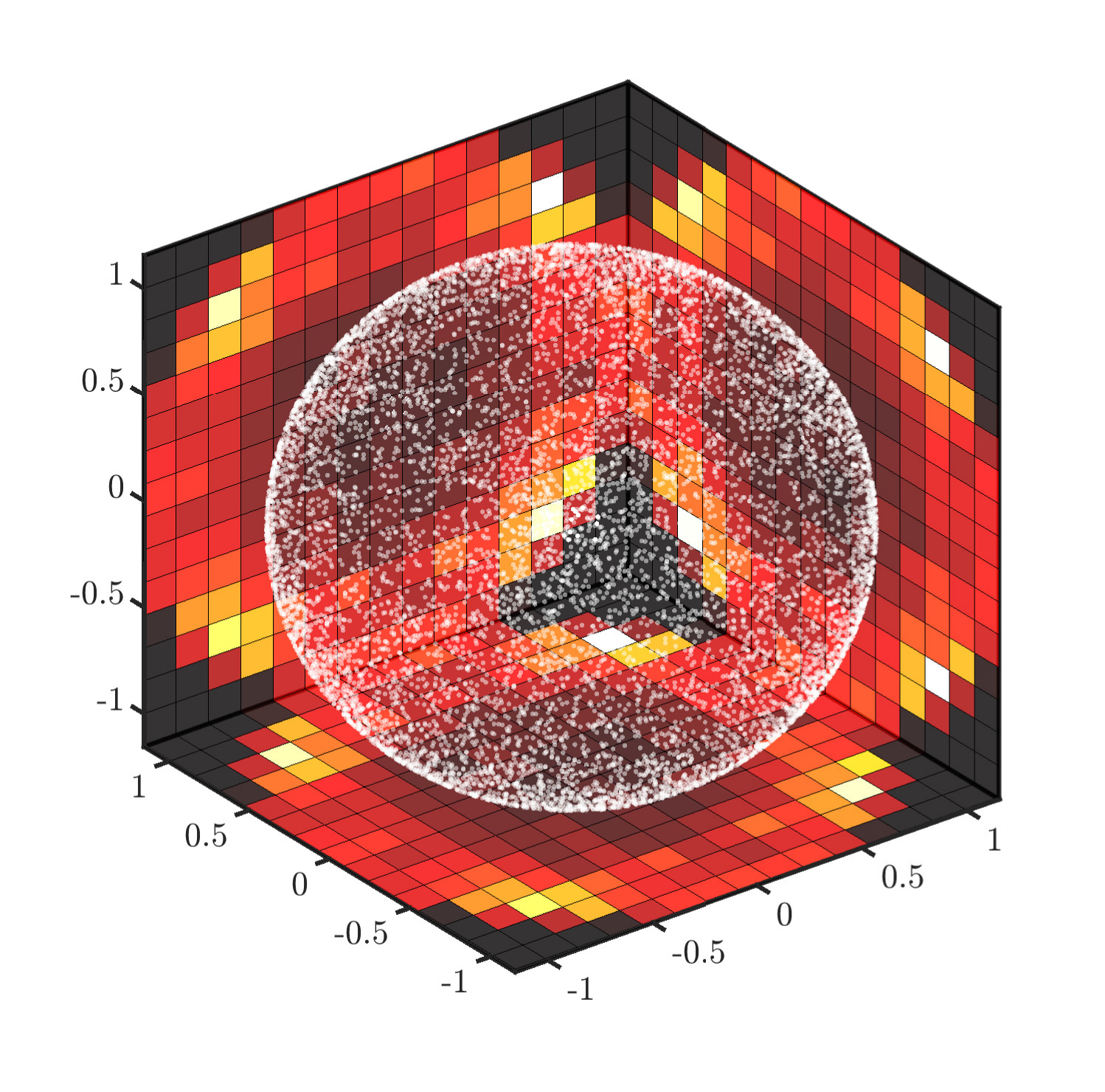}\\
    \begin{subfigure}[b]{0.25\textwidth}
    \includegraphics[width =\textwidth]{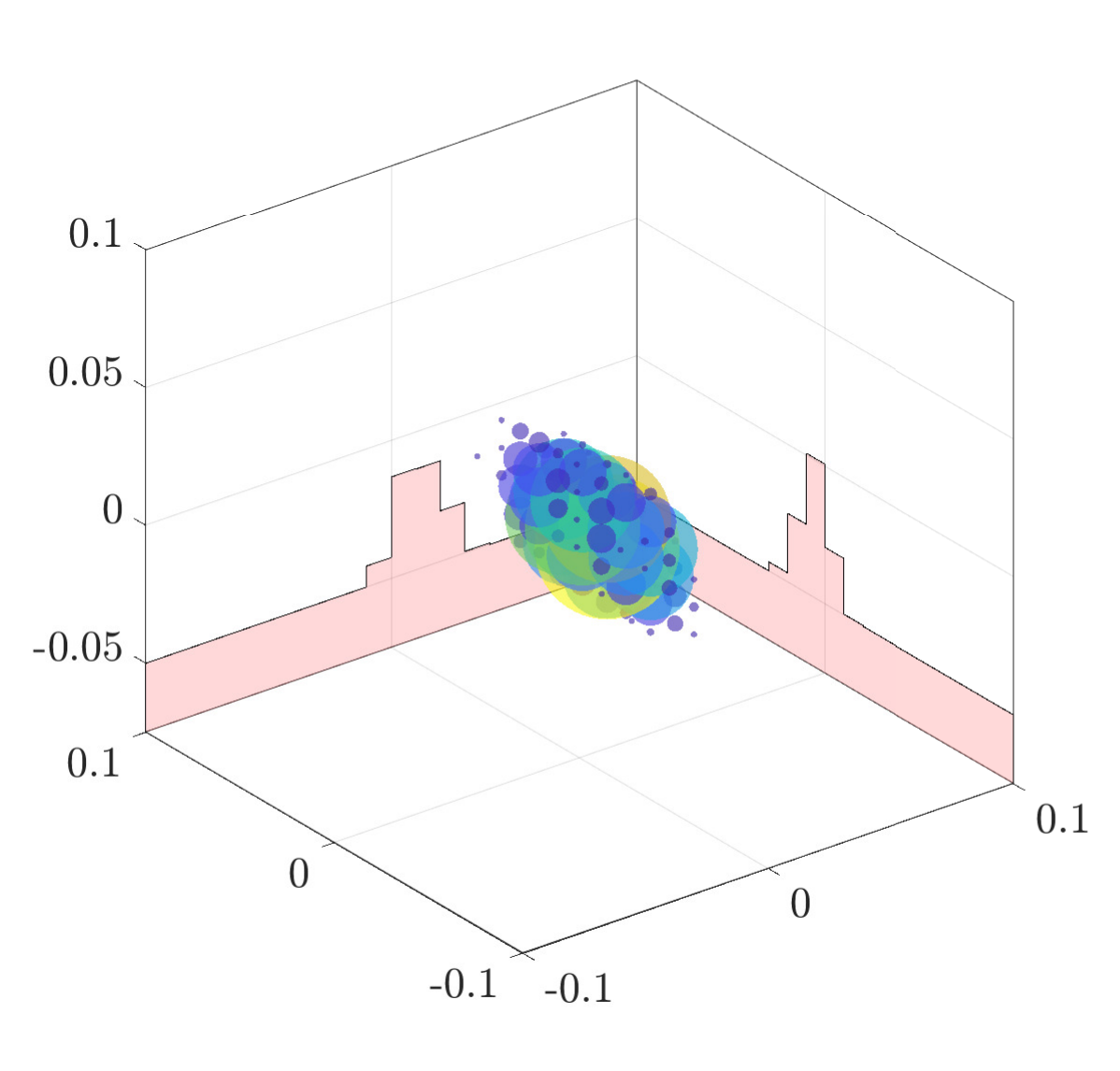}\\
    \includegraphics[width =\textwidth]{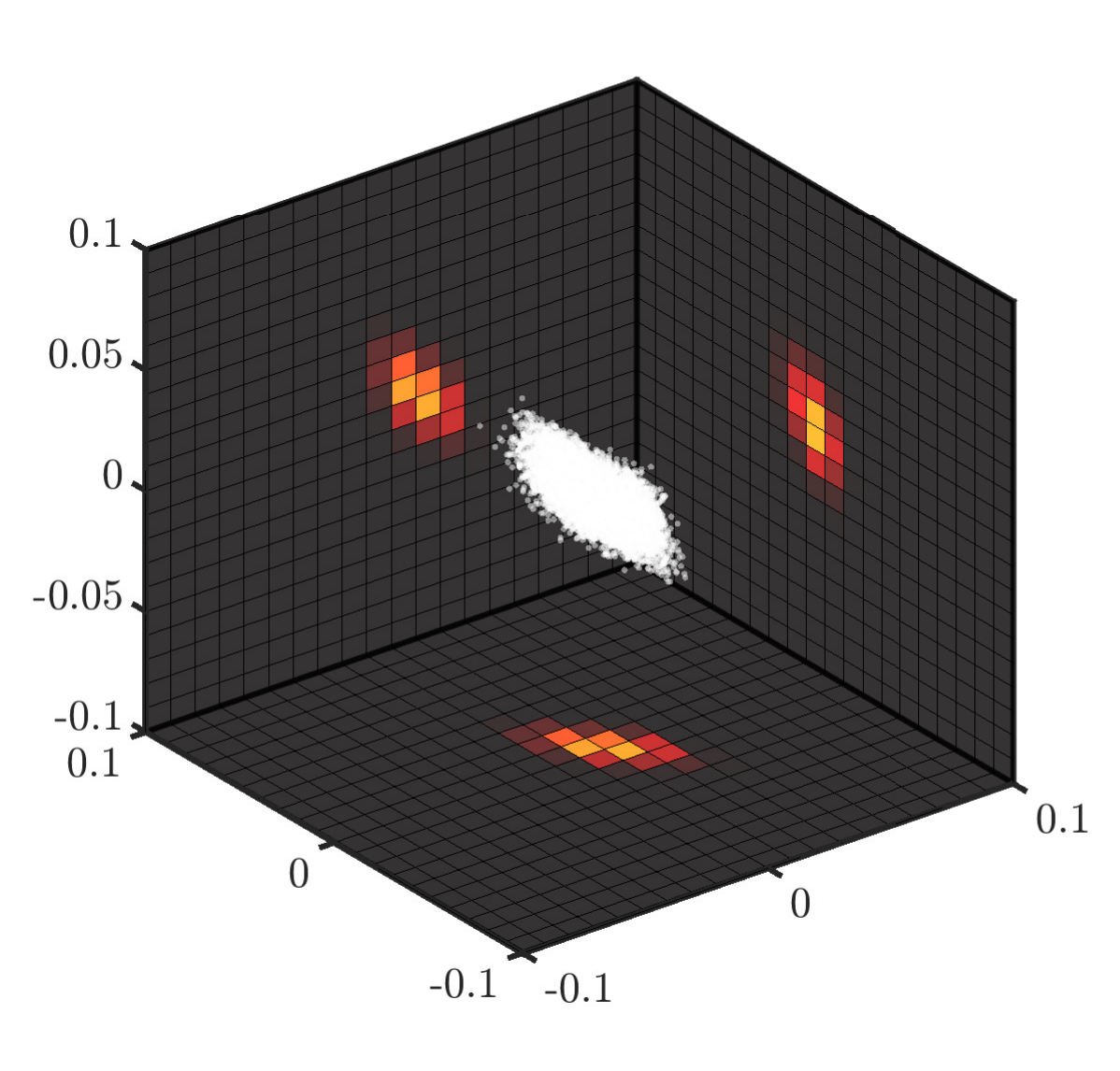}
    \caption{${\bs'}^{RNN}_\theta$}
    \end{subfigure}\hspace{-0.3cm}
    \begin{subfigure}[b]{0.25\textwidth}
    \includegraphics[width =\textwidth]{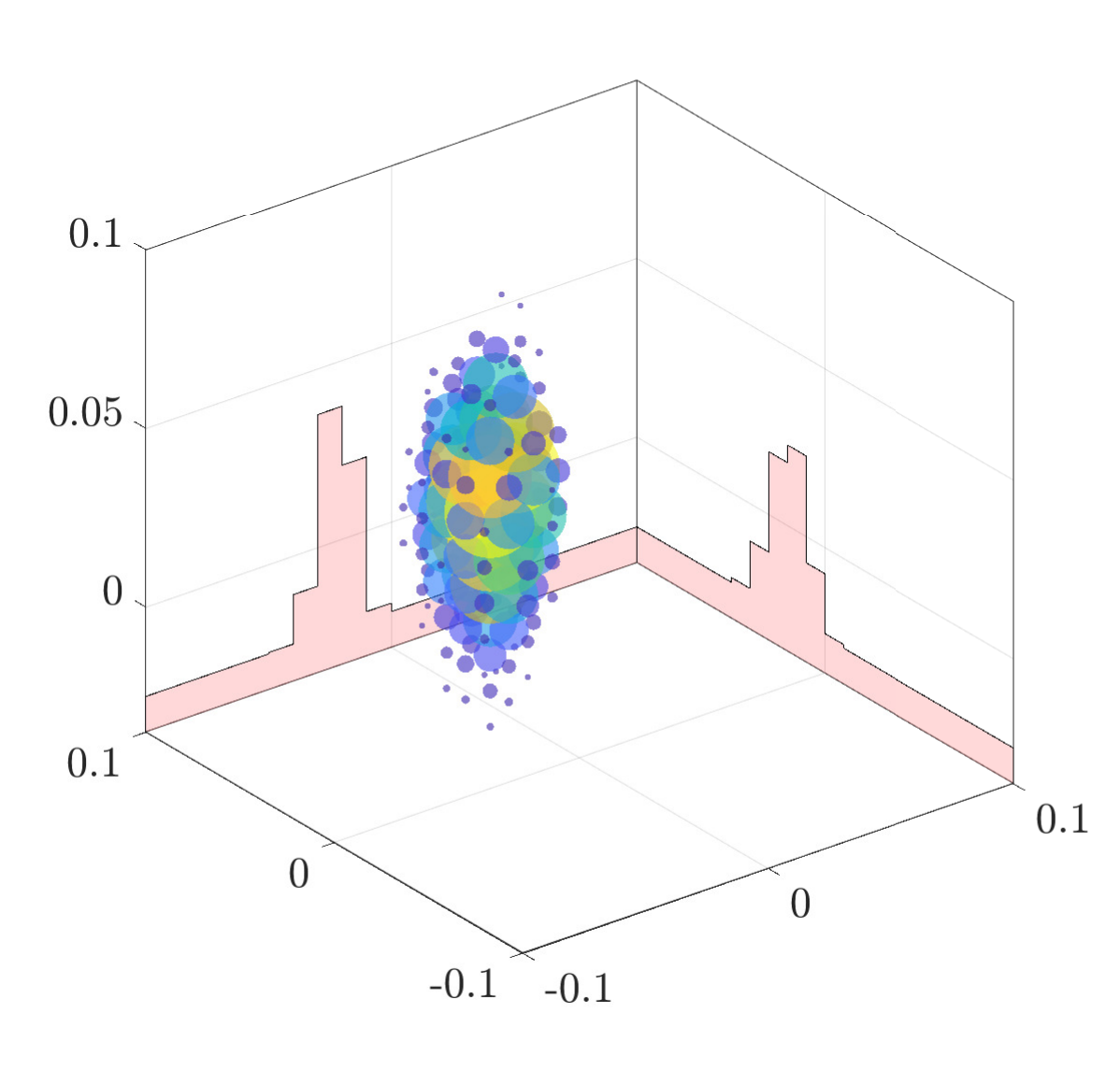}\\
    \includegraphics[width =\textwidth]{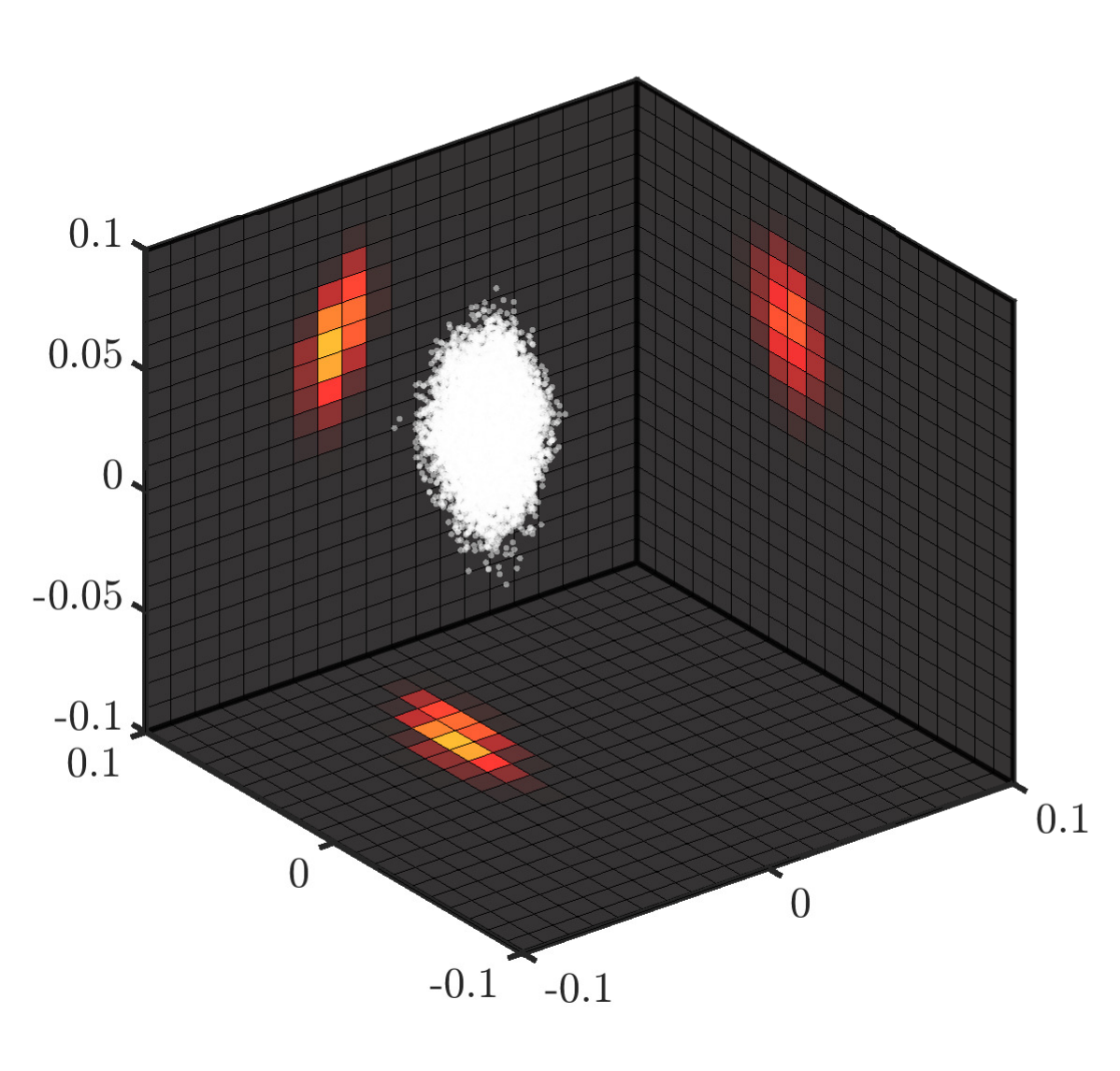}
    \caption{${\bs'}^{FNN}_\theta$}
    \end{subfigure}\hspace{-0.3cm}
    \begin{subfigure}[b]{0.25\textwidth}
    \includegraphics[width =\textwidth]{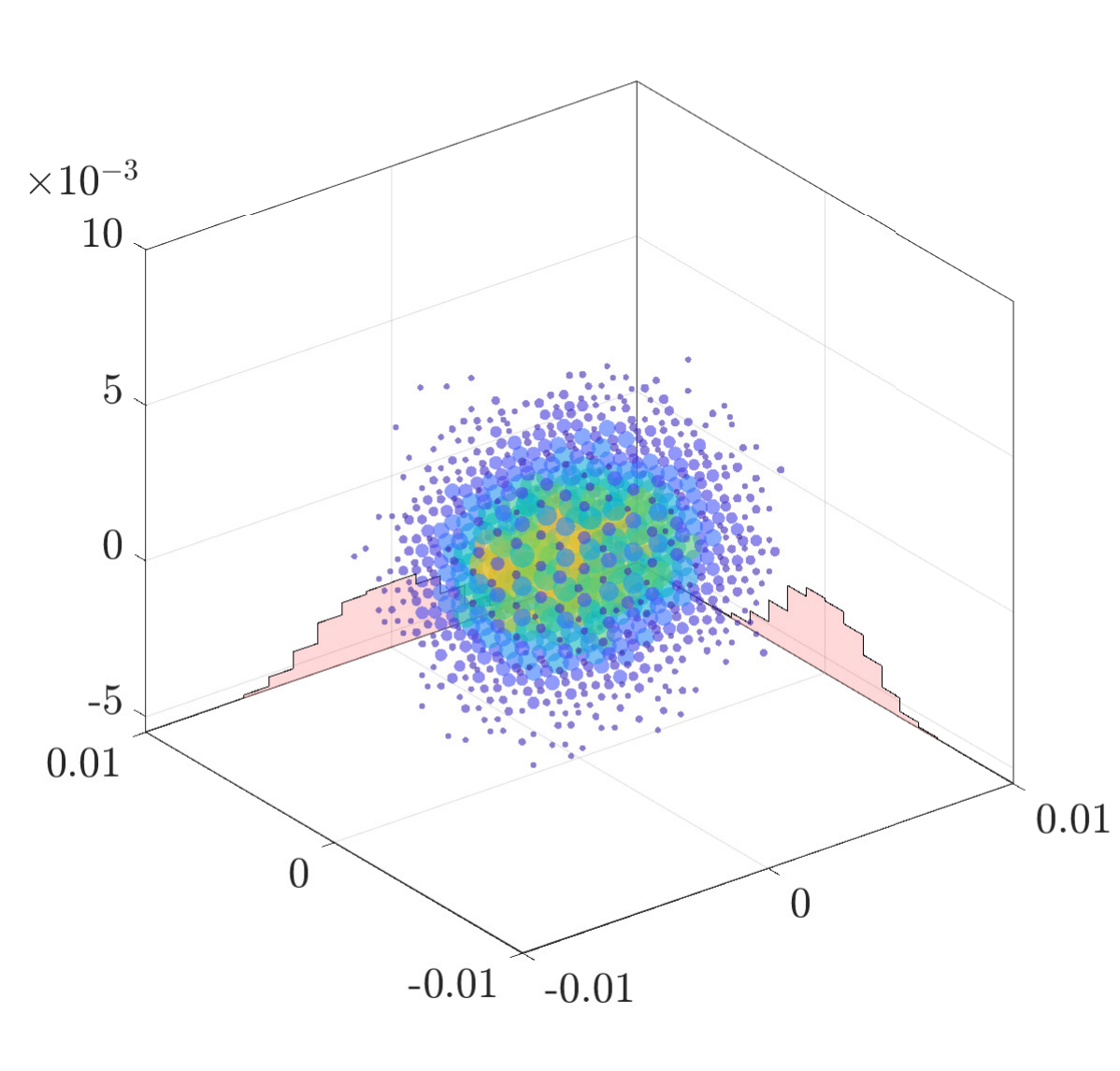}\\
    \includegraphics[width =\textwidth]{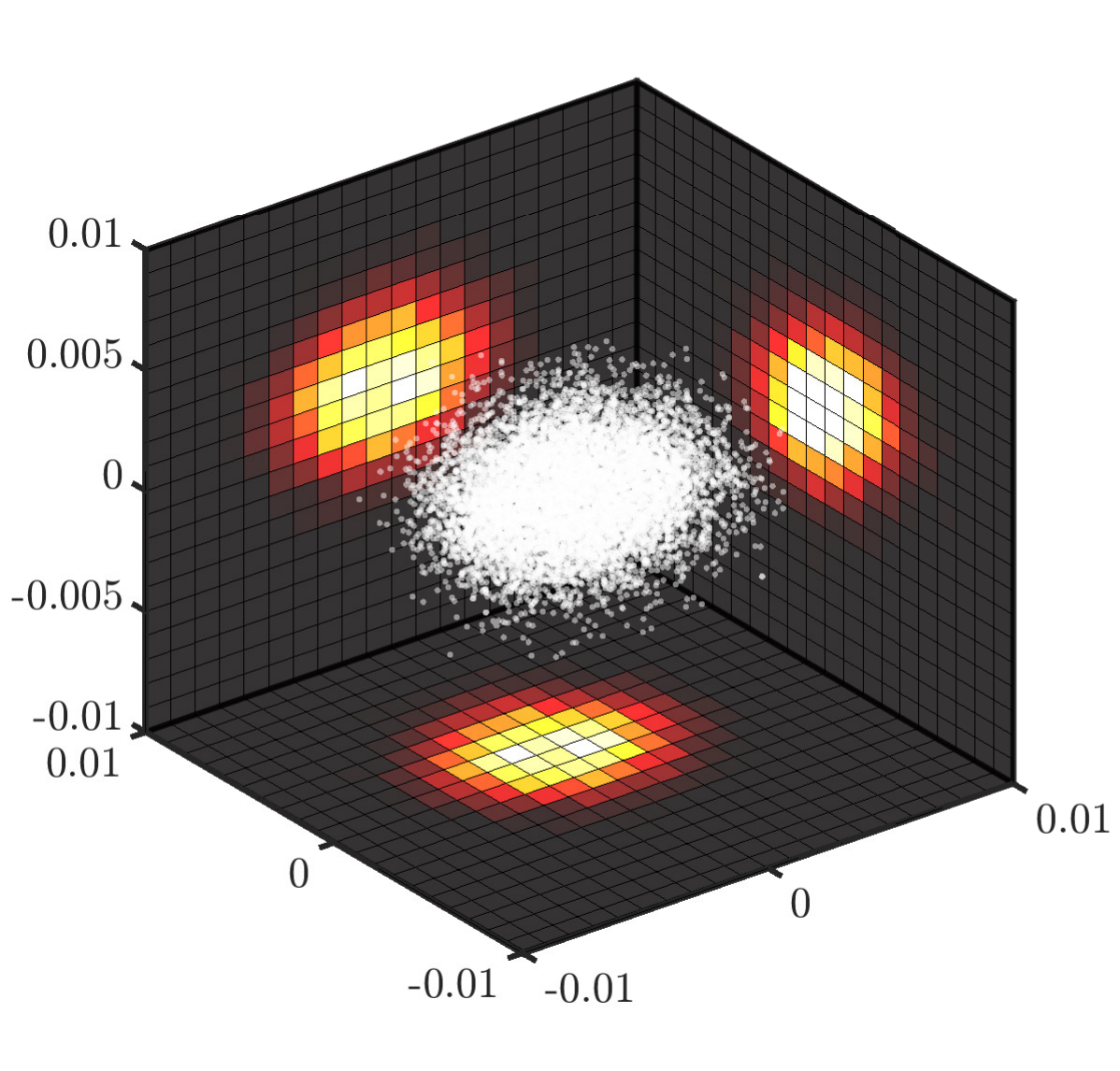}
    \caption{$\bu^{RNN}_\theta$}
    \end{subfigure}\hspace{-0.3cm}
    \begin{subfigure}[b]{0.25\textwidth}
    \includegraphics[width =\textwidth]{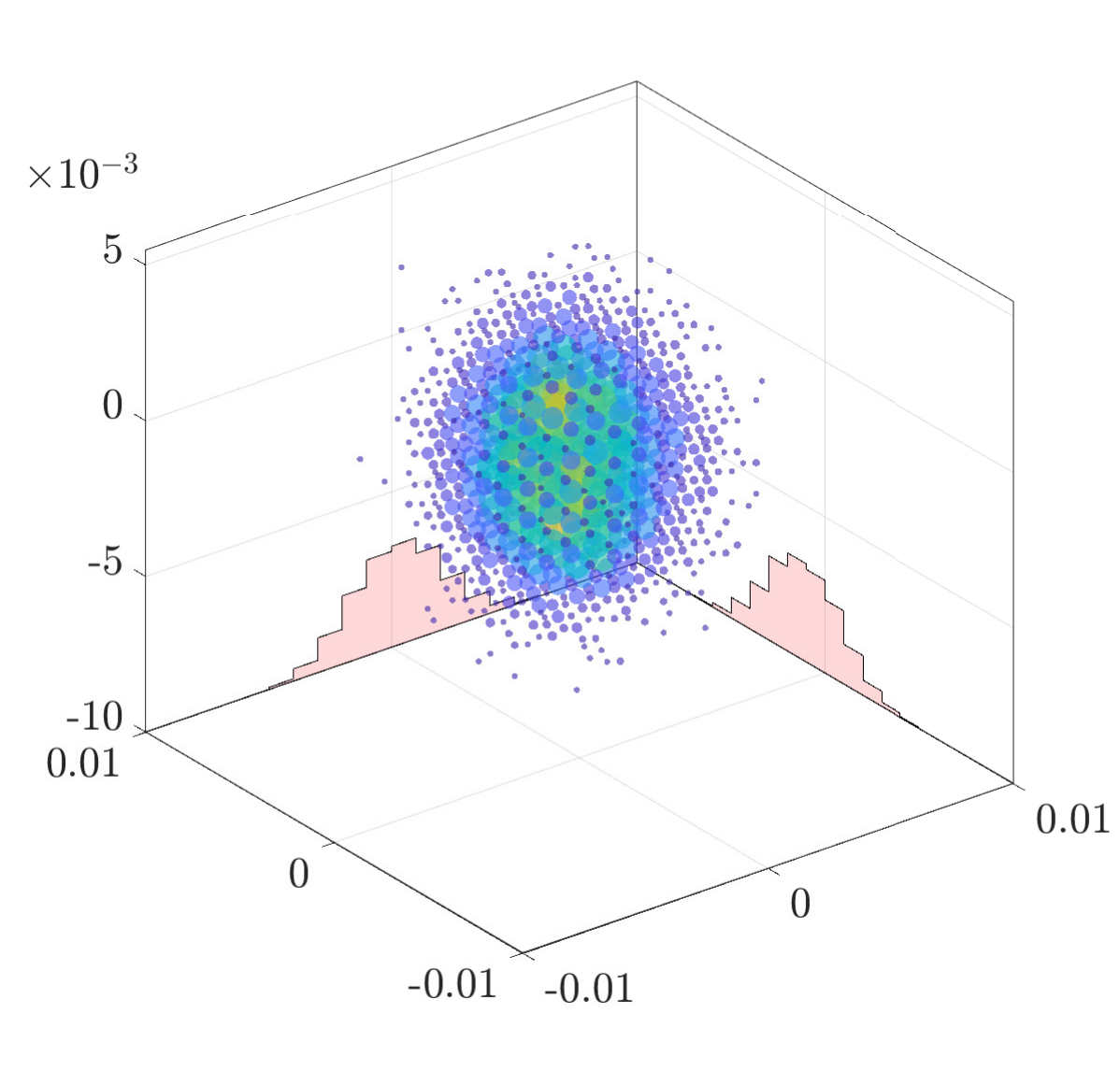}\\
    \includegraphics[width =\textwidth]{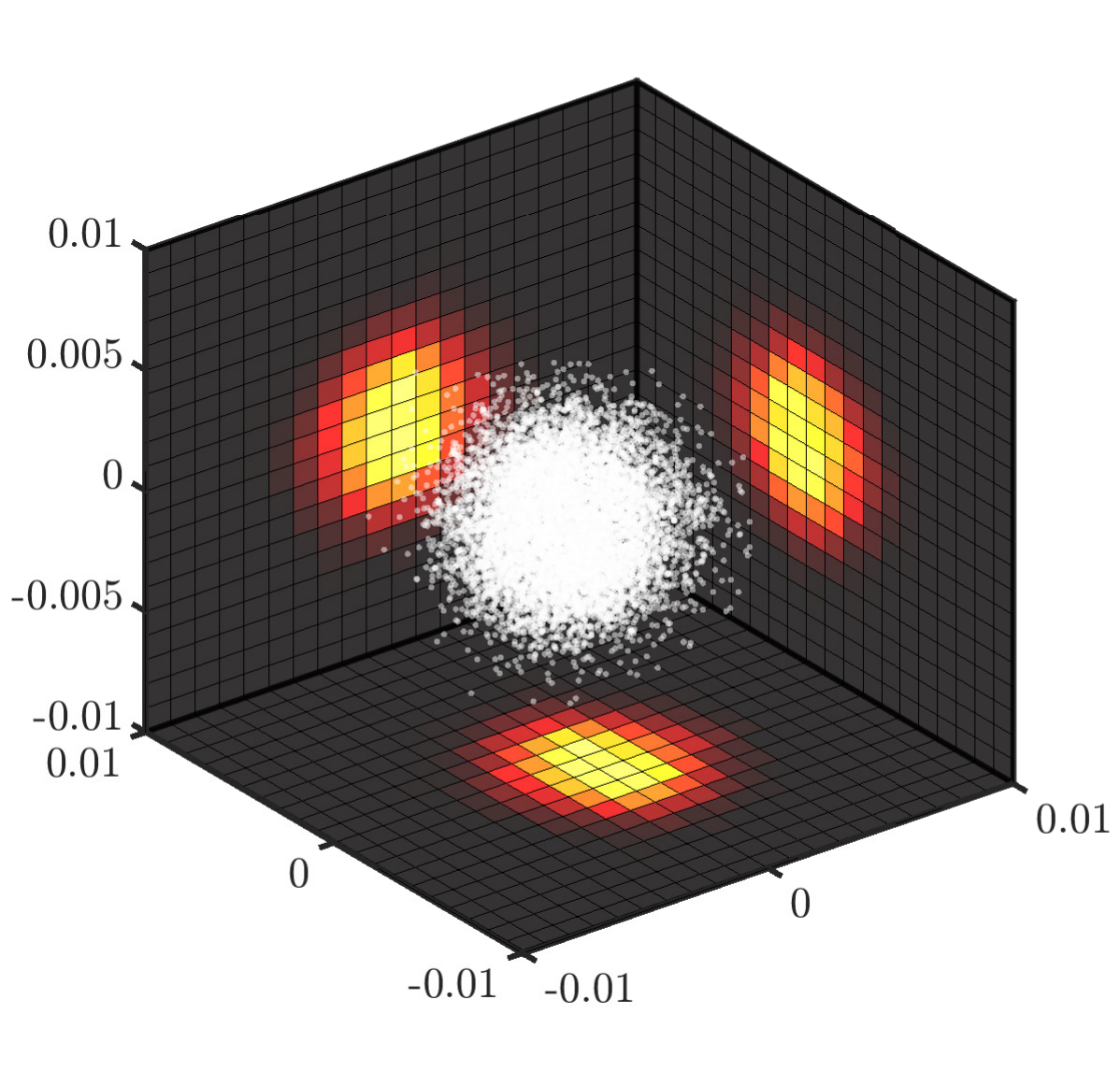}
    \caption{$\bu^{FNN}_\theta$}
    \end{subfigure}
    \caption{Density of the system in the velocity space.  First row: density of the uncontrolled system at initial and final times $t \in \{0s,2s\}$. Second row: controlled system configuration at  time $T = 2s$, with projection of partial densities along the first and second dimension. Third row: configuration of the controlled system at the final time $T = 2s$, with projection of partial densities along all the three dimensions.}
    \label{fig:morse_density}
\end{figure}

\begin{figure}
    \centering
    \begin{subfigure}[b]{0.35\textwidth}
    \includegraphics[width =\textwidth]{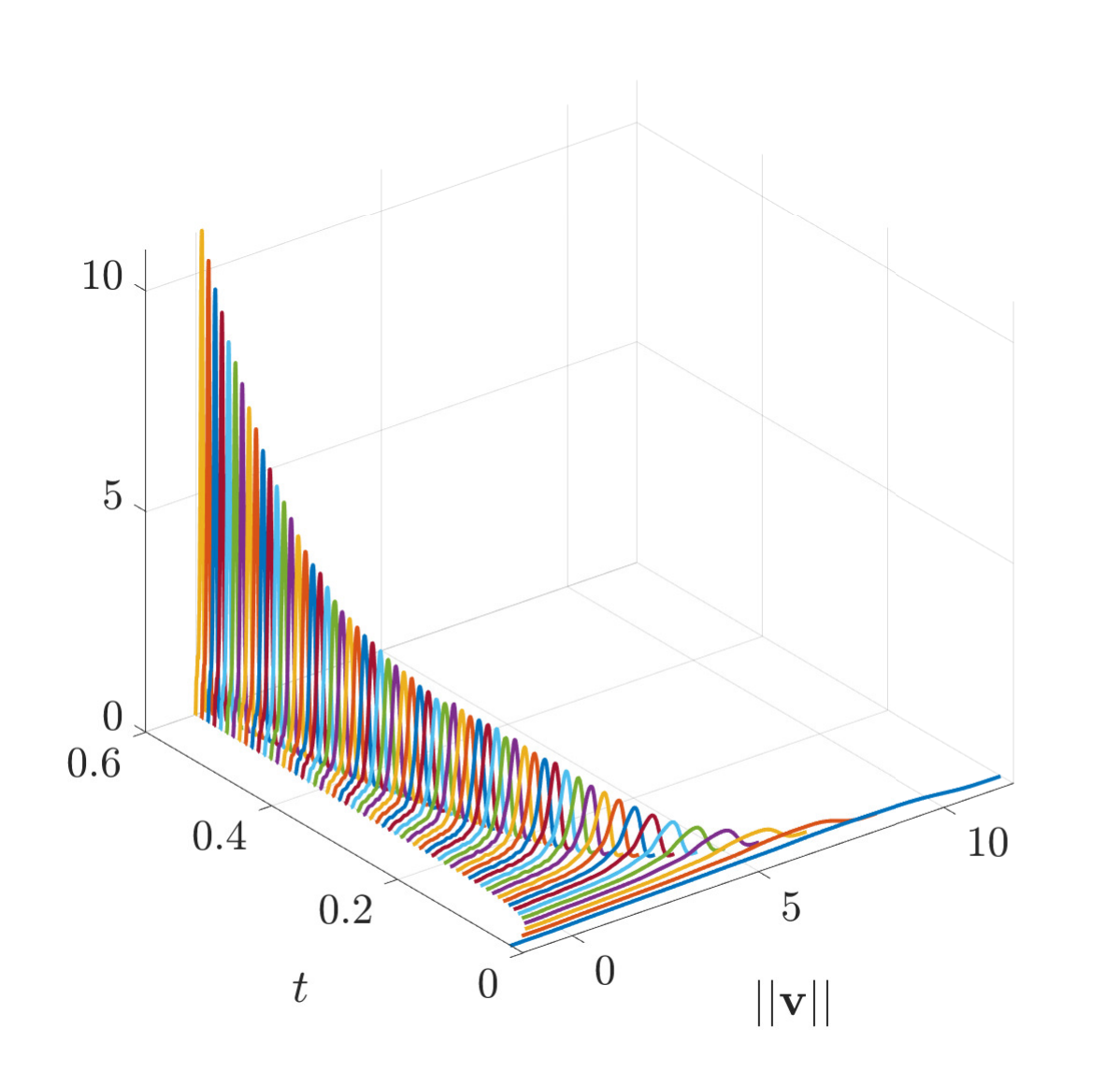}
    \caption{$\bu^{RNN}_\theta$}
    \end{subfigure}
     \begin{subfigure}[b]{0.35\textwidth}
    \includegraphics[width =\textwidth]{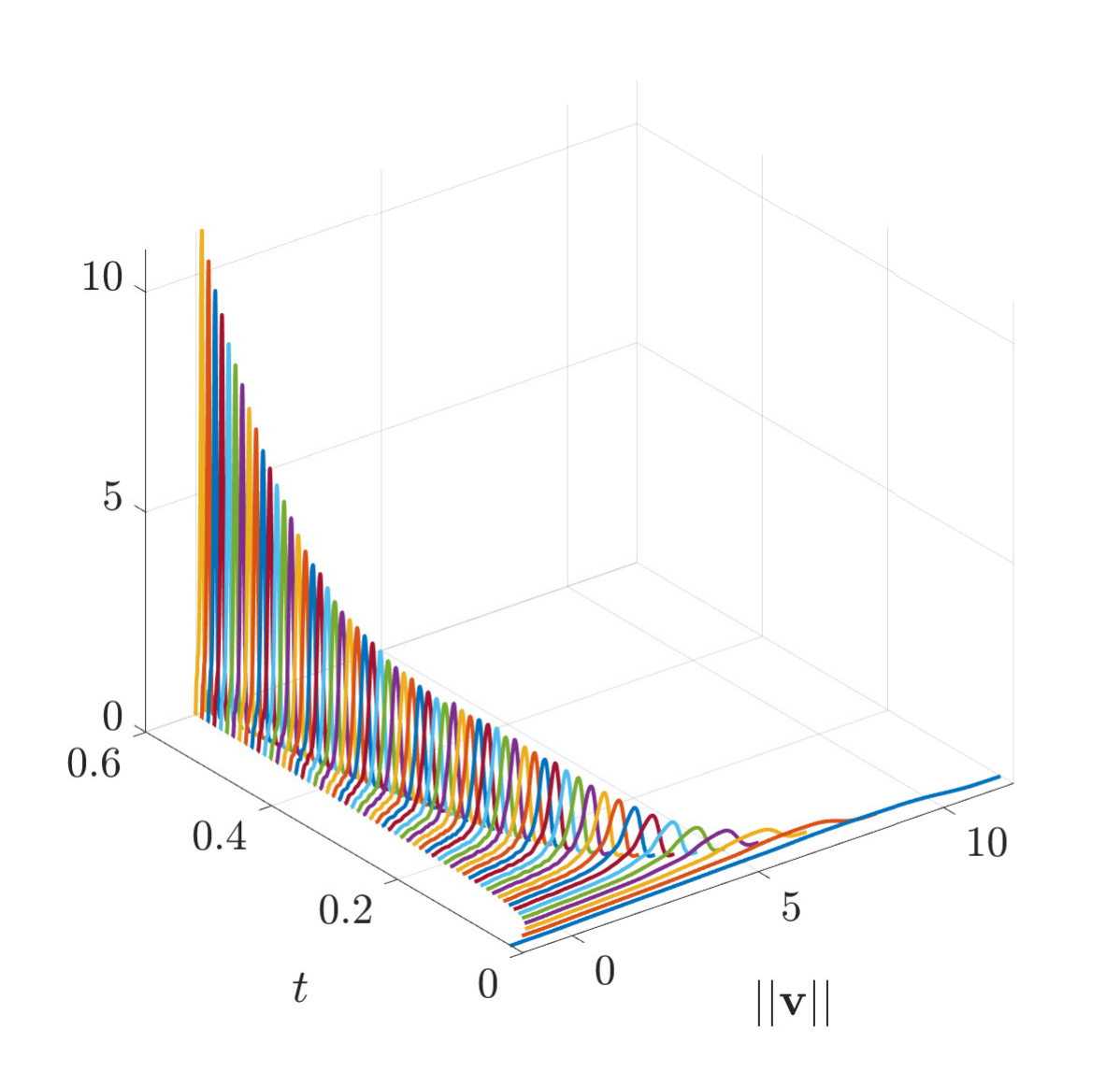}
    \caption{$\bu^{FNN}_\theta$}
    \end{subfigure}\\
     \begin{subfigure}[b]{0.35\textwidth}
    \includegraphics[width =\textwidth]{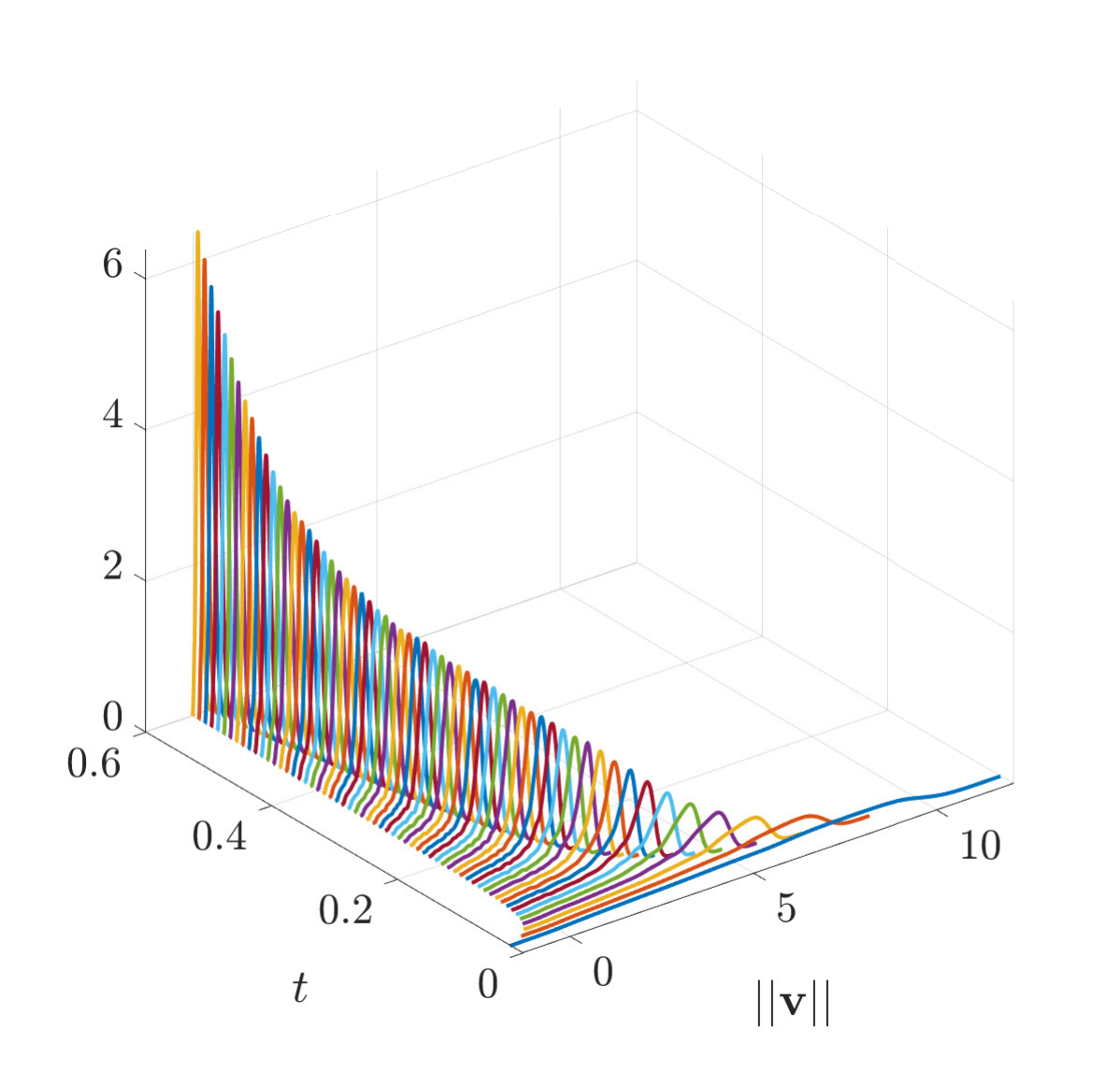}
    \caption{${\bs'}^{RNN}_\theta$}
    \end{subfigure}
     \begin{subfigure}[b]{0.35\textwidth}
    \includegraphics[width =\textwidth]{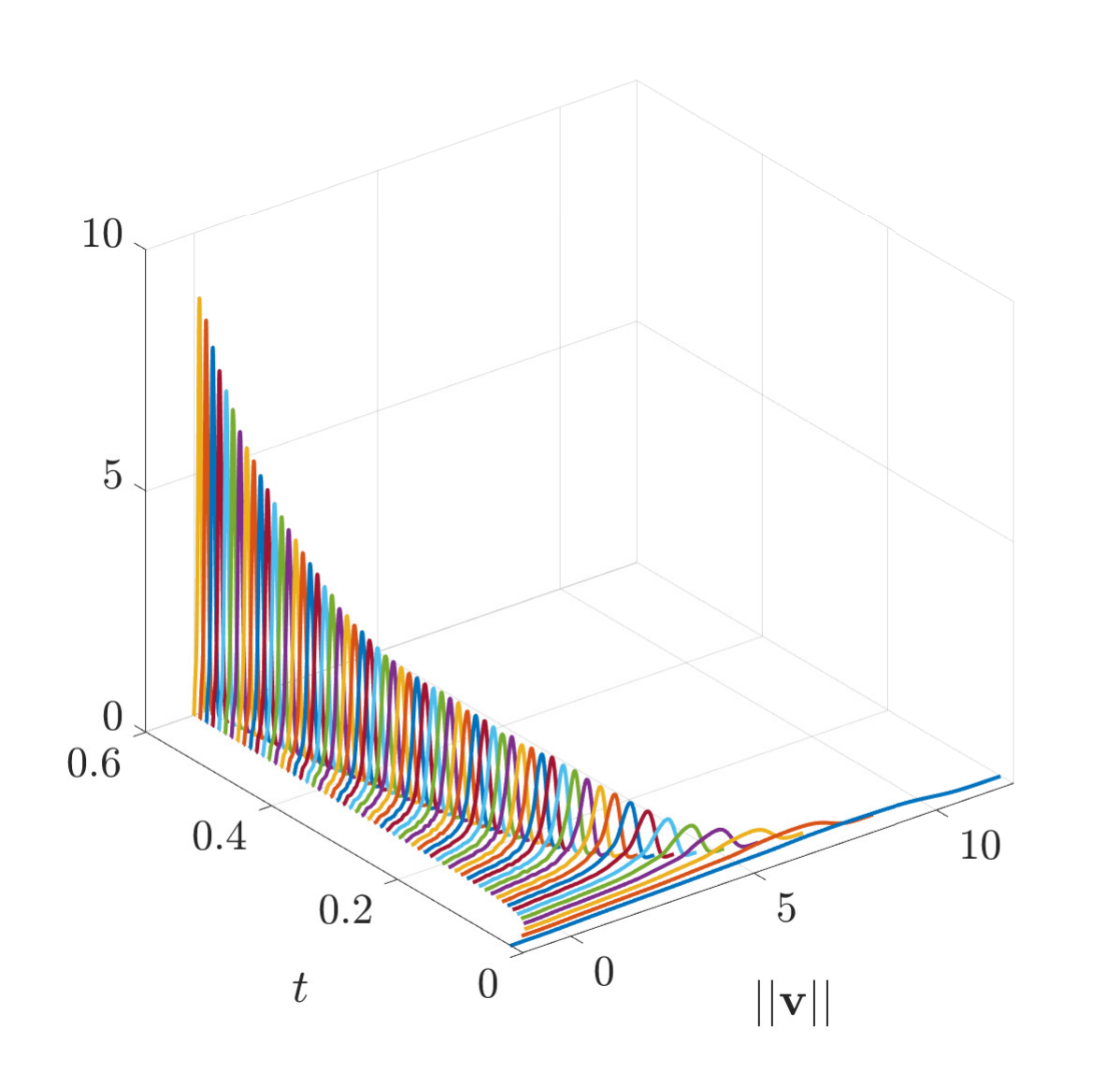}
    \caption{${\bs'}^{FNN}_\theta$}
    \end{subfigure}
    \begin{subfigure}[b]{0.35\textwidth}
    \includegraphics[width =\textwidth]{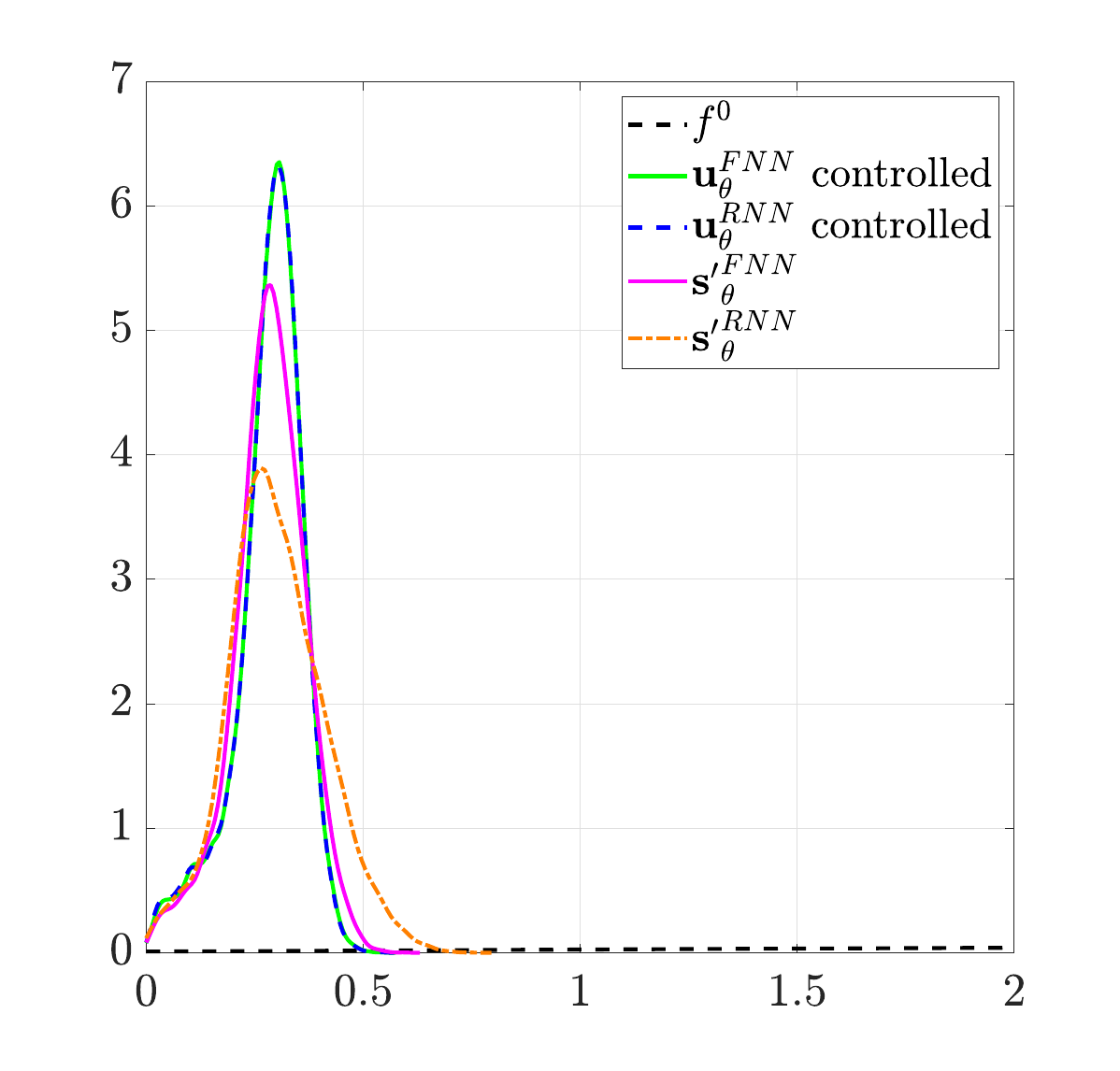}
    \caption{$t_{5} = 0.5s$}
    \end{subfigure}
    \begin{subfigure}[b]{0.35\textwidth}
    \includegraphics[width =\textwidth]{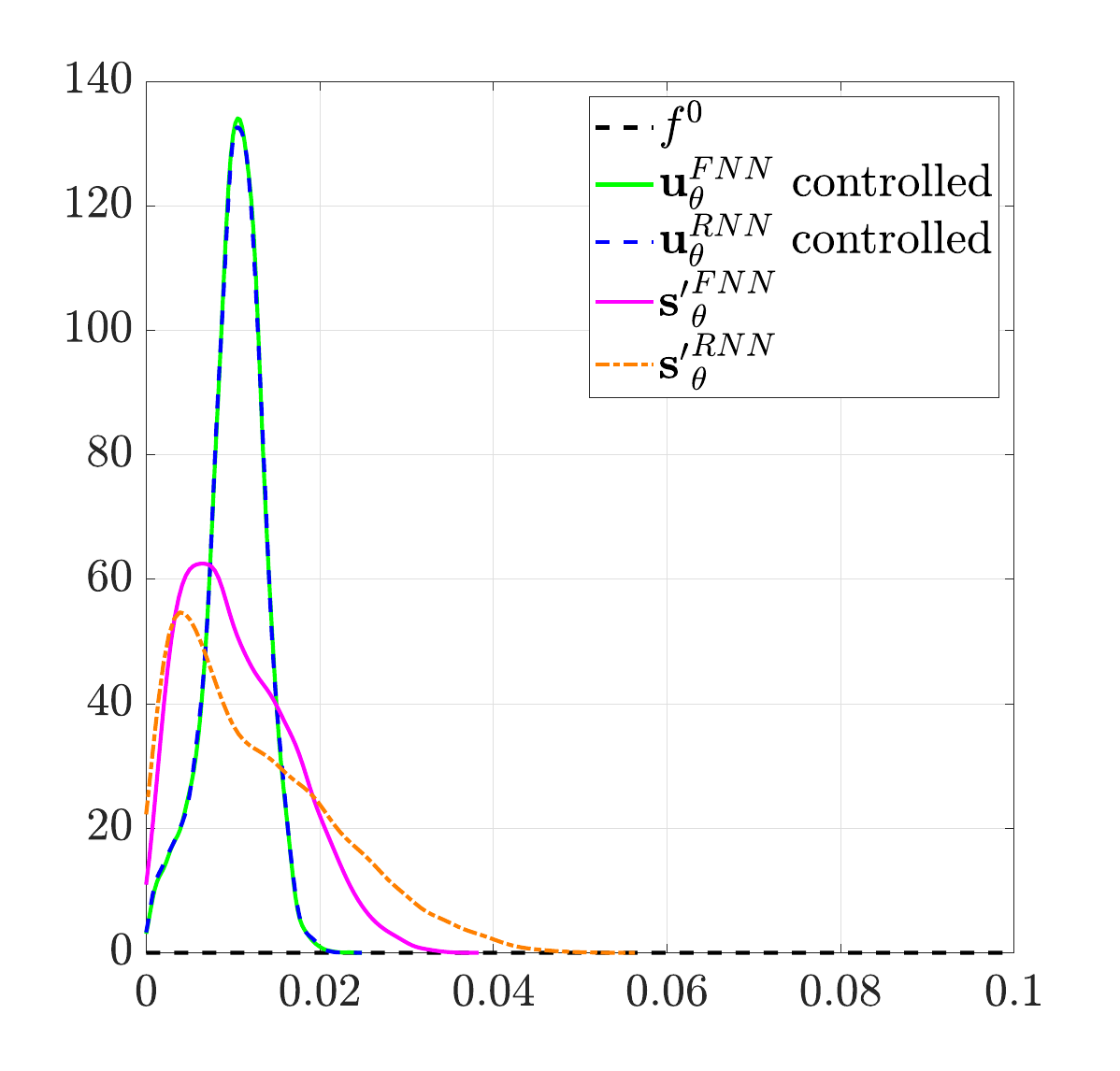}
    \caption{$t_{100} = 1s$}
    \end{subfigure}
    \caption{Density of agents' velocities in norm for test 2. Comparison of the time evolution of the MC pdf obtained considering different approximation approaches for the controlled binary interactions. In the bottom row, we display a comparison of the different approximated densities at two different discrete times.}
    \label{fig:morse_simulation}
\end{figure}

The goodness of fit of these models outside $\mathcal{T}\cup\mathcal{T}_v$ is displayed in table \ref{tab:goodness}, whilst a comparison of the approximately controlled binary system is displayed in Figure \ref{fig:morse_binary}, where the dynamics of a couple of agents with states randomly sampled in $[-4,4]^{12}$ evolve throughout $40$ discrete time intervals of length $\Delta t = 0.02s$. Whilst all the NN models lead to a similar final state, those approximating the feedback design perform closer to the reference DSDRE solution.
As shown in Figure \ref{fig:morse_density}, the controlled evolution converges to concentrated density profiles, exhibiting slight differences in mean and variance among the various models.  This can also be noticed from the \ref{fig:morse_simulation}, which depicts the density of agents along the first two dimensions, together with the associated vector field of the velocities.

\section{Conclusions}

We have developed a novel computational method for mitigating the curse of dimensionality arising in the optimal control of large-scale, high-dimensional, agent-based models. The key ingredients of the proposed approach are: the use of a kinetic model to reduce the complexity associated to the particle ensemble to the sampling of two-agent subproblems, and the synthesis of control laws for the two-agent system by means of deep neural networks, supervised learning, and a discrete-time state-dependent Riccati equation approach.

Our numerical experiments validate the effectiveness of the approach in controlling consensus and attraction-repulsion dynamics in high-dimensional interacting particle systems. The use of neural network approximation models for fast feedback synthesis allows for a speedup of 2 to 3 orders of magnitude compared to solving a binary optimal control problem every time two interacting agents are sampled. 
\color{black}
The proposed framework can be extended to construct a deep neural network that directly predicts the post-interaction states, or to different control problems, provided that a solver for synthetic data generation is available; moreover, the same approach can be applied to other types of dynamics, as long as they can be simulated or approximated through suitable numerical methods
\color{black}

Potential future research directions include the generalization of the proposed approach to the $\mathcal{H}_\infty$ robust control framework, which could further enhance the applicability of the method to real-world scenarios where robustness to uncertainties in the agents' interaction forces is crucial. Additionally, exploring the integration of adaptive time-stepping techniques with the neural network approximation models could lead to faster and more accurate simulations of the controlled agent-based system. Finally, investigating the scalability of the proposed approach to even higher-dimensional problems and larger agent populations could enable applications in consensus-based optimization and mean field limits of neural networks.

\section*{Acknowledgments}
GA is member of Indam GNCS and thanks the support of MUR-PRIN Project 2022 PNRR No. P2022JC95T, ``Data-driven discovery and control of multi-scale interacting artificial agent systems", and of MUR-PRIN Project 2022 No. 2022N9BM3N  ``Efficient numerical schemes and optimal control methods for time-dependent partial differential equations" financed by the European Union - Next Generation EU. This research was supported by the UK Engineering and Physical Sciences Research Council (EPSRC) grant EP/T024429/1
\appendix
\section{Proof of Theorem 2.1}\label{AppendixA}

\paragraph{Step 1}
Let us introduce a test function  {$\varphi \in C_0^2\big(\mathbb{R}^d\times \mathbb{R}^d\big)$}, is the state space hosting the $N_a$ agents of the system. We consider the following weak formulation of the Boltzmann equation \eqref{boltzmann}
\begin{equation}\label{weak_form}
\begin{aligned}
       \dfrac{d}{dt}\langle f,\varphi\rangle + \langle f,v\cdot\nabla_x\varphi\rangle&= \lambda\langle \mathcal{Q}_{\eta,u}(f,f),\varphi\rangle  \\ &= \lambda\iint\limits_{\mathbb{R}^{2d}\times \mathbb{R}^{2xd}} \big(\varphi(x,v')-\varphi(x,v)\big)f(t,x_*,v_*)\,dx_*\, dv_*\, dx\, dv 
\end{aligned}
\end{equation}

\paragraph{Step 2}
By definition of the binary post-interaction dynamics \eqref{post_int_dyn}, we can explicitly write $v'-v$ as 
\begin{equation}
 \begin{aligned}
    v'-v &= \eta\, P(x,x_*)(v_*-v)+\eta\, u_\eta(x,v,x_*,v_*)=: \eta F_\eta(x,v,x_*,v_*)
    \end{aligned}
\end{equation}
where we introduce the function $F_\eta(x,v,x_*,v_*)$ for the controlled dynamics binary dynamics. Moreover, we expand $\varphi(v')$ inside the operator \eqref{weak_form} in Taylor series of  $v'-v$ up to the second order, obtaining 
\begin{equation}
\begin{aligned}\label{taylor}
    \varphi(x,v') - \varphi(x,v) &=
     \eta\, F_\eta(x,v,x_*,v_*)\cdot\nabla_v\,\varphi (x,v)+ R^{\varphi}(\eta^2)
\end{aligned}
\end{equation} where the term $R^{\varphi}(\eta^2)$ represents the reminder of the Taylor expansion as follows
\begin{equation}\label{eq:app_reminder1}
\begin{aligned}
R^{\varphi}(\eta^2) &= \sum_{i,j=1}^d\left(\left(\partial^{(i,j)}_v\varphi(x,v)-\partial^{(i,j)}_v\varphi(x,\tilde v)\right)(v'-v)_i(v'-v)_j\right)\cr
& = \eta^2\sum_{i,j=1}^d\left(\left(\partial^{(i,j)}_v\varphi(x,v)-\partial^{(i,j)}_v\varphi(x,\tilde v)\right)(F_\eta)_i(F_\eta)_j\right)=:\eta^2 \tilde R^\varphi_\eta(x,v,x_*,v_*)
\end{aligned}
\end{equation}
where $\tilde v = (1-\theta)v'+\theta v$ for some $\theta\in[0,1]$ and we use the multi-index notation for the second order partial derivatives of $\varphi$.
\paragraph{Step 3}
Embedding \eqref{taylor} into the interacting operator for the weak form \eqref{weak_form}, and introducing the quasi invariant scaling \eqref{eq:scaling}, i.e.  $\lambda=1/\varepsilon$, $\eta=\varepsilon$ we have 
\begin{equation}
\begin{aligned}\label{eq:app_reminder2}
      \frac{1}{\varepsilon}\langle \mathcal{Q}_{\varepsilon,u}(f,f),\varphi\rangle &= \iint\limits_{\mathbb{R}^{2d}\times \mathbb{R}^{2d}} F_\varepsilon(x,v,x_*,v_*)\cdot\nabla_v\,\varphi(x,v)\, f\,f_*\,dx_*\, dv_*\, dx\, dv + \varepsilon\mathcal{R}_\varepsilon^\varphi,
      \end{aligned}
\end{equation}
with reminder term
\begin{equation}\label{eq:app_reminder3}
\begin{aligned}
      \mathcal{R}_\varepsilon^\varphi =\iint\limits_{\mathbb{R}^{2d}\times \mathbb{R}^{2d}} \tilde R^\varphi_\varepsilon(x,v,x_*,v_*)\, ff_*\,dx_*\, dv_*\, dx\, dv. 
\end{aligned}
\end{equation}
For $\varepsilon\to0$, assuming that the reminder vanishes to zero, and integrating by parts the scaled weak form of the Povzner-Boltzmann model \eqref{weak_form} we have
\[
 \bigg\langle \partial_t f +v\cdot \nabla_x f+\nabla_v \cdot \bigg[f\int\limits_{\mathbb{R}^{2d}} \left(P(x,x_*)(v_*-v)+ u(x,v,x_*,v_*)\right)\,f(t,x_*,v_*)\,dx_*\, dv_* \bigg],\varphi\bigg\rangle =0
\]
for every $\varphi\in C^2(\mathbb{R}^{d}\times \mathbb{R}^{d})$. {  Hence, in strong form, we retrieve consistency of the scaled Povzner-Boltzmann model \eqref{boltzmann} with controlled mean field model \eqref{eq:FP}.}
 Hence, in the next step we show that the reminder is bounded and the previous limit hold true for $\varepsilon\to 0$.
\paragraph{Step 4} Finally we provide estimate on the reminder of the Taylor expansion, to show that such term vanishes for $\varepsilon\to 0$. We observe that
\begin{align*}
\|\mathcal{R}^\varphi_\varepsilon\|\leq \frac{C}{2}\iint_{\mathbb{R}^{2d}\times \mathbb{R}^{2d}} \|F_\varepsilon(x,v,x_*,v_*)\|^2\, f(t,x,v) f(t,x_*,v_*) \, dx \, dv \, dx_*\,dv_*
\end{align*}
where $C$ is the bounding constant of $C^2_0(\mathbb{R}^{d}\times \mathbb{R}^{d})$, and since $F_\varepsilon(x,v,x_*,v_*) \in L^2_{loc}$ we can conclude that the quasi-invariant limit holds true.
\begin{flushright}
    $\square$
\end{flushright}
\begin{rmk}
 {The proof provided here can be readapted seamlessly in the case of the full state binary dynamics \eqref{eq:genbinary} by assuming that 
$F_\eta(s,s_*):=G({s},{s}_*)+ H u({s},{s}_*)\in L^2_{loc}(\R^{2d}\times\R^{2d})$. More precisely, in Step 1, we consider the test function $\Phi\in C^2_0(\R^{2d})$ and the weak form of \eqref{boltzmann_fullstate} as follows
\begin{align}\label{boltzmann_fullstate_weak}
       \dfrac{d}{dt}\langle f,\Phi\rangle &= \lambda\iint\limits_{\R^{2d}\times\R^{2d}} \big(\Phi({s}')-\Phi({s})\big)f(t,{s}_*)f(t,{{s} })\,d{s}_* d{s}, 
\end{align}
where $s' = (x',v')$ is the post-collisional state in \eqref{eq:genbinary}. In Step 2, by  Taylor's expansion around $s=(x,v)$, we retrieve 
\begin{equation}
 \begin{aligned}\label{eq:taylor}
    \Phi({s}')-\Phi({s}) &= \eta\left( G({s},{s}_*)+ H u({s},{s}_*)\right)\cdot \nabla_{{s}}\Phi({s}) + \eta^2 \tilde R^{\Phi}_\eta,
    \end{aligned}
\end{equation}
with reminder term defined now as
\begin{equation}\label{eq:app_reminder4}
\begin{aligned}
\tilde R^\Phi_\eta(s,s_*):=& \sum_{i,j=1}^{2d}\left(\left(\partial^{(i,j)}_v\Phi(s)-\partial^{(i,j)}_v\Phi(\tilde s)\right)(F_\eta(s,s_*))_i(F_\eta(s,s_*))_j\right),
\end{aligned}
\end{equation}
where $\tilde s= (1-\theta) s + \theta s'$, $\theta\in[0,1]$, similarly to \eqref{eq:app_reminder1}.
Then in Step 3, introducing the scaling \eqref{eq:scaling} and using \eqref{eq:taylor} in \eqref{boltzmann_fullstate_weak} we have
\begin{align}
       \dfrac{d}{dt}\langle f,\Phi\rangle &= \iint\limits_{\R^{2d}\times\R^{2d}} \big(\left( G({s},{s}_*)+ H u({s},{s}_*)\right)\cdot \nabla_{{s}}\Phi({s})\big)f(t,{s}_*)f(t,{{s} })\,d{s}_* d{s}+\varepsilon\mathcal{R}^{\Phi}_\varepsilon,
\end{align}
where the reminder $\mathcal{R}^\Phi_\varepsilon$ is defined as in \eqref{eq:app_reminder3}, as integrating $\tilde R_\varepsilon^\Phi(s,s_*)$.
Thus, integrating by parts, we have that for any test function $\Phi$ is
\begin{align}
       \left\langle \partial_t f + \nabla_{{s}}\cdot \left(f\int_{\R^{2d}}\left( G({s},{s}_*)+ H u({s},{s}_*)\right)f({s}_*)\,d{s}_*\right),\Phi\right\rangle = 0.
\end{align}
Finally, the reminder vanishes in the limit $\varepsilon\to 0$ following the same argument of Step 4.}
\end{rmk}
 
\section{Asymptotic constrained symmetric Nanbu algorithm} \label{sec:Appendix_B}
\begin{algorithm}[H] 
\caption{ }\label{alg:pseudoDSMC}
$\left\{(x^{0}_i, v^{0}_i)\right\}_{i=1}^{N_s}\sim f^0\;\;i.i.d.$\Comment*[r]{\small{$N_s$ samples from the initial distribution}}
\For{$h=0,...,N_T-1$, $t_n = n\cdot \Delta t$}{
select $\{(i_k,j_k)\}_{k=1}^{N_s/2}\;$ random pairs of agents without repetitions\\
\For{$k = 1,...,N_s/2$}{
Compute $(v^{n+1}_{i_k}, v^{n+1}_{j_k})$ according to \eqref{binary_euler} \Comment*[r]{\small{interaction}}
}
\For{$i = 1,...,N_s$}{
Compute $x_{i}^{n+1}$ according to \eqref{transport} \Comment*[r]{\small{transport}}
}
}
\end{algorithm}\vspace{0.5cm}

For the generalized binary dynamics \eqref{eq:genbinary} an analogous stochastic simulation technique can be designed, where in this case \texttt{interaction} and \texttt{transport} are updated simultaneously considering the particle states at time $t_n$ as follows $s_{i}^{n}=(x_i^n,v_i^n)^\top$, $s_{j}^{n}:=(x_j^n,v_j^n)^\top$.
 


\end{document}